\documentclass[11pt,reqno]{amsart}

\usepackage{amscd}
\usepackage{amsfonts}
\usepackage{amsmath}
\usepackage{amssymb}
\usepackage{amsthm}
\usepackage{fancyhdr}
\usepackage{latexsym}
\usepackage[colorlinks=true, pdfstartview=FitV, linkcolor=blue,
            citecolor=blue, urlcolor=blue]{hyperref}
\usepackage{enumitem}

\usepackage{mathtools}            
\usepackage{upgreek}

\usepackage{cancel}
\input epsf
\input texdraw
\input txdtools.tex
\input xy
\xyoption{all}

\usepackage{tikz}
\usepackage{tikz-cd}

\usetikzlibrary{patterns}
\usetikzlibrary{calc,matrix, arrows}
\usetikzlibrary{snakes}
\usetikzlibrary{arrows.meta}
\usepackage{subfigure}
\usetikzlibrary{decorations.pathreplacing}
\tikzstyle{level 1}=[level distance=3.5cm, sibling distance=8cm]
\tikzstyle{level 2}=[level distance=3.5cm, sibling distance=3cm]
\tikzstyle{level 3}=[level distance=5cm, sibling distance=1.4cm]

\usepackage{color}

\definecolor{Red}{rgb}{1.00, 0.00, 0.00}
\definecolor{DarkGreen}{rgb}{0.00, 1.00, 0.00}
\definecolor{Blue}{rgb}{0.00, 0.00, 1.00}
\definecolor{Cyan}{rgb}{0.00, 1.00, 1.00}
\definecolor{Magenta}{rgb}{1.00, 0.00, 1.00}
\definecolor{DeepSkyBlue}{rgb}{0.00, 0.75, 1.00}
\definecolor{DarkGreen}{rgb}{0.00, 0.39, 0.00}
\definecolor{dgreen}{RGB}{0,200,100}
\definecolor{ddgreen}{RGB}{0,170,0}
\definecolor{SpringGreen}{rgb}{0.00, 1.00, 0.50}
\definecolor{DarkOrange}{rgb}{1.00, 0.55, 0.00}
\definecolor{OrangeRed}{rgb}{1.00, 0.27, 0.00}
\definecolor{DeepPink}{rgb}{1.00, 0.08, 0.57}
\definecolor{DarkViolet}{rgb}{0.58, 0.00, 0.82}
\definecolor{SaddleBrown}{rgb}{0.54, 0.27, 0.07}
\definecolor{Black}{rgb}{0.00, 0.00, 0.00}
\definecolor{dark-magenta}{rgb}{.5,0,.5}
\definecolor{myblack}{rgb}{0,0,0}
\definecolor{darkgray}{gray}{0.5}
\definecolor{lightgray}{gray}{0.75}
\newcommand{\nn}{\nonumber}
\newcommand{\p}{\partial}

\renewcommand{\b}{\mathcolor{blue}}
\usepackage{scalerel}
\usepackage{stackengine}
\stackMath
\def\widehatgap{-5.5pt}
\def\subdown{-3.3pt}
\newcommand\what[2][]{%
\renewcommand\stackalignment{l}%
\stackon[\widehatgap]{#2}{%
\stretchto{%
    \scalerel*[\widthof{$#2$}]{\kern-.6pt\text{\textasciicircum}\kern-1.1pt}%
    {\rule[-0.8\textheight]{1ex}{\textheight}}%WIDTH-LIMITED BIG WEDGE
}{2ex}% THIS SQUEEZES THE WEDGE TO 0.5ex HEIGHT
_{\smash{\belowbaseline[\subdown]{\scriptscriptstyle#1}}}%
}}

\newcommand\reallywidehat[1]{%
\savestack{\tmpbox}{\stretchto{%
  \scaleto{%
    \scalerel*[\widthof{\ensuremath{#1}}]{\kern-.6pt\bigwedge\kern-.6pt}%
    {\rule[-\textheight/2]{1ex}{\textheight}}%WIDTH-LIMITED BIG WEDGE
  }{\textheight}% 
}{0.6ex}}%
\stackon[1pt]{#1}{\tmpbox}%
}
\newcommand{\rr}{\mathbb{R}}
\renewcommand{\p}{\partial}

\newcommand{\cc}{\mathbb{C}}
\newcommand{\NN}{\mathbb{N}}

\def\medno{\medskip \noindent}

\def\nin{\noindent}

\makeatletter
\def\mathcolor#1#{\@mathcolor{#1}}
\def\@mathcolor#1#2#3{%
  \protect\leavevmode
  \begingroup
    \color#1{#2}#3%
  \endgroup
}
\makeatother

\DeclareMathSizes{12}{12}{7}{5}

\def\refer #1\par{\noindent\hangindent=\parindent\hangafter=1 #1\par}

\theoremstyle{plain}  % default
\newtheorem{theorem}{Theorem}[section]

\newtheorem{lemma}{Lemma}[section]

\usepackage[breakable, theorems, skins]{tcolorbox}
\tcbset{enhanced}

\theoremstyle{definition}

\newenvironment{Proof}[1][\proofname]
{\proof[\textnormal{\textbf{#1.}}]}{\endproof}
\newcommand{\bp}{\begin{Proof}}
\newcommand{\ep}{\end{Proof}}

\numberwithin{figure}{section}
\numberwithin{equation}{section}

\usepackage[framemethod=TikZ]{mdframed}

\makeatletter
\@namedef{subjclassname@2020}{%
  \textup{2020} Mathematics Subject Classification}
\makeatother
\begin{document}
\title{
A higher dispersion  
K\MakeLowercase{d}V 
equation on the half-line}
\author{A. Alexanddrou Himonas  \& Fangchi Yan}
\date{March 8, 2021/Revised April 4, 2022
 \mbox{}$^*$\!\textit{Corresponding author}:
himonas.1@nd.edu}



%\tableofcontents

\begin{abstract}
The initial-boundary value problem (ibvp) for  the $m$-th order 
dispersion Korteweg-de Vries (KdV)  equation on the half-line with rough 
data and solution in restricted Bourgain spaces is studied using  the Fokas Unified Transform Method (UTM). Thus, this work advances the implementation of the Fokas method,
 used earlier  for the KdV on the half-line with smooth data 
 and solution in the classical Hadamard space, consisting
 of function that are continuous in time and Sobolev
 in the spatial variable, to the more general Bourgain spaces 
 framework of dispersive equations with rough data
 on the half-line.
 The spaces needed and the estimates required arise 
 at the linear level and in particular in the estimation
of the  linear pure ibvp,  which has forcing
 and initial data zero but non-zero boundary data.
 Using the iteration map defined by the Fokas solution formula 
 of the forced linear ibvp in combination with the bilinear estimates in 
 modified Bourgain spaces introduced by this map, well-posedness 
 of the nonlinear ibvp is established 
 for rough initial and boundary data belonging
 in Sobolev spaces
of the same optimal regularity
 as in the case of the  initial value problem for this equation on the whole line.
 \end{abstract}

\maketitle

\markboth{
A higher dispersion KdV equation on the half-line
}{
Alex Himonas  \& Fangchi Yan
}

%\tableofcontents
%
%
%%%%%%%%%%%%%%%%%%%%%
%
%			Introduction
%
%%%%%%%%%%%%%%%%%%%%%
%
%

%
\section{Introduction and Results}
In this work, we study the  initial-boundary value problem (ibvp) for the $m$-th order dispersion Korteweg-de Vries  (KdVm)  equation  on the half-line,
that is
\begin{subequations}
\label{KdVm}
\begin{align}
\label{KdVm eqn}
&\p_tu+(-1)^{j+1}\p_x^{m}u +uu_x
=
0,
\quad
0<x<\infty,
\,\,
0<t<T,
\\
\label{KdVm ic}
&u(x,0)
=
u_0(x),
\quad
0<x<\infty,
\\
\label{KdVm bc}
&
u(0,t)
=
g_0(t),
\cdots ,
\p_x^{j-1}u(0,t)
=
g_{j-1}(t),
\quad
0<t<T,
\end{align}
\end{subequations}
where  $m=2j+1$,  $j\in \mathbb{N}$. Note that $m=3$ ($j=1$) gives
the celebrated Korteweg-de Vries  (KdV)  equation
(\cite{kdv1895}, \cite{b1877}),
and $m=5$ ($j=2$) gives a Kawahara equation \cite{kawahara1972}.
The initial data belong in the Sobolev spaces $H^s(0,\infty)$
and the boundary data $g_\ell(t)\in H_t^{(s+j-\ell)/{m}}(0,T)$,
$\ell=0,1,\dots,j-1$, arise in the estimation
of the linear pure ibvp and also reflect the time regularity
of the solution to the linear KdVm initial value problem (ivp).
Here, we show that  if 
$s>-j+\frac 14$,
then  ibvp \eqref{KdVm} is well-posed
with solution in an appropriately modified Bourgain space restricted
to $(0, \infty)\times (0, T)$ for some lifespan $T>0$ depending 
on the size of the data.  
Thus,  the optimality  $s>-j+\frac 14$
of our  ibvp well-posedness result  here is exactly the same 
with the one obtained in \cite{fhy2020} for the Cauchy problem
of this equation on the whole line.
The starting point is to define the iteration map
via the Fokas solution formula of the forced linear ibvp and then 
after deriving appropriate bilinear estimates for the nonlinearity
in the modified Bourgain spaces we show that this map is 
a contraction in a ball. For the KdV equation with smooth data ($s>3/4$) 
well-posedness on the half-line was proved in \cite{fhm2016},
and for KdVm with smooth data ($s>m/4$) well-podsedness 
of the ibvp was proved in \cite{y2020}. Both works are based 
on the Fokas method but the solution spaces used,
which are  subsets of the space of functions that are 
Sobovev in $x$ and continuous in $t$,
are motivated by the work on the KdV Cauchy problem
by Kenig, Ponce and Vega in \cite{kpv1991}.
Our work here advances the implementation of the Fokas method
for solving ibvp with smooth data
(see \cite{fhm2016, fhm2017, hm2015, hm2020, hmy2019, hmy2019b, y2020})
to the Bourgain spaces framework of dispersive 
equations with rough data on the half-line.

Next, we define the spaces needed for stating our main result precisely.
We recall that  for any $s$ and $b$ real numbers, the Bourgain space $X^{s,b}(\rr^2)$ corresponding to the linear part of  KdVm 
is defined by the norm
\begin{equation}
\label{def-Bourgain}
\|u\|_{X^{s,b}}^2
=
\|u\|_{s,b}^2
=
\int_{-\infty}^\infty\int_{-\infty}^\infty
(1+|\xi|)^{2s}
(1+|\tau-\xi^m|)^{2b}
|\widehat{u}(\xi,\tau)|^2
d\xi
d\tau,
\end{equation}
where  $\widehat{u}$ denotes the space-time Fourier transform
\begin{equation}
\label{Fourier-transform}
\widehat{u}(\xi,\tau)
=
\int_{-\infty}^\infty\int_{-\infty}^\infty
e^{-i(\xi x +\tau t)} u(x, t) dx dt.
\end{equation}
Also, we shall need the following modification 
of the Bourgain norm
\begin{align}
\label{Bourgain-like-norms-def}
\|u\|_{X^{s,b,\alpha}(\rr^2)}^2
=&
\|u\|_{s,b,\alpha}^2
\doteq
\|u\|_{X^{s,b}(\rr^2)}^2
+
\Big[
\int_{-\infty}^\infty
\int_{-1}^1
(1+|\tau|)^{2\alpha}
|\widehat{u}(\xi,\tau)|^2
d\xi
d\tau 
\Big]
\\
\simeq&
\int_{-\infty}^\infty
\int_{-\infty}^\infty
\Big[
(1+|\xi|)^{s}
(1+|\tau-\xi^m|)^{b}
+
\chi_{|\xi|< 1}(1+|\tau|)^{\alpha}
\Big]^2
|\widehat{u}(\xi,\tau)|^2
d\xi
d\tau.
\notag
\end{align}
A similar modification was introduced first 
by Bourgain \cite{b1993-kdv}, where
the norm used for the periodic case with $b=1/2$
was modified by the $\alpha$-part in order to prove
well-posedness of the KdV equation  on the line in $H^s$, $s\ge 0$.  
This idea was also utilized later by 
 Colliander, Keel, Staffilani, Takaoka and Tao
in \cite{ckstt2003} for the global well-posedness of KdV on the circle for $s\ge -1/2$.
The Bourgain norm with $b>1/2$ was used 
by Kenig, Ponce and Vega in \cite{kpv1996}
in order to extend the well-posedness of 
KdV in $H^s(\rr)$ to $s>-3/4$.

Finally, we shall need  the  restriction space 
$X_{\rr^+\times(0, T)}^{s,b,\alpha}$,
which is defined as follows
\begin{equation}
\label{sb-restrict}
X_{\rr^+\times(0, T)}^{s,b,\alpha}
\doteq
\{
u;
u(x,t)
=
v(x,t)
\;\;\mbox{on}\;\; \mathbb{R^+}\times (0, T)
\,\,
\text{with}
\,\,
v\in X^{s,b,\alpha}(\rr^2)
\},
\end{equation}
and which is equipped  with the norm
\begin{equation}
\label{sb-restrict-norm}
\| u\|_{X_{\rr^+\times(0, T)}^{s,b,\alpha}}
\doteq
\inf\limits_{v\in X^{s,b,\alpha}}\left\{
\| v\|_{s,b,\alpha};\;v(x,t)
=
u(x,t)
\;\;\mbox{on}\;\; \mathbb{R^+}\times (0, T)
\right\}.
\end{equation}
Now, we are ready to state our first result precisely. It reads as follows. 
\begin{theorem}
[Well-posedness of KdVm on the half-line]
\label{thm-kdvm-half-line}
If $-j+\frac14< s\le j+1$, $s\neq \frac12,\frac32,\dots,j-\frac12$, then there exists $b\in(0,\frac12)$ depending on $s$ such that 
for any   initial data $u_0 \in H^s(0,\infty)$ 
and  boundary data $g_\ell\in H^{\frac 1m(s+j-\ell)}(0,T)$, $\ell=0,1,\dots,j-1$,  there is a lifespan $T_0>0$ such that
the ibvp 
\eqref{KdVm} with 
compatibility condition:
\begin{align}
\label{first-com-con}
\p_x^\ell u_0(0)
=
g_\ell(0),
\quad
\text{for } \ell \text{ such that }
\,\,
\frac 1m (s+j-\ell)
>
\frac12
\,
\text{ or }
\,
s>\ell+\frac12,
\end{align}
admits 
a unique solution $u
\in
X_{\rr^+\times(0, T_0)}^{s,b,\alpha}
$,  satisfying the size estimate
\begin{equation}
\label{Bourgain-like-norms-est}
\|u\|_{X_{\rr^+\times(0, T_0)}^{s,b,\alpha}}
\leq
C
\Big(
\|u_0\|_{H^s(\rr^+)}
+
\sum\limits_{\ell=0}^{j-1}
\|g_\ell\|_{H^{\frac{s+j-\ell}{m}}(0,T)}
\Big),
\end{equation}
for some $\alpha\in (\frac 12, 1)$.
Furthermore, an estimate for the lifespan is given by
\begin{equation}
\label{lifespan-est}
T_{0}
=
c_0
\Big(
1+\|u_0\|_{H^s(\rr^+)}
+\sum\limits_{\ell=0}^{j-1}
\|g_\ell\|_{H^{\frac{s+j-\ell}{m}}(0,T)
}
\Big)^{-4/\beta},
\quad
c_0=c_0(s,b,\alpha),
\end{equation}
for some $\beta>0$  depending on $s$ and $m$
(see   \eqref{beta-choice}).
Finally, the solution depends Lip-continuously on the data $u_0$
and $g_\ell$, $\ell=0,1,\dots,j-1$.
\end{theorem}
As we have mentioned earlier, the optimality $s>-j+\frac 14$
for the data in our well-posedness Theorem \ref{thm-kdvm-half-line} on the half-line is exactly the same with the optimality
of the data in the well-posedness result on the whole line
of KdVm obtained in \cite{fhy2020}.
%

%%%%%%%%%%%%%%%%%%%%%
%
%
%    Fokas method solution formula
%
%
%%%%%%%%%%%%%%%%%%%%%
%
We prove Theorem \ref{thm-kdvm-half-line}  
by showing that the iteration map defined via the
solution formula of the forced linear KdVm ibvp,
which is obtained by the Fokas method, is a contraction  
in the solution space $ X^{s,b,\alpha}$.
Therefore, we begin with the  linear KdVm ibvp
with forcing, that is
\begin{subequations}
\label{LKdVm}
\begin{align}
\label{LKdVm eqn}
&\p_tu+(-1)^{j+1}\p_x^{2j+1}u 
=
f(x,t),
\quad
0<x<\infty,
\,\,
0<t<T,
\\
\label{LKdVm ic}
&u(x,0)
=
u_0(x),
\hskip1.2in
0<x<\infty,
\\
\label{LKdVm bc}
&
u(0,t)
=
g_0(t),
\cdots,
\p_x^{j-1}u(0,t)
=
g_{j-1}(t),
\quad
0<t<T,
\end{align}
\end{subequations}
where $T>0$ is any given time.
Using the Fokas method,  also referred in the literature as the Unified Transform Method (UTM),  we get the following solution
formula to the problem \eqref{LKdVm} 
(see  \cite{y2020} or Section \ref{kdvm-sln-derivation} 
for an outline of the derivation):
\begin{align}
\label{UTM-sln-compact}
u(x,t)
=
S[u_0,g_0,\dots,g_{j-1};f]
\doteq&
\frac{1}{2\pi}\int_{-\infty}^\infty e^{i\xi x+i\xi ^mt}
[\widehat u_0(\xi)+F(\xi,t)]d\xi
\\
+&
\sum\limits_{p=1}^{j}
\sum\limits_{n=1}^{j+1}
C_{p,n}
\int_{\p D_{2p}^+}e^{i\xi x+i\xi ^mt}[\widehat u_0(\alpha_{p,n}\xi)+F(\alpha_{p,n}\xi,t)]d\xi
\nonumber
\\
+&
\sum\limits_{p=1}^j
\sum\limits_{\ell=0}^{j-1}
C'_{p,\ell}
\int_{\p D_{2p}^+}e^{i\xi x+i\xi ^mt}(i\xi )^{2j-\ell}\tilde g_{\ell}(\xi^m,T)d\xi,
\nonumber
\end{align}
where $C_{p,n}$ and $C'_{p,\ell}$ are constants and the rotation numbers $\alpha_{p,n}$ are given by
\begin{align}
\label{a-rot-angles}
\alpha_{p,n}
\doteq
e^{i[m-(2p+1)+2n]\frac{\pi}{m}},
\quad
p=1,2,\dots,j
\quad
n=1,2,\dots, j+1.
\end{align}
Also, 
$\widehat{u}_0(\xi)$ is the Fourier transform of $u_0(x)$ on the half-line, which is defined by the formula
\begin{equation}
\label{FT-halfline}
\widehat{u}_0(\xi)
\doteq
\int_0^\infty e^{-i\xi x}u_0(x)dx,
\quad
\text{Im}(\xi)\leq 0,
\end{equation}
and $F(\xi,t)$ is the following time integral of the half-line Fourier transform of the forcing $f(\cdot, t)$
\begin{equation}
\label{F-time-transform}
F(\xi,t)
\doteq
\int_0^t e^{-i\xi^m\tau}\hat f(\xi,\tau)d\tau
=
\int_0^t e^{-i\xi^m\tau}\int_0^\infty e^{-i\xi x}f(x,\tau)dxd\tau,
\quad
\text{Im}(\xi)\leq 0.
\end{equation}
Furthermore, $\tilde{g}_\ell$ is the temporal Fourier transform 
of $g_\ell$ over the interval $[0, t]$
\begin{equation}
\label{g-time-transform}
\tilde g_\ell(\xi,t)
\doteq
\int_0^te^{-i\xi\tau}g_\ell(\tau) d\tau.
\end{equation}
Finally,  the domains $D_{2p}^+$
in the upper half-plane are as in the two figures below:
\noindent
\begin{figure}[h!]
\begin{minipage}{0.45\linewidth}
\begin{center}
\begin{tikzpicture}[scale=1]
\fill [line width=1pt, opacity=0.2 ,gray] (2.5,{2.5*tan(180/7)})--(0,0)--(1.7,{1.7*tan(360/7)});
\fill [line width=1pt, opacity=0.2 ,gray] (-0.5,{0.5*tan(540/7)})--(0,0)--(0.5,{0.5*tan(540/7)});
\fill [line width=1pt, opacity=0.2 ,gray] (-1.7,{1.7*tan(360/7)})--(0,0)--(-2.5,{2.5*tan(180/7)});
\draw[line width=1pt, black]
(2.5,{2.5*tan(180/7)})--(-2.5,{-2.5*tan(180/7)});
\draw[line width=1pt, black]
(-2.5,{2.5*tan(180/7)})--(2.5,{-2.5*tan(180/7)});
\draw[line width=1pt, black]
(-3,0)--(3,0);
\draw[line width=1pt, black]
(1.7,{1.7*tan(360/7)})--(-1.7,{-1.7*tan(360/7)});
\draw[line width=1pt, black]
(-1.7,{1.7*tan(360/7)})--(1.7,{-1.7*tan(360/7)});
\draw[line width=1pt, black]
(0.5,{0.5*tan(540/7)})--(-0.5,{-0.5*tan(540/7)});
\draw[line width=1pt, black]
(-0.5,{0.5*tan(540/7)})--(0.5,{-0.5*tan(540/7)});
\draw[line width=1pt, black]
(-3,0)--(-1.5,0);
\draw[line width=1pt, black]
(0,0)--(1.5,0);
\draw[line width=1pt, black,->]
({-0.9},{0.9*tan(180/7)})--({-0.8},{0.8*tan(180/7)});
\draw[line width=1pt, black,->]
({-0.3},{0.3*tan(540/7)})--({-0.2},{0.2*tan(540/7)});
\draw[line width=1pt, black,->]
({0.6},{0.6*tan(360/7)})--({0.5},{0.5*tan(360/7)});
\draw[line width=1pt, black,->]
({0.8},{0.8*tan(180/7)})--({0.9},{0.9*tan(180/7)});
\draw[line width=1pt, black,->]
({0.1},{0.1*tan(540/7)})--({0.2},{0.2*tan(540/7)});
\draw[line width=1pt, black,->]
({-0.5},{0.5*tan(360/7)})--({-0.6},{0.6*tan(360/7)});
\draw (1.3,0)
node[xshift=0.2cm,yshift=0.3cm]
{\it\fontsize{11}{11} $\frac{\pi}{m}$}
arc (0:{180/7}:1.3);
\draw [] (0,-1.5) node {\it\fontsize{11}{11} $D^-_{j+1}$};
\draw [] (-1.6,-1.3) node {\it\fontsize{11}{11} $D_{2j}^-$};
\draw [] (1.5,-1.3) node {\it\fontsize{11}{11} $D_2^-$};
\draw [] (-2.2,-0.6) node {\it\fontsize{11}{11} $D_{2j+1}^-$};
\draw [] (2.2,-0.6) node {\it\fontsize{11}{11} $D_1^-$};
\draw [] (0.7,-1.5) node {\it\fontsize{11}{11} $\cdots$};
\draw [] (-0.8,-1.5) node {\it\fontsize{11}{11} $\cdots$};
\draw [] (0,1.5) node {\it\fontsize{11}{11} $D^+_{j+1}$};
\draw [] (-1.6,1.3) node {\it\fontsize{11}{11} $D_{2j}^+$};
\draw [] (1.5,1.3) node {\it\fontsize{11}{11} $D_2^+$};
\draw [] (-2.2,0.6) node {\it\fontsize{11}{11} $D_{2j+1}^+$};
\draw [] (2.2,0.6) node {\it\fontsize{11}{11} $D_1^+$};
\draw [] (0.7,1.5) node {\it\fontsize{11}{11} $\cdots$};
\draw [] (-0.8,1.5) node {\it\fontsize{11}{11} $\cdots$};
\end{tikzpicture}
\caption{Domains for $j$ odd }\label{D-plane}
\end{center}
\end{minipage}
\hskip0.2in
\begin{minipage}{0.45\linewidth}
\begin{center}
\begin{tikzpicture}[scale=1]
\fill [line width=1pt, opacity=0.2 ,gray] (2.5,{2.5*tan(180/7)})--(0,0)--(1.7,{1.7*tan(360/7)});
%\fill [line width=1pt, opacity=0.2 ,gray] (-0.5,{0.5*tan(540/7)})--(0,0)--(0.5,{0.5*tan(540/7)});
\fill [line width=1pt, opacity=0.2 ,gray] (-1.7,{1.7*tan(360/7)})--(0,0)--(-2.5,{2.5*tan(180/7)});
\draw[line width=1pt, black]
(2.5,{2.5*tan(180/7)})--(-2.5,{-2.5*tan(180/7)});
\draw[line width=1pt, black]
(-2.5,{2.5*tan(180/7)})--(2.5,{-2.5*tan(180/7)});
\draw[line width=1pt, black]
(-3,0)--(3,0);
\draw[line width=1pt, black]
(1.7,{1.7*tan(360/7)})--(-1.7,{-1.7*tan(360/7)});
\draw[line width=1pt, black]
(-1.7,{1.7*tan(360/7)})--(1.7,{-1.7*tan(360/7)});
\draw[line width=1pt, black]
(0.5,{0.5*tan(540/7)})--(-0.5,{-0.5*tan(540/7)});
\draw[line width=1pt, black]
(-0.5,{0.5*tan(540/7)})--(0.5,{-0.5*tan(540/7)});
\draw[line width=1pt, black]
(-3,0)--(-1.5,0);
\draw[line width=1pt, black]
(0,0)--(1.5,0);
\draw[line width=1pt, black,->]
({-0.9},{0.9*tan(180/7)})--({-0.8},{0.8*tan(180/7)});
\draw[line width=1pt, black,->]
({-0.3},{0.3*tan(540/7)})--({-0.2},{0.2*tan(540/7)});
\draw[line width=1pt, black,->]
({0.6},{0.6*tan(360/7)})--({0.5},{0.5*tan(360/7)});
\draw[line width=1pt, black,->]
({0.8},{0.8*tan(180/7)})--({0.9},{0.9*tan(180/7)});
\draw[line width=1pt, black,->]
({0.1},{0.1*tan(540/7)})--({0.2},{0.2*tan(540/7)});
\draw[line width=1pt, black,->]
({-0.5},{0.5*tan(360/7)})--({-0.6},{0.6*tan(360/7)});
\draw (1.3,0)
node[xshift=0.2cm,yshift=0.3cm]
{\it\fontsize{11}{11} $\frac{\pi}{m}$}
arc (0:{180/7}:1.3);
\draw [] (0,-1.5) node {\it\fontsize{11}{11} $D^-_{j+1}$};
\draw [] (-1.6,-1.3) node {\it\fontsize{11}{11} $D_{2j}^-$};
\draw [] (1.5,-1.3) node {\it\fontsize{11}{11} $D_2^-$};
\draw [] (-2.2,-0.6) node {\it\fontsize{11}{11} $D_{2j+1}^-$};
\draw [] (2.2,-0.6) node {\it\fontsize{11}{11} $D_1^-$};
\draw [] (0.7,-1.5) node {\it\fontsize{11}{11} $\cdots$};
\draw [] (-0.8,-1.5) node {\it\fontsize{11}{11} $\cdots$};
\draw [] (0,1.5) node {\it\fontsize{11}{11} $D^+_{j+1}$};
\draw [] (-1.6,1.3) node {\it\fontsize{11}{11} $D_{2j}^+$};
\draw [] (1.5,1.3) node {\it\fontsize{11}{11} $D_2^+$};
\draw [] (-2.2,0.6) node {\it\fontsize{11}{11} $D_{2j+1}^+$};
\draw [] (2.2,0.6) node {\it\fontsize{11}{11} $D_1^+$};
\draw [] (0.7,1.5) node {\it\fontsize{11}{11} $\cdots$};
\draw [] (-0.8,1.5) node {\it\fontsize{11}{11} $\cdots$};
\end{tikzpicture}
\caption{Domains for $j$ even }
\label{D-plane}
\end{center}
\end{minipage}
\end{figure}

\noindent
Below we show the Fokas solution formula \eqref{UTM-sln-compact} for the KdV ($j=1$) and its domain $D^+=D^+_2$ 
\begin{align}
\label{KdV-UTM}
u(x,t)
&=
\frac{1}{2\pi}\int_{-\infty}^\infty e^{i\xi x+i\xi^3t}
[\widehat u_0(\xi)+F(\xi,t)]d\xi
\nonumber
\\
&+
\frac{1}{2\pi}\int_{\p D^+} e^{i\xi x+i\xi^3t}\{
e^{i\frac{2\pi}{3}}[\widehat u_0( e^{i\frac{2\pi}{3}} \xi)+F(e^{i\frac{2\pi}{3}} \xi,t)]+e^{i\frac{4\pi}{3}}[\widehat u_0(e^{i\frac{4\pi}{3}} \xi)+F(e^{i\frac{4\pi}{3}} \xi,t)]\}d\xi
\nonumber
\\
&-
\frac{3}{2\pi}\int_{\p D^+} e^{i\xi x+i\xi^3t}\xi^2 \tilde g_0(\xi^3,t)d\xi.
\end{align}
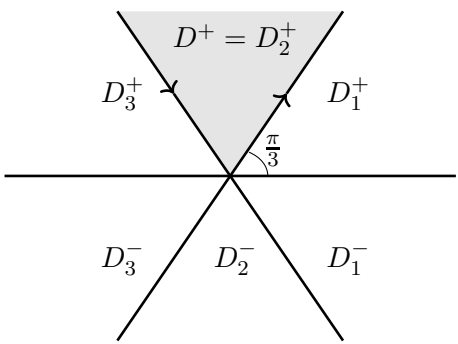
\begin{minipage}{0.5\linewidth}
\begin{center}
\begin{tikzpicture}[yscale=0.73]
\fill [line width=1pt, opacity=0.2 ,gray] (-1.5,3)--(0,0)--(1.5,3);
\draw[line width=1pt, black,->]
(-0.75,1.5)--(0,0)
(-1.5,3)--(-0.75,1.5);
\draw[line width=1pt, black]
(-1.5,0)--(-3,0)
(0,0)--(-1.5,0);
\draw[line width=1pt, black]
(-0.75,-1.5)--(0,0)
(-1.5,-3)--(-0.75,-1.5);
\draw[line width=1pt, black,->]
(0.75,1.5)--(1.5,3)
(0,0)--(0.75,1.5);
\draw[line width=1pt, black]
(1.5,-3)--(0.75,-1.5)
(0,0)--(0.75,-1.5);
\draw[line width=1pt, black]
(0,0)--(1.5,0)
(3,0)--(1.5,0);
\draw (0.5,0)
node[above]
{\it\fontsize{11}{11} $\frac{\pi}{3}$}
arc (0:60:0.5);
\draw [] (0,2.5) node {\it\fontsize{11}{11} $D^+=D_2^+$};
\draw [] (0,-1.5) node {\it\fontsize{11}{11} $D_2^-$};
\draw [] (-1.5,-1.5) node {\it\fontsize{11}{11} $D_3^-$};
\draw [] (1.5,-1.5) node {\it\fontsize{11}{11} $D_1^-$};
\draw [] (-1.5,1.5) node {\it\fontsize{11}{11} $D_3^+$};
\draw [] (1.5,1.5) node {\it\fontsize{11}{11} $D_1^+$};
\end{tikzpicture}
\end{center}
\end{minipage}
\begin{minipage}{0.5\linewidth}
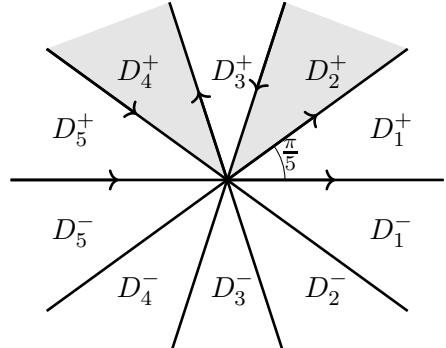
\begin{center}
\begin{tikzpicture}[scale=0.96]
\fill [line width=1pt, opacity=0.2 ,gray] (2.5,{2.5*tan(36)})--(0,0)--(0.8,{0.8*tan(72)});
\fill [line width=1pt, opacity=0.2 ,gray] (-2.5,{2.5*tan(36)})--(0,0)--(-0.8,{0.8*tan(72)});
\draw[line width=1pt, black]
(2.5,{2.5*tan(36)})--(-2.5,{-2.5*tan(36)});
\draw[line width=1pt, black]
(-2.5,{2.5*tan(36)})--(2.5,{-2.5*tan(36)});
\draw[line width=1pt, black]
(-3,0)--(3,0);
\draw[line width=1pt, black]
(0.8,{0.8*tan(72)})--(-0.8,{-0.8*tan(72)});
\draw[line width=1pt, black]
(-0.8,{0.8*tan(72)})--(0.8,{-0.8*tan(72)});
\draw[line width=1pt, black,->]
(-3,0)--(-1.5,0);
\draw[line width=1pt, black,->]
(-2.5,{2.5*tan(36)})--(-1.25,{1.25*tan(36)});
\draw[line width=1pt, black,->]
(0,0)--(1.5,0);
\draw[line width=1pt, black,->]
(0,0)--(1.25,{1.25*tan(36)});
\draw[line width=1pt, black,->]
(0,0)--(-0.4,{0.4*tan(72)});
\draw[line width=1pt, black,->]
(0.8,{0.8*tan(72)})--(0.4,{0.4*tan(72)});
\draw (0.8,0)
node[above]
{\it\fontsize{11}{11} $ \frac{\pi}{5}$}
arc (0:36:0.8);
\draw [] (0,-1.5) node {\it\fontsize{11}{11} $D^-_3$};
\draw [] (-1.3,-1.5) node {\it\fontsize{11}{11} $D_4^-$};
\draw [] (1.3,-1.5) node {\it\fontsize{11}{11} $D_2^-$};
\draw [] (-2.2,-0.7) node {\it\fontsize{11}{11} $D_5^-$};
\draw [] (2.2,-0.7) node {\it\fontsize{11}{11} $D_1^-$};
\draw [] (0,1.5) node {\it\fontsize{11}{11} $D^+_3$};
\draw [] (-1.3,1.5) node {\it\fontsize{11}{11} $D_4^+$};
\draw [] (1.3,1.5) node {\it\fontsize{11}{11} $D_2^+$};
\draw [] (-2.2,0.7) node {\it\fontsize{11}{11} $D_5^+$};
\draw [] (2.2,0.7) node {\it\fontsize{11}{11} $D_1^+$};
\end{tikzpicture}
\end{center}
\end{minipage}
 Also, here we show the solution formula  \eqref{UTM-sln-compact} 
 for KdV5 (Kawahara equation) together with its domains
  $D^+_2$ and   $D^+_4$   
\begin{align}
\label{KdV5-UTM}
u(x,t)
&=
\frac{1}{2\pi}\int_{-\infty}^\infty e^{i\xi x+i\xi^5t}
[\widehat u_0(\xi)+F(\xi,t)]d\xi
\\
&+
\sum\limits_{p=1}^2\sum\limits_{n=1}^3
C_{p,n} \int_{\p D_{2p}^+}e^{i\xi x+i\xi^5t}[\widehat u_0(\alpha_{p,n}\xi)+F(\alpha_{p,n}\xi,t)]d\xi
\nonumber
\\
&+
\sum\limits_{\ell=0}^1\sum\limits_{p=1}^2
C_{p,\ell}'
\int_{\p D_{2p}^+}e^{i\xi x+i\xi^5t}(i\xi)^{4-\ell}\tilde g_\ell(\xi^5,t)d\xi.
\nonumber
\end{align}
For the solution \eqref{UTM-sln-compact} to  ibvp \eqref{LKdVm}, we have the following basic estimate.
\begin{theorem}
[\textcolor{blue}{Forced linear KdVm estimate on the half-line}]
\label{forced-linear-kdvm-thm}
Suppose that $-j-\frac12< s\le j+1$, $s\neq \frac12,\frac32,\dots,j-\frac12$, $0<T<\frac12$. 
Then for some $0<b<\frac12$ and $\alpha>\frac12$ the Fokas  formula 
\eqref{UTM-sln-compact} defines a solution
 to the 
forced linear KdVm ibvp  \eqref{LKdVm} with 
compatibility condition \eqref{first-com-con},
which    satisfies the estimate
\begin{align}
\label{forced-linear-kdvm-est}
&\|S[u_0,g_0,\dots,g_{j-1};f]\|_{X^{s,b,\alpha}(\rr^+\times(0,T))}
\\
\leq&
\begin{cases}
c_{s,b,\alpha}\Big[
\|u_0\|_{H_x^s(0,\infty)}
+
\sum\limits_{\ell=0}^{j-1}\|g_\ell\|_{H_t^\frac{s+j-\ell}{m}(0,T)}
+
\|f\|_{X^{s,-b,\alpha-1}(\rr^+\times(0,T))}
\Big],
\quad
-1
\le
s
\le
\frac12,
\\
c_{s,b,\alpha}\Big[
\|u_0\|_{H_x^s(0,\infty)}
\hskip-0.05in
+
\hskip-0.05in
\sum\limits_{\ell=0}^{j-1}\|g_\ell\|_{H_t^\frac{s+j-\ell}{m}(0,T)}
\hskip-0.1in
+
\|f\|_{X^{s,-b,\alpha-1}(\rr^+\times(0,T))}
\hskip-0.05in
+
\|f\|_{Y^{s,-b}(\rr^+\times(0,T))}
\Big],
\,\,
s
\not\in
[-1,\frac12],
\end{cases}
\nn
\end{align}
where $Y^{s,b}$  is a  ``temporal" Bourgain space defined  by the norm
\begin{equation}
\label{Ysb-def}
\|
u
\|_{Y^{s,b}}^2
\doteq
\int_{-\infty}^\infty
\int_{-\infty}^\infty
(1+|\tau|)^{\frac{2s}{m}}
(1+|\tau-\xi^m|)^{2b}
|\widehat{u}(\xi,\tau)|^2
d\xi d\tau.
\end{equation}
\end{theorem}
\vskip-0.1in
We will prove our well-posedness result by showing that 
the iteration map defined by the Fokas solution formula
is a contraction on a ball of the space 
$X^{s, b,\alpha}(\rr^+\times \rr)$.
The key ingredient is the basic linear estimate 
\eqref{forced-linear-kdvm-est} where the forcing  $f$
is replaced by the KdVm nonlinearity   $\p_x(u^2)$,
which is quadratic.
To apply this linear estimate we will extend 
each one of its two factors  from $\rr^+\times \rr$ to $\rr^2$ appropriately, and 
 to show that the iteration map is a contraction  we  shall need the following bilinear estimate.
\begin{theorem} [Bilinear estimates]
\label{kdv-bilinear-estimate-thm}
If  $s>-j+\frac{1}{4}$, then there is $0<b<\frac12$ and 
 $\frac12<\alpha \le 1-b$  such that
the following estimate holds
\begin{equation}
\label{bilinear-est}
\|\p_x(f\cdot g)  \|_{X^{s,-b,\alpha-1}(\rr^2)}
\le
c_{s,b, \alpha}
\| f\|_{X^{s,b',\alpha'}(\rr^2)} \| g \|_{X^{s,b',\alpha'}(\rr^2)},
\quad
f, g \in X^{s,b',\alpha'}(\rr^2),
\end{equation}
where   $b'$  and  $\alpha'$ are chosen as follows
\begin{equation}
\label{b-b'-alpha-cond}
\frac12-\beta
\le
b'
\le
b
<
\frac12
<
\alpha'
\le
\alpha
\le
\frac12+\beta,
\end{equation}
and $\beta$ is given by 
\begin{equation}
\label{beta-choice}
\beta
=
\begin{cases}
\min\{\dfrac{1}{12m},\frac{m-s}{3m}\}
\quad
s
\ge
0,
\\
\min\Big\{
\dfrac{j-\frac34}{32m},
\dfrac{s+1/2}{2m}
\Big\},
\quad
 -\frac12< s<0,
\\
\dfrac{1}{32m}\big[s-(-j+\frac14)\big],
\quad
-j+\frac14<s\le-\frac12.
\end{cases}
\end{equation}
\end{theorem}
 We note that this is a more general result than  is needed here and  is of  interest in its own right.  Furthermore, it can be shown that for any $s<-j+\frac{1}{4}$ and $b<\frac12$ the estimates \eqref{bilinear-est} fails. 
 A similar estimate in $X^{s,b}$ with $b>1/2$ for the KdVm
 equation on the line was proved in \cite{fhy2020}.
 For the KdV bilinear estimates in the Bourgain space  $X^{s,\frac 12, \alpha}$ were first proved in  \cite{b1993-kdv}, and  in $X^{s,b}$ with $b>1/2$ were proved in \cite{kpv1996}.

Also, we shall need the following  
 bilinear estimate in ``temporal" Bourgain space,
 which is used when we estimate 
 the Sobolev norm of the solution of the forced ivp
 in the time variable if $s>1/2$ or $s<-1$ 
 (see \eqref{forced-ivp-te}).
 For the KdV equation, estimate \eqref{bi-est-Y} was proved by Holmer
in \cite{h2006}.
 %
 %
% \textcolor{red}{This is
% instead  of  the alternative $L^2([0, T]; H_x^s(\mathbb R))$,
%  which leads to the ``crazy" $\Lambda$-norm.}
%

%%%%%%%%%%%%%%%%%
%
%
% Temporal bilinear estimate
%
%%%%%%%%%%%%%%%%%%
%
\begin{theorem}
[Bilinear estimate in $Y^{s,b}$]
\label{bi-est-Y-thm}
If  $b$, $b'$ and $\alpha'$ satisfies \eqref{b-b'-alpha-cond}, then we have 
\begin{equation}
\label{bi-est-Y}
\|
\p_x(fg)
\|_{Y^{s,-b}(\rr^2)}
\le
c_{s,b,\alpha}
\|
f
\|_{X^{s,b',\alpha'}(\rr^2)}
\|
g
\|_{X^{s,b',\alpha'}(\rr^2)},
\quad
-j+1/4<s<m,
\end{equation}
\begin{equation}
\label{bi-est-Y-1}
\|
\p_x(fg)
\|_{Y^{s,-b}}
\le
\|
\p_x(fg)
\|_{X^{s,-b}}
+
c_{s,b}
\|
f
\|_{X^{s,b'}}
\|
g
\|_{X^{s,b'}},
\quad
-j+1/4<s<m.
\end{equation}
\end{theorem}
Besides the method we follow here, there are two other approaches to the study of initial-boundary value problems for KdV type equations.
In the first approach, which has been initiated by Colliander and Kenig
in \cite{ck2002} and by Holmer in \cite{h2006} the forced linear ibvp is written as a superposition of ivps on the line. Then,  the modern harmonic analysis techniques developed earlier for proving well-posedness of the nonlinear equation in Bourgain spaces are utilized. 
In the second approach, which  for KdV has been initiated 
by Bona, Sun and Zhang in \cite{bsz2002, bsz2003, bsz2008},
the forced linear ibvp is solved via a Laplace transform in the time variable and then by deriving appropriate estimates
well-posedness of the ibvp is established.

%{\bf Kenig approach.}
%Instead of solving the reduced pure initial boundary value problem (ibvp) of the KdV equation, they solve a forced initial value problem (ivp) of KdV with some appropriate forcing term. 
%By using the Riemann-Liouville fractional integration operator,  they get a solution formula,  which satisfies the boundary data in the  reduced pure ibvp, to the forced initial value problem. Hence, they 
%get an extension, which is from the half-line   to the whole line, of the solution to the reduced pure ibvp. And such an extension depends on the selection of the forcing terms.
%In fact, in \cite{ck2002} the authors select the forcing term to be $\delta(x)g(t)$, where $g(t)$ depends on the boundary data.   And in \cite{h2006}, the forcing term was chosen as $\frac{x_{-}^\lambda}{\Gamma(\lambda)}g(t)$, where $x_{-}\doteq
%\begin{cases}
%-x,
%\quad
%x
%<
%0,
%\\
%0,
%\quad
%x
%\ge 
%0,
%\end{cases}$, $\lambda$ is a parameter, which depends on $s$ and $g(t)$ relies on the boundary data.

As we have mentioned before  KdVm  includes the KdV equation ($m=3$), which is integrable, and the Kawahara equation ($m=5$), which is not integrable. The literature 
about the KdV is very extensive. It begins with 
Scott Russell's observation of the ``great wave of translation" 
\cite{jsr1845} and continues with the derivation of the KdV model by  Boussinesq in 1877 \cite{b1877} and  Korteweg and de Vries  in 1895 \cite{kdv1895}.
Then, in 1965, Zabusky and Kruskal \cite{zk1965} observed numerically  that soliton solutions of KdV interact almost linearly by preserving their shape and speed after a collision. Soon after the KdV initial  value problem  on the line with data 
   of sufficient smoothness and decay was solved  by
    Gardner, Greene, Kruskal and  Miura
    \cite{ggkm1967}
   via the  inverse scattering transform (IST) which is based 
   on its integrability, that is 
   its  Lax pair formulation \cite{lax1968}.
   The KdV ivp  in Sobolev spaces  $H^s$ using methods 
 from partial differential equations has been studied extensively
 by many authors. For well-posedness  results on the line
 using Bourgain spaces we refer the reader to 
 Bourgain \cite{b1993-kdv} when $s\ge 0$,
  Kenig, Ponce and Vega \cite{kpv1996} when  $s> -3/4$,
  Guo \cite{g2009} when $s=-3/4$,
  and to Colliander, Keel, Staffilani, Takaoka, Tao
  \cite{ckstt2003} for its global well-posedness when $s> -3/4$.
  For additional results we refer the reader to
  \cite{Bona1975, st1976, cs1988,  cks1992, b1993-nls,
  kato1983, kpv1993,kpv2001,kpv1989,sjoberg1970, BSaut2003, lpbook,ckt2003,ghhp2013}
  and the references therein.

Concerning the Fokas method for solving ibvp, whose initial motivation came from integrable equations
(in particular the  KdV and the cubic NLS equations),
we refer the reader
\cite{dtv2014, f1997,  f2002,fp2005, fi2004, fis2005, 
fl2012, lenells2013, fpbook2014,  fs2012, oy2019} 
and the references therein.
Also, for a detailed introduction to this method 
we refer to the book \cite{f2008}.
 For additional work on solving ibvp we refer the reader to
 \cite{et2016, h2005, fa1988,bw1983,ton1977} 
 and the references therein.

 Our work  is structured as follows.
In Section 2, we study the reduced pure ibvp,
 which is the homogeneous linear problem with initial data zero but non-zero boundary data and derive a key estimate
 for its solution in modified Bourgain spaces.
 In Section 3, guided by the reduced pure ibvp,
  we decompose the forced linear ibvp into four simpler 
  linear sub-problems and derive appropriate estimates
  for their solutions. Then, combining these estimates
  we prove Theorem \ref{forced-linear-kdvm-thm},
  which provides the basic estimate for the Fokas solution
  formula and the iteration map of the nonlinear problem.
   In Section 4, we prove the needed bilinear estimate for the nonlinearity
    in  modified Bourgain spaces $X^{s,b, \alpha}$.
   Then, in section 5, we prove the temporal bilinear estimate
   in  $Y^{s,b}$ spaces.
   In Section 6, we prove our KdVm  well-posedness result, Theorem \ref{thm-kdvm-half-line}, using  estimates \eqref{forced-linear-kdvm-est} for the forced linear ibvp and the bilinear estimates \eqref{bilinear-est}.
   Finally, in Section 7 we provide a brief outline of the Fokas solution formula
   for the forced linear ibvp on the half-line.

%
%%%%%%%%%%%%%%%%%%%%%%%%%%%%%
%
%
%
%    Proof of  estimate for  Reduced ibvp part 1
%
%
%
%
%
%%%%%%%%%%%%%%%%%%%%%%%%%%%%
%
%
%
\section{Reduced pure ibvp}
\label{linear-kdv-reduced}
We begin with the most basic linear KdVm ibvp on the half-line.
This is the homogeneous ibvp with zero initial data  and 
nonzero boundary data. Furthermore, we assume that 
the boundary data $h_\ell$ are test functions of time which
are  {\bf compactly supported} in the interval $[0, 2]$.
This problem, which we call the {\bf reduced pure ibvp},
reads as follows:
\begin{subequations}
\label{LKdVm-reduced}
\begin{align}
\label{LKdVm-reduced eqn}
&\p_tv+(-1)^{j+1}\p_x^{2j+1}v=0,
\quad
0<x<\infty,
\,\,
0<t<2,
\\
\label{LKdVm-reduced ic}
&v(x,0)
=
0,
\\
\label{LKdVm-reduced bc}
&v(0,t)
=
h_0(t),
\cdots,
\p_x^{j-1}v(0,t)
=
h_{j-1}(t),
\quad
0<t<2.
\end{align}
\end{subequations}
In this situation we have that
\begin{equation}
\label{h-l-tilde-recall}
\tilde h_\ell(\xi,2)
=
\int_0^2e^{-i\xi\tau}h_\ell(\tau)d\tau
=
\int_\rr e^{-i\xi\tau}h_\ell(\tau)d\tau
=
\widehat h_\ell(\xi),
\end{equation}
 and the Fokas solution formula
 of our reduced pure ibvp \eqref{LKdVm-reduced} takes the simple form
\begin{align}
\label{pure-kdvm-sln}
v(x,t)
=&
S[0,h_0,\cdots,h_{j-1};0]
\\
=&
\sum\limits_{p=1}^j
\sum\limits_{\ell=0}^{j-1}
C'_{p,\ell}\int_{\p D_{2p}^+}e^{i\xi x+i\xi ^mt}(i\xi )^{2j-\ell}\tilde h_{\ell}(\xi^m,2)d\xi
=
\sum\limits_{p=1}^j
\sum\limits_{\ell=0}^{j-1}
C'_{p,\ell}\,  v_{p\ell}
,
\nonumber
\end{align}
where 
\begin{equation}
\label{def-v-p-l}
v_{p\ell}(x,t)
\doteq
\int_{\p D_{2p}^+}e^{i\xi x+i\xi ^mt}(i\xi )^{2j-\ell}
\,
\widehat h_{\ell}(\xi^m)d\xi
\quad
x\in\rr^+,
\quad
t\in[0,2].
\end{equation}
In the next  result we estimate this solution 
in the Hadamard and the Bourgain spaces.
Our objective  is to obtain the optimal 
bounds in temporal Sobolev spaces.
In fact, 
%it comes as no surprise that 
%
in both cases the bounds that arise naturaly
are the ones suggested by the time regularity
of  the solution  to the linear homogeneous
Cauchy problem with data  in $H^s(\rr)$
(for KdV, see \cite{kpv1991-reg}, \cite{h2006}, \cite{fhm2016})
More precisely, we have the following result.
\begin{theorem}
[\b{Estimates for pure ibvp  on the half-line}]
\label{reduced-pure-ibvp-thm}
For boundary data   test functions that are compactly supported 
in the interval $[0,2]$,
the solution for the reduced pure ibvp \eqref{LKdVm-reduced} satisfies the  following
Hadamard space estimate
\begin{align}
\label{reduced-pure-ibvp-H-es}
\sup_{t\in \rr} \|
S[0,h_0,\cdots,h_{j-1};0]
\|_{H_x^s(0,\infty)} 
\leq
c_{s} \sum\limits_{\ell=0}^{j-1}\|h_\ell\|_{H_t^{\frac{s+j-\ell}{m}}(\rr)}
,
\quad
s\ge 0.
\end{align}
In addition, 
for $b\in[0,\frac12)$  and 
$
\frac12<\alpha
\le
\frac{1}{2}
+
\frac 1m (s+j+\frac12)
$ 
it satisfies the Bourgain spaces
 estimate
\begin{align}
\label{reduced-pure-ibvp-B-es}
\|
S[0,h_0,\cdots,h_{j-1};0]
\|_{X^{s, b,\alpha}(\rr^+\times (0,2))} 
\leq
c_{s,b,\alpha}
\sum\limits_{\ell=0}^{j-1}\|h_\ell\|_{H_t^{\frac{s+j-\ell}{m}}(\rr)},
\quad
s>-j-\frac12.
\end{align}
\end{theorem}
The proof of the  Hadamard space estimate 
\eqref{reduced-pure-ibvp-H-es}
can be found in \cite{y2020} for KdVm and  in \cite{fhm2016}
for KdV. The restriction $s\ge 0$ comes from the 
use of the physical space description of the Sobolev norm.
Here, we focus on the Bourgain spaces estimate 
\eqref{reduced-pure-ibvp-B-es}
which is new and useful.
It is the basic ingredient in the proof of well-posedness 
of KdVm on the half-line.

\medno
{\bf Proof of Theorem \ref{reduced-pure-ibvp-thm}.}
Here we present the proof of the estimate \eqref{reduced-pure-ibvp-B-es}.
Using the parametrization
$[0,\infty)\ni \xi\to \gamma\xi$
 for the right side of 
the domain $D_{2p}^+$,
and the parametrization
$[0,\infty)\ni \xi\to \gamma'\xi$ for the left side of 
$D_{2p}^+$,
we obtain the following decomposition 
$
v_{p\ell}(x,t)
=
V_{1}(x,t)+V_{2}(x,t),
$
where
\begin{align}
\label{def-v1}
V_{1}(x,t)
=
\int_0^\infty e^{i\gamma \xi x-i\xi ^mt}(\gamma \xi )^{2j-\ell}
\,
\widehat h_\ell(-\xi ^m)d\xi 
\simeq
\int_0^\infty e^{-i\xi ^mt}e^{i\gamma_R\xi x}e^{-\gamma_I \xi x} \xi ^{2j-\ell}
\,
\widehat h_\ell(-\xi ^m)d\xi,
\end{align}
\begin{align}
\label{def-v2}
V_{2}(x,t)
=
\int_0^\infty e^{i\gamma' \xi x+i\xi^mt}(\gamma' \xi )^{2j-\ell}
\,
\widehat h_\ell(\xi^m)d\xi 
\simeq
\int_0^\infty e^{i\xi^mt}e^{i\gamma'_R\xi x}e^{-\gamma'_I \xi x} \xi^{2j-\ell}
\,
\widehat h_\ell(\xi^m)d\xi,
\quad
\end{align}
and
\begin{align*}
\gamma
=
e^{i(2p-1)\cdot\frac{\pi}{m}}
=
\cos\left(\frac{2p-1}{m}\pi\right)+i\sin\left(\frac{2p-1}{m}\pi\right),
\,\,\,
\gamma'
=
e^{i2p\cdot\frac{\pi}{m}}
=
\cos\left(\frac{2p}{m}\pi\right)+i\sin\left(\frac{2p}{m}\pi\right).
\end{align*}
Note, that  the imaginary parts of $\gamma$ and $\gamma'$
are positive, which is crucial for our estimates.

Here we only estimate $V_{1}$. The estimation  of $V_{2}$ is similar. 
Also, we split $V_1$ as the sum of two functions, 
one for $\xi$ near $0$ and the other away from $0$, that is
$
V_1
=
v_0
+
v_1,
$
where
\begin{equation}
\label{def-v0}
v_0(x,t)
\doteq
\int_0^1 e^{-i\xi ^mt}e^{i\gamma_R\xi x}e^{-\gamma_I \xi x} \xi ^{2j-\ell} \, \widehat h_\ell(-\xi ^m)d\xi,
\quad
x\in\rr^+,
\quad
t\in[0,2],
\end{equation}
and
\begin{equation}
\label{def-v1}
v_1(x,t)
\doteq
\int_1^\infty e^{-i\xi ^mt}e^{i\gamma_R\xi x}
e^{-\gamma_I \xi x} 
\xi ^{2j-\ell}
\,
\widehat h_\ell(-\xi ^m)d\xi,
\quad
x\in\rr^+,
\quad
t\in[0,2].
\end{equation}
The estimate of the Bourgain norm for $v_0$ 
follows from the boundedness of the Laplace transform in $L^2$ and we will do it later. 
Next, we  estimate  the Bourgain norm for $v_1$.
Using the identity
\begin{align*}
\xi^j
[
e^{i\gamma_R\xi x}
e^{-\gamma_I\xi x}
]
=
\frac{\p_x^j
[
e^{i\gamma_R\xi x}
e^{-\gamma_I\xi x}
]}{
(i\gamma)^j
}
\end{align*}
and the fact that  $e^{-\gamma_I \xi x}$ is exponentially 
decaying in $\xi$ for $x>0$ we can
take the $\p_x^j$-derivative outside the integral sign 
in \eqref{def-v1}
to rewrite $v_1(x,t)$ as follows
\begin{equation}
\label{def-v1-modified}
v_1(x,t)
=
\frac{1}{(i\gamma)^j}
\p_x^j
\int_1^\infty e^{-i\xi^mt}e^{i\gamma_R\xi x}e^{-\gamma_I \xi x} \xi^{j-\ell}
\,
\widehat h_\ell(-\xi ^m)d\xi,
\quad
x\in\rr^+,
\quad
t\in[0,2].
\end{equation}
Next, we extend $v_1$  from $\rr^+\times[0,2]$ to $\rr\times\rr$ 
by using  the 
the one-sided  cutoff function $\rho(x)$, which satisfies $0\le\rho(x)\le1, x\in\rr$, and is as follows
\vskip.01in
\noindent
\begin{minipage}{0.6\linewidth}
\begin{center}
\begin{tikzpicture}[yscale=0.5, xscale=0.8]
%%%%%%%%%%%%%%%%%%%
%
%Variables defined
%
%%%%%%%%%%%%%%%%%%%%
\newcommand\X{0};
\newcommand\Y{0};
\newcommand\FX{11};
\newcommand\FY{11};
\newcommand\R{0.6};
\newcommand*{\TickSize}{2pt};
%%%%%%%%%%%%%%%%%%
%
%End
%
%%%%%%%%%%%%%%%%%%%%%
\draw[black,line width=1pt,-{Latex[black,length=2mm,width=2mm]}]
(-5,0)
--
(5,0)
node[above]
{\fontsize{\FX}{\FY}\bf \textcolor{black}{$x$}};

\draw[black,line width=1pt,-{Latex[black,length=2mm,width=2mm]}]
(0,0)
--
(0,3)
node[right]
{\fontsize{\FX}{\FY}\bf \textcolor{black}{$\rho$}};

\draw[line width=1pt, yscale=2,domain=-1.5:-4.3,smooth,variable=\x,red]  plot ({\x},{0});

\draw[line width=1pt, yscale=2,domain=0:4.3,smooth,variable=\x,red]  plot ({\x},{1});

\draw[smooth,line width=1pt, red]
(0,2)
to[out=-170,in=10]
(-1.5,0)
;

\draw[red]
(2,2.5)
node[]
{\fontsize{\FX}{\FY}$\rho(x)$}

(0,0)
node[yshift=-0.2cm]
{\fontsize{\FX}{\FY}$0$}

(-1.5,0)
node[yshift=-0.2cm]
{\fontsize{\FX}{\FY}$-1$};

\end{tikzpicture}
\end{center}
\end{minipage}
\hskip-0.88in
\begin{minipage}{0.5\linewidth}
\begin{equation}
\label{rho-def}
\rho(x)
=
\begin{cases}
1,
&\quad
x\ge0,
\\
0,
&\quad
x\le -1.
\end{cases}
\end{equation}
\end{minipage}

\nin
Using it, we extend  $v_{1}$ via the formula below (keeping the same notation $v_1$ for it)
\begin{align}
\label{v-p-l-ext-1}
v_{1}(x,t)
\doteq&
\frac{1}{(i\gamma)^j}
\p_x^j
\int_1^\infty 
e^{-i\xi ^mt}e^{i\gamma_R\xi x}e^{-\gamma_I \xi x} \rho(\gamma_I\xi x)
\xi ^{j-\ell}
\,
\widehat h_\ell(-\xi ^m)d\xi,
\quad
x\in\rr,
\quad
t\in\rr,
\end{align}
where $\gamma_I$ is the imaginary part of 
 $\gamma$  and  $\gamma_R$ is its real part.
 Also, we could have localized $v_1$ in $t$ further by multiplying it
 by  the standard cutoff function $\psi(t)$ 
in $C^{\infty}_0(-1, 1)$ such that 
$0\le \psi \le 1$ and $\psi(t)=1$ for $|t|\le 1/2$.
 Then using the estimate
 \begin{align}
\label{mult-by-cutoff}
\|\psi(t)v_1\|_{X^{s, b,\alpha}(\rr^2)} 
\leq
c_{\psi} \|v_1\|_{X^{s, b,\alpha}(\rr^2)},
\end{align}
we are reduced in estimating $\|\cdot\|_{X^{s, b,\alpha}(\rr^2)}$,
which we do next.  Notice that the quantities 
under the integral defining $v_1$ make sense for all $t$
since $t$ appears in oscillatory terms.

\noindent
Extension \eqref{v-p-l-ext-1} is good since
$ \rho(\gamma_I\xi x)=1$ for $x>0$. Also,
$e^{-\gamma_I \xi x} \rho(\gamma_I\xi x)$ is bounded for all $x$ and $t$ since $e^{-\gamma_I \xi x} \le e$ and $\rho\le 1$, that is
$
|
e^{-\gamma_I \xi x} \rho(\gamma_I\xi x)
|
\leq
e^{1}
\cdot
1,
$
$
x\in\rr,
\,\,
t\in\rr.
$
Making the change of variables $\tau=-\xi^m$ and defining
\begin{equation}
\label{def-eta}
\eta(x)
\doteq
e^{i\frac{\gamma_R}{\gamma_I}x}
e^{-x} \rho(x)
,
\end{equation}
we write $v_1$ in the form
\begin{equation}
\label{eqn-V1-1}
v_1(x,t)
\simeq
\p_x^j
\int_{-\infty}^{-1} e^{i\tau t}
\eta(-\gamma_I\tau^{1/m}x)
\tau^{-(\ell+j)/m}
\,
\widehat h_\ell(\tau)
d\tau,
\quad
x\in\rr,
\quad
t\in\rr, 
\end{equation}
and prove the following result  for it.

\begin{lemma}
[Bourgain  space estimate for reduced ibvp]
\label{reduced-pibvp-thm1}
For any $\varepsilon>0$, if $s\geq -j-\frac12-\varepsilon$, 
and $b\ge 0$,
then  the function $v_{1}(x,t)$, which is part of the solution $v=S\big[0, h_0,\dots,h_{j-1}; 0\big]$ to pure ibvp \eqref{LKdVm-reduced} and  defined by \eqref{v-p-l-ext-1}
satisfies the  space estimate
\begin{align}
\label{reduced-ibvp-Bourgain-est}
\|v_1\|_{X^{s, b}} 
\leq
c_{s,b} \|h_\ell\|_{H_t^{\frac{s+mb-\frac12-\ell+\varepsilon}{m}}(\mathbb R)}.
\end{align}
Moreover, if we chose  $\varepsilon$ such that 
$
\varepsilon
\le
m(\frac12-b),
$
which is possible if $b<1/2$, then
we have the estimate (needed in our 
well-posedness theorem)
\begin{align}
\label{reduced-ibvp-Bourgain-est-1}
\|v_1\|_{X^{s, b}} 
\le
c_{s,b}
\|h_\ell\|_{H_t^{\frac{s+j-\ell}{m}}(\mathbb R)}.
\end{align}
\end{lemma}
\nin
{\bf Proof of Lemma \ref{reduced-pibvp-thm1}.} 
In order to estimate the $\|v_1\|_{s,b}$, we need to calculate the Fourier transform of $v_1(x,t)$.
Using the inverse Fourier transform, we get
\begin{align}
\label{v1-t-FT}
\widehat{v}_1^t(x,\tau)
\simeq
\begin{cases}
\p_x^j\eta(-\gamma_I\tau^{1/m}x)
\tau^{-(\ell+j)/m}
\,
\widehat h_\ell(\tau),
&\quad
\tau<-1,
\\
0,
&\quad
\tau
\ge-1.
\end{cases}
\end{align}
In addition, taking the Fourier transform with respect to $x$, we get
\begin{align}
\label{v1-xt-FT}
\widehat{v}_1(\xi,\tau)
\simeq
\begin{cases}
\xi^j
F(\xi,\tau)
\tau^{-(\ell+j)/m}
\,
\widehat h_\ell(\tau),
&\quad
\tau<-1,
\\
0,
&\quad
\tau
\ge-1,
\end{cases}
\end{align}
where $F(\xi,\tau)$ is given by
$
F(\xi,\tau)
\doteq
\int_{x\in\rr}
e^{-i\xi x}\eta(-\gamma_I\tau^\frac1m x)
dx.
$
Also, using the fact that $\eta$ is a Schwarz function and making a change of variables, we get  the following result.
\begin{lemma}
\label{F-bound-lem}
For any 
$n\ge 0$,  $\tau<-1$ and $\xi\in\rr$, we have
\begin{align}
\label{F-bound-ine}
|F(\xi,\tau)|
\leq
c_{\rho,\gamma,n}
\cdot
\frac{1}{|\tau|^{1/m}}
\cdot
\left(
\frac{|\tau|^{1/m}}{|\xi|+|\tau|^{1/m}}
\right)^n,
\end{align}
where $c_{\rho,\gamma,n}$ is a constant depending on 
$\gamma$, $n$ and $\rho$, which is described in \eqref{rho-def}.
\end{lemma}
\nin
Furthermore,  using \eqref{v1-xt-FT} we get
\begin{align}
\label{X-norm-v1}
\|v_1\|_{X^{s,b}}^2
\simeq
\int_{-\infty}^{-1}
\left[
\int_{-\infty}^{\infty}
(1+|\xi|)^{2s}
(1+|\tau-\xi^m|)^{2b}
|
\xi^j
F(\xi,\tau)
|^2
d\xi
\right]
\cdot
\left|
\tau^{-(\ell+j)/m}
\widehat h_\ell(\tau)
\right|^2
d\tau.
\end{align}
Next we will estimate the $d\xi$ integral in \eqref{X-norm-v1}. In fact, we have the following estimate:
\begin{lemma}
\label{v1-mult-lema}
For any $\varepsilon>0$, if $s\geq -j-\frac12-\varepsilon$, $b\ge 0$ then we have
\begin{equation}
\label{v1-mult-ine-1}
\int_{-\infty}^{\infty}
(1+|\xi|)^{2s}
(1+|\tau-\xi^m|)^{2b}
\left|
\xi^j
F(\xi,\tau)
\right|^2
d\xi
\leq
c_{s,b,m}
|\tau|^{\frac{2(s+j)+2mb-1+2\varepsilon}{m}},
\end{equation}
where $c_{s,b,m}$ is a constant depending on $s$, $b$ and $m$. 
\end{lemma}
\noindent
We prove  \eqref{v1-mult-ine-1} below. Next,
combining  it with \eqref{X-norm-v1},
 we obtain
\begin{align*}
\|v_1\|_{X^{s,b}}^2
\le
c_{s,b,m}
\int_{-\infty}^{-1}
\Big|
\tau^{\frac{s+mb-\frac12-\ell+\varepsilon}{m}}
\widehat h_\ell(\tau)
\Big|^2
d\tau
\le
c_{s,b,m} \|h_\ell\|_{H_t^{\frac{s+mb-\frac12-\ell+\varepsilon}{m}}}^2,
\end{align*}
which is the desired estimate \eqref{reduced-ibvp-Bourgain-est} for $v_1$.
\,\,
$\Box$

\nin
{\bf Proof of Lemma \ref{v1-mult-lema}.} 
Using the estimate  \eqref{F-bound-ine},
denoting the integrand by
\begin{align*}
I(\xi,\tau)
\doteq
(1+|\xi|)^{2s}
(1+|\tau-\xi^m|)^{2b}
\left|
\xi^j
F(\xi,\tau)
\right|^2
\lesssim
(1+|\xi|)^{2s}
(1+|\tau-\xi^m|)^{2b}
\frac{\xi^{2j}|\tau|^{\frac{2n-2}{m}}}{|\xi|^{2n}+|\tau|^{2n/m}},
\end{align*}
and using our assumption $b\ge0$, $|\tau|>1$, we have
\begin{equation}
\label{mult-b-est}
1
\leq
|\xi|+|\tau|^\frac1m
\Rightarrow
1+|\tau-\xi^m|
\leq
2(|\tau|+|\xi|^m)
\Rightarrow
(1+|\tau-\xi^m|)^b
\leq
c_{b,m}(
|\tau|^b
+
|\xi|^{mb}
),
\end{equation}  
where $c_{b,m}$ is a constant depending on $b$ and $m$. Hence  we obtain
\begin{equation}
\label{L2-est-mult}
I(\xi,\tau)
\le
c_{s,b,m}
\left[
(1+|\xi|)^{2s}|\tau|^{2b}+(1+|\xi|)^{2s}|\xi|^{2mb}
\right]
\cdot
\frac{\xi^{2j}|\tau|^{\frac{2n-2}{m}}}{|\xi|^{2n}+|\tau|^{2n/m}},
\end{equation}
where $c_{s,b,m}$ is a constant depending on $s$, $b$ and $m$.
Now we shall consider the following cases:

\vskip0.05in
\nin
$\bullet$ $|\xi|> |\tau|^\frac1m$
\quad
and 
\quad
$\bullet$ $|\xi|\le |\tau|^\frac1m$

\vskip0.05in
\nin
{\bf Case $|\xi|> |\tau|^\frac1m$.} Here we have $|\xi|\ge 1$ and choosing $n=(s+j)+mb+\frac12+\varepsilon$, from \eqref{L2-est-mult} we have
$$
I(\xi,\tau)
\lesssim
c_{s,b,m}
(1+|\xi|)^{2s}|\xi|^{2mb}
\cdot
\frac{\xi^{2j}|\tau|^{\frac{2n-2}{m}}}{|\xi|^{2n}}
\lesssim
c_{s,b,m}
|\xi|^{-1-2\varepsilon}
|\tau|^{\frac{2(s+j)+2mb-1+2\varepsilon}{m}}.
$$
Thus,  for the integral  \eqref{v1-mult-ine-1} we have the following inequality
$
\int_{|\xi|\ge 1}
I(\xi,\tau)
d\xi
\le
$
$
c_{s,b,m}
|\tau|^{\frac{2(s+j)+2mb-1+2\varepsilon}{m}},
$
which is the desired estimate \eqref{v1-mult-ine-1} in this case.

\nin
{\bf Case $|\xi|\le |\tau|^\frac1m$.}
 Then $1+|\xi|\lesssim|\tau|^{1/m}$. Choosing $n=(s+j)+\frac12+\varepsilon$, from \eqref{L2-est-mult} we have
$$
I(\xi,\tau)
\lesssim
c_{s,b,m}
(1+|\xi|)^{2s}|\tau|^{2b}
\cdot
\frac{\xi^{2j}|\tau|^{\frac{2n-2}{m}}}{|\xi|^{2n}+1}
\lesssim
(1+|\xi|)^{2s}|\tau|^{2b}
\cdot
\frac{|\tau|^{\frac{2n-2}{m}}}{(|\xi|+1)^{2n}}
\lesssim
(1+|\xi|)^{-1-2\varepsilon}
|\tau|^{\frac{2(s+j)+2mb-1+2\varepsilon}{m}}
$$
Therefore, integrating $\xi$ we get
$
\int_{|\xi|\le |\tau|^{1/m}}
I(\xi,\tau)
d\xi
\lesssim
|\tau|^{\frac{2(s+j)+2mb-1+2\varepsilon}{m}},
$
which is the desired estimate \eqref{v1-mult-ine-1} in this case.
This  completes the proof of Lemma \ref{v1-mult-lema}.\,\,
$\Box$

\nin
{\bf Estimation of $D^\alpha$ norm.}
Now we estimate the second part in the modified Bourgain norm $\|\cdot \|_{s,b,\alpha}$. More precisely we have the result.
\begin{lemma}
\label{alpha-est-lem} 
If 
$s>-j-\frac12$
and $\frac12<\alpha\le\frac{1}{2}+\frac1m(s+j+\frac12)$, then we have
\begin{equation}
\label{alpha-est-I2}
\|v_1\|_{D^\alpha}^2
\doteq
\int_{-\infty}^{\infty}
\int_{-1}^1
(1+|\tau|)^{2\alpha}
|\widehat{v}_1(\xi,\tau)|^2
d\xi
d\tau 
\lesssim
\|h_\ell\|_{H^{\frac{s+j-\ell}{m}}(\rr)}^2.
\end{equation}
\end{lemma}
\nin
{\bf Proof of Lemma \ref{alpha-est-lem}.} First, we recall that the Fourier transform of $v_1$ is given by \eqref{v1-xt-FT}.
Thus we have
\begin{align*}
\|v_1\|_{D^\alpha}^2
=
\int_{-\infty}^{-1}
\int_{-1}^1
(1+|\tau|)^{2\alpha}
|\xi^j
F(\xi,\tau)
\tau^{-(\ell+j)/m}
\widehat h_\ell(\tau)|^2
d\xi
d\tau,
\end{align*}
where 
$
F(\xi,\tau)
=
\int_{x\in\rr}
e^{-i\xi x}\eta(-\gamma_I\tau^\frac1m x)
dx.
$
Also, applying Lemma \ref{F-bound-lem} with $n=0$, we get
the following bound for $F$, that is
$
|F(\xi,\tau)|
\le
|\tau|^{-1/m}.
$
Hence, after integrating $\xi$, we have
\begin{align*}
\|v_1\|_{D^\alpha}^2
\lesssim
\int_{-\infty}^{-1}
(1+|\tau|)^{2\alpha}
|
\tau^{-\frac{1+j+\ell}{m}}
\widehat h_\ell(\tau)|^2
d\tau
\lesssim
\int_{-\infty}^{-1}
(1+|\tau|)^{\frac{2m\alpha-2-2j-2\ell}{m}}
|
\widehat h_\ell(\tau)|^2
d\tau
\le
\|h_\ell\|_{H_t^{\frac{m\alpha-j-1-\ell}{m}}}^2.
\end{align*}
Choosing $\alpha$ such that 
$
\alpha
\le
\frac{1}{2}+\frac1m\big(s+j+\frac12\big),
$
we get the desired estimate \eqref{alpha-est-I2}
.
\,\,
$\square$

\medno
{\bf Bound near $\xi=0$.}
Next, we estimate the Bourgain norm for $v_0$. 
We begin with extending  it from $\rr^+\times(0,2)$ to $\rr\times\rr$ (keeping the same notation)
\begin{equation}
\label{def-v0-extension}
v_0(x,t)
\doteq
\int_0^1 e^{i\gamma\xi \varphi_1(x)}e^{-i\xi ^mt}\xi ^{2j-\ell}\,\widehat {h_\ell}(-\xi ^m)d\xi,
\quad
x\in\rr,
\,\,
t\in\rr,
\end{equation}
where $\varphi_1(x)$ is a smooth version of 
$|x|$. More precisely 

\vskip-0.05in
\begin{minipage}{0.4\linewidth}
\begin{center}
\begin{tikzpicture}[yscale=0.6, xscale=1]
%%%%%%%%%%%%%%%%%%%
%
%Variables defined
%
%%%%%%%%%%%%%%%%%%%%
\newcommand\X{0};
\newcommand\Y{0};
\newcommand\FX{11};
\newcommand\FY{11};
\newcommand\R{0.6};
\newcommand*{\TickSize}{2pt};
%%%%%%%%%%%%%%%%%%
%
%End
%
%%%%%%%%%%%%%%%%%%%%%
\draw[black,line width=1pt,-{Latex[black,length=2mm,width=2mm]}]
(-2.5,0)
--
(2.5,0)
node[above]
{\fontsize{\FX}{\FY}\bf \textcolor{black}{$x$}};

\draw[black,line width=1pt,-{Latex[black,length=2mm,width=2mm]}]
(0,-0.5)
--
(0,3)
node[right]
{\fontsize{\FX}{\FY}\bf \textcolor{black}{$y$}};

\draw[line width=1pt, yscale=1,domain=-2.1:-1.2,smooth,variable=\x,red]  plot ({\x},{-\x});

\draw[line width=1pt, yscale=1,domain=0:2.1,smooth,variable=\x,red]  plot ({\x},{\x});

\draw[smooth,line width=1pt, red]
(-0.4,-0.2)
to[out=5,in=-135]
(0,0)
;

\draw[smooth,line width=1pt, red]
(-1.2,1.2)
to[out=-45,in=175]
(-0.4,-0.2)
;

\draw[red,dashed, line width=0.5pt]
(2,1.3)
node[]
{\fontsize{\FX}{\FY}$\varphi_1(x)$}

(0,0)
node[yshift=-0.2cm,xshift=0.2cm]
{\fontsize{\FX}{\FY}$0$}

(-1.5,1.5)
--
(-1.5,0)
node[yshift=-0.2cm]
{\fontsize{\FX}{\FY}$-1$};

\end{tikzpicture}
\end{center}
\end{minipage}
\begin{minipage}{0.6\linewidth}
\begin{equation}
\label{def-phi3}
\varphi_1(x)
=
\begin{cases}
x,
\quad
x\geq
0
\\
-x,
\quad
x\leq
-1
\\
\text{ smooth on } \rr.
\quad
\end{cases}
\end{equation}
\end{minipage}

\nin
For $v_0$,  we have  the $L^2$ estimate
\begin{align}
\label{v0-L2-est}
\|\psi(t) v_0\|^2_{s,b,\alpha}
\hskip-0.05in
\lesssim
\|\psi(t) v_0\|^2_{L^2_{x,t}}
+
\|\p_x^{n_1}[\psi(t) v_0]\|^2_{L^2_{x,t}}
+
\|\p_t^{n_2}[\psi(t) v_0]\|^2_{L^2_{x,t}}
+
\|\p_x^{n_1}\p_t^{n_2}[\psi(t) v_0]\|^2_{L^2_{x,t}},
\end{align}
where 
$n_1
=
n_1(s,b)
\doteq
2m\lfloor |b|\rfloor+2\lfloor |s|\rfloor+2
$
and
$
n_2
=
n_2(b,\alpha)
\doteq
2\lfloor |b|\rfloor+2\lfloor|\alpha|\rfloor+2.
$
Using the $L^2$ boundedness of  Laplace transform (see Lemma 2.3 in  \cite{fhm2016} or  \cite{hardy1933})  we get   
\begin{align}
\label{v0-L2-Laplace-est}
\|\p_x^{n_1}\p_t^{n_2}[\psi\cdot v_0]\|_{L^2_{x,t}}^2
\le
C_{n_1, n_2}
\int_0^{\infty}
|
\widehat{h_\ell}(\tau)
|^2
d\tau
\lesssim
\|h_\ell\|_{H_t^{\frac{s+j-\ell}{m}}}^2
\quad
\forall n_1,n_2\in \NN_0.
\end{align}
\nin
{\bf End of Proof for Theorem \ref{reduced-pure-ibvp-thm}.}
Combining Lemma \ref{reduced-pibvp-thm1}, Lemma \ref{alpha-est-lem} and estimate \eqref{v0-L2-est} with \eqref{v0-L2-Laplace-est}, we get estimate \eqref{reduced-pure-ibvp-B-es}. This completes the proof of Theorem \ref{reduced-pure-ibvp-thm}.
\,\,
$\Box$

%
%
%%%%%%%%%%%%%%%%%%%%%%%%%%%%%%%%%
%
%  
%
%     Decomposition and proof of forced linear  ibvp estimates 
%
%
%
%%%%%%%%%%%%%%%%%%%%%%%%%%%%%%%%%
%
%
%
\section{
Proof of forced linear  ibvp estimates
 }
\label{proof-half-line}
In this section, we prove the basic linear estimate \eqref{forced-linear-kdvm-est}.
We begin by decomposing  the forced linear ibvp \eqref{LKdVm}
into a homogeneous ibvp (A) and  an inhomogeneous ibvp 
with zero data  (B). 
Then, we decompose both problems further in a  convenient  way
simplifying both their Fokas solution formula and its estimation.

\medno
{\bf A. The homogeneous linear ibvp:} 
\begin{subequations}
\label{homo-ibvp}
\begin{align}
&\p_tu+(-1)^{j+1}\p_x^{2j+1}u
=
0, 
%&&x\in(0,\infty),\ t\in(0,T),
\\
\label{homo-ibvp:ic}
&u(x,0) = u_0(x)\in H_x^s(0,\infty),  
%&&x\in(0,\infty),
\\
\label{homo-ibvp:bc}
&\p_x^\ell u(0,t) = g_\ell(t)
\in H_t^{\frac{s+j-\ell}{m}}(0,T), 
&&
\ell
=
0,1,\dots,j-1,
%\quad
%t\in(0,T),
\end{align}
\end{subequations}
with solution denoted by  
$
u(x,t)
\doteq
S\big[u_0, g_0,\dots,g_{j-1}; 0\big](x,t)
$ 
and which is defined in \eqref{UTM-sln-compact}.
We decompose it  further into the following two problems.

\nin
%\advance\itemsep 2mm
%
{\bf A$_1$. The homogeneous linear ivp:}
\begin{subequations}
\label{homo-ivp}
\begin{align}
\label{homo-ivp:eq}
&\p_tU+(-1)^{j+1}\p_x^{2j+1}U=0,  
%&& x\in \mathbb R,\ t\in (0, T),
\\
\label{homo-ivp:ic}
&
U (x,0)
= 
U_0(x)\in H_x^s(\mathbb R), 
%\quad 
%&& x\in \mathbb R, 
\end{align}
\end{subequations}
where $U_0\in H_x^s(\mathbb R)$ is an extension of the initial datum $u_0\in H_x^s(0, \infty)$ such that
\begin{equation}
\label{assum-U0}
\|U_0\|_{H_x^s(\mathbb R)} 
\leqslant 
2\|u_0
\|_{H_x^s(0, \infty)}
\end{equation}
with its solution given by
\begin{equation}
\label{homo-ivp:sln}
U (x,t) 
= 
S\big[U_0; 0\big] (x, t) 
= 
\frac{1}{2\pi} \int_{\rr} e^{i\xi x+i\xi^mt}\, \widehat{U}_0(\xi) d\xi,
\end{equation}
where 
$
\widehat{U}_0(\xi) 
=
\int_{\rr} e^{-i\xi x}\, U_0(x) dx,\quad \xi \in \mathbb R
$.

\nin
{\bf A$_2$. The homogeneous linear ibvp with
 zero initial data:}
\begin{subequations}
\label{pure-ibvp}
\begin{align}
\label{pure-ibvp:eqn}
&
\p_tu+(-1)^{j+1}\p_x^{2j+1}u
=
0, 
 \quad 
&&
%x\in(0,\infty),
%\,\,
%t\in(0,T),
\\
\label{pure-ibvp:ic}
&u(x,0)= 0, 
%&&x\in (0,\infty),  
\\
\label{pure-ibvp:bc}
&
\p_x^\ell u(0,t) = g_\ell(t) - \p_x^\ell U (0,t)
 \doteq G_\ell(t)
 \in H_t^{\frac{s+j-\ell}{m}}(0,T), , 
 \quad 
&&
\ell
=
0,1,\dots, j-1,
%\quad
%t\in (0, T), 
\end{align}
\end{subequations}
with solution 
$
u(x, t) 
\doteq
S\big[0, G_0,\dots,G_{j-1}; 0\big](x, t)
$, which is defined in \eqref{UTM-sln-compact}.

\medno
{\bf B. The forced linear ibvp with zero data:}
\begin{subequations}
\label{forced-ibvp}
\begin{align}
&
\p_tu+(-1)^{j+1}\p_x^{2j+1}u
=
f(x,t),
%\quad  &&x\in(0,\infty),\ t\in(0,T),
\\
\label{forced-ibvp:ic}
&u(x,0) = 0, 
%\quad &&x\in(0,\infty),
\\
\label{forced-ibvp:bc}
&\p_x^\ell u(0,t) =0, 
\quad
&&
\ell
=
0,1,\dots,j-1,
%\quad t\in(0,T),
\end{align}
\end{subequations}
whose solution 
$
u(x,t)
\doteq
S\big[0, 0,\dots,0; f\big](x,t)
$  is defined in \eqref{UTM-sln-compact}.
This problem can be further decomposed into the following two problems (B$_1$) and (B$_2$).

\nin
%\advance\itemsep 2mm
%
{\bf B$_1$. The forced linear ivp with zero initial data:}
\begin{subequations}
\label{forced-ivp}
\begin{align}
\label{forced-ivp:eqn}
&\p_tW+(-1)^{j+1}\p_x^{2j+1}W
= 
w(x, t),\, 
%&& x\in \mathbb R,\ t\in (0, T),
\\
\label{forced-ivp:ic}
& W(x,0)
=
0, 
%&& x\in \mathbb R, 
\end{align}
\end{subequations}
where $w$ is an {\it extension of the forcing $f$} such that
\begin{align}
\label{assum-F-extension}
&\|w\|_{X^{s,-b,\alpha-1}(\rr^2)}
\leqslant
2\|f\|_{X^{s,-b,\alpha-1}(\rr^+\times(0,T))},
\hskip.2in
-1
\le
s
\le
\frac12,
\\
&\|w\|_{X^{s,-b,\alpha-1}(\rr^2)}
+
\|w\|_{Y^{s,-b}(\rr^2)}
\leqslant
2(
\|f\|_{X^{s,-b,\alpha-1}(\rr^+\times(0,T))}
+
\|f\|_{Y^{s,-b}(\rr^+\times(0,T))}
),
s\not\in[-1,1/2],
\nn
\end{align}
where $Y^{s,b}$ is defined in \eqref{Ysb-def}. The solution of this problem is given by Duhamel's formula
\begin{align}
\label{Duhamel-1}
W(x,t) 
\doteq
S\big[0; w\big](x, t)
=&
-\frac{i}{2\pi} \int_{\rr} \int_{0}^t  e^{i\xi x+i\xi^m(t-t')}
\widehat w(\xi, t') dt' d\xi,
\\
\label{Duhamel-2}
=&
-i \int_{0}^t  S\big[w(\cdot, t'); 0\big](x, t-t') dt',
\end{align}
where $\widehat{w}$ is the Fourier transform of $w$ with respect to $x$, 
and  $S\big[w(\cdot, t'); 0\big]$ in the Duhamel representation \eqref{Duhamel-2} denotes the solution \eqref{homo-ivp:sln} of 
 ivp \eqref{homo-ivp} (that is Problem A$_1$) with $w(x, t')$ in place of the initial data  and zero forcing.

\nin
{\bf B$_2.$ The homogeneous linear ibvp with  zero initial data:}
\begin{subequations}
\label{forced-pure-ibvp}
\begin{align}
\label{forced-pure-ibvp:eqn}
&\p_tv+(-1)^{j+1}\p_x^{2j+1}v
= 
0, 
%&&x\in (0,\infty), \ t\in (0,T),
\\
\label{forced-pure-ibvp:ic}
&v(x,0)
= 
0, 
%&&x\in (0,\infty),   
\\
\label{forced-pure-ibvp:bc}
&
\p_x^\ell v(0,t) = -\p_x^\ell W(0, t) \doteq -W_\ell(t)
&&
\ell
=
0,1,\dots,j-1,
%\quad t\in (0, T), 
\end{align}
\end{subequations}
whose solution  
$
v(x, t) 
\doteq
S\big[0, -W_0,\dots,-W_{j-1}; 0\big](x, t)
$
 is defined in \eqref{UTM-sln-compact}.

Next we describe the estimates for each one of the above 
sub-problems.
\begin{theorem}
[Estimates for homogeneous ivp  A$_1$]
\label{homo-ivp-thm}
The solution   $U=S\big[U_0; 0\big]$ to ivp   \eqref{homo-ivp} defined by formula \eqref{homo-ivp:sln} 
satisfies the  space estimate
\begin{align}
\label{homo-ivp-se}
\sup_{t\in [0, T]} \|S\big[U_0; 0\big](t)\|_{H_x^s(\mathbb R)} = \|U_0\|_{H_x^s(\mathbb R)}, &&s\in\mathbb R,
\end{align}
and the time estimate for its $\ell$-th derivative (\it needed to have boundary data in desired space)
\begin{align}
\label{homo-ivp-te}
\sup_{x\in\mathbb R}\|\psi(t)\p_x^\ell S\big[U_0; 0\big](x)\|_{H_t^{\mu_\ell}(\rr)} 
\leqslant 
c_s \|U_0\|_{H_x^s(\mathbb R)},
\,\,
s\in\rr,
\end{align}
where $\ell=0,1,\dots,j-1$ and  
$\mu_\ell=\frac 1m(s+j-\ell)$.
Also, it satisfies the following estimate in modified Bourgain spaces 
\begin{equation}
\label{homo-ivp-se-bourgain}
\|\psi(t)S\big[U_0; 0\big](x,t)\|_{X^{s,b,\alpha}}
\leq
c_{\psi}\|U_0\|_{H^s(\rr)},
\quad
\forall s, b,\alpha\in\rr,
\end{equation}
where $c_\psi$ is a constant depending only on $\psi$.
Here and elsewhere in this paper
$\psi$ is a cutoff function
in $C^{\infty}_0(-1, 1)$ such that 
$0\le \psi \le 1$ and $\psi(t)=1$ for $|t|\le 1/2$.
\end{theorem}
\nin
{\bf Proof of Theorem \ref{homo-ivp-thm}.}
The proof of   the space  estimate \eqref{homo-ivp-se}
is straightforward.  The proof  of  the time  estimate 
\eqref{homo-ivp-te}  is similar to that for KdV, which can be
 found in Holmer \cite{h2006}, and Colliander and Kenig \cite{ck2002}.
Finally, estimate  \eqref{homo-ivp-se-bourgain}
 follows from inequality
$$
\int_{-\infty}^\infty
\int_{-\infty}^\infty
\chi_{|\xi|< 1}(1+|\tau|)^{2\alpha}
|\widehat{u}(\xi,\tau)|^2
d\xi
d\tau
\lesssim
\|u\|^2_{X^{s,\alpha}},
$$
and  the estimate
\begin{equation}
\label{pure-ivp-bourgain-est-no2}
\|\psi(t)S\big[U_0; 0\big](x,t)\|_{X^{s,b}}
\leq
c_{b}\|U_0\|_{H^s},
\quad
\text{where} 
\quad
c_b=\|\psi\|_{H^b},
\quad
s, b \in \rr,
\end{equation}
whose proof  can be found in \cite{fhy2020}.
\, $\square$

%
%%%%%%%%%%%%%%%%%%%%%
%
%
%
%    Estimates for pure ibvp  A$_2$
%
%
%
%%%%%%%%%%%%%%%%%%%%%
%
%
%
\begin{theorem}
[Estimates for pure ibvp  on the half-line]
\label{pure-ibvp-thm}
Let $s\neq \frac12,\frac32,\dots, j-\frac12$.
The solution of the pure ibvp \eqref{pure-ibvp}  satisfies the  space estimate
\begin{align}
\label{pure-ibvp-thm:se}
\sup_{t\in \rr} \|S\big[0, G_0,\dots,G_{j-1}; 0\big]\|_{H_x^s(0,\infty)} 
\leq
c_{s,m} \sum\limits_{\ell=0}^{j-1}\|G_\ell\|_{H_t^{\mu_\ell}(0,T)}
,
\quad
s\ge 0.
\end{align}
Also, for $b\in[0,\frac12)$, 
$\frac12<\alpha\le\frac{1}{2}
+
\frac 1m (s+j+\frac12)
$
its solution satisfies the estimate in Bourgain spaces
\begin{align}
\label{pure-ibvp-thm:se-bourgain}
\|S\big[0, G_0,\dots,G_{j-1}; 0\big]\|_{X^{s, b,\alpha}(\rr^+\times (0,T))} 
\le
c_{s,b,\alpha}
\sum\limits_{\ell=0}^{j-1}\|G_\ell\|_{H_t^{\frac{s+j-\ell}{m}}(0,T)},
\quad
s>-j-\frac12.
\end{align}
\end{theorem}
\noindent
{\bf 
Proof of Theorem \ref{pure-ibvp-thm}
}
The proof of   the space  estimate \eqref{pure-ibvp-thm:se}
can be found in \cite{y2020} for KdVm,
and in \cite{fhm2016} for KdV.
Here, we prove estimate 
\eqref{pure-ibvp-thm:se-bourgain}, which is new. 
We do this by transforming problem A$_2$ to the reduced pure ibvp \eqref{LKdVm-reduced}.
For this, we extend  $G_\ell$ from $(0,T)$ to
 a function $h_\ell$ on $\rr$ supported in $[0, 2]$ 
and such that
$
\|h_\ell\|_{H_t^{\mu_\ell}(\mathbb R)}
 \lesssim
\|G_\ell\|_{H_t^{\mu_\ell}(0,T)}
$
for any $\ell=0,1,\dots,j-1$,
via the following result, whose proof can be found 
in \cite{lmbook1972,WMbook2000,fhm2017,y2020}.
\begin{lemma}
\label{Extension-lemma}
For a general function $h^*(t)\in H_t^s(0,2)$, $s\ge 0$,  let 
the extension
\begin{align*}
\tilde h^*(t)
\doteq
\begin{cases}
h^*(t),
\quad
t\in(0,2),
\\
0,
\quad
\text{elsewhere}.
\end{cases}
\end{align*} 
If $0\leq s<\frac12$,  then  the extension
$\tilde h^*\in H^s(\rr)$ and 
for some $c_s>0$ we have
\begin{align}
\label{char-mult-est}
\|\tilde h^*\|_{H_t^s(\rr)}
\leq
c_s\|h^*\|_{H_t^s(0, 2)}.
\end{align}
If  $\frac12<s\leq1$, then for estimate  \eqref{char-mult-est} to hold
we must have the condition  
\begin{align}
\label{vanishing-condition}
h^*(0)=h^*(2)=0.
\end{align}
\end{lemma}
\vskip-0.08in
Also, we shall need the following multiplier by a characteristic
estimate from Holmer  \cite{h2006}.
\begin{align}
\label{ext-est-neg-new}
\|\chi_{(0,\infty)}g\|_{H^s(\rr)}
\le
c_s\|g\|_{H^s(\rr)},
\,\,
 g\in H^s(\rr), \,\,  -\frac12<s<\frac12.
\end{align}

\nin
Since $h_\ell$ extends $G_\ell$ from $(0,T)$ to $\rr$,
 we have 
$
S\big[0, G_0,\dots,G_{j-1}; 0\big] (x,t)
=
S\big[0, h_0,\dots,h_{j-1}; 0\big] (x,t),
$
 for $x\in\rr^+,
t\in(0,T)$.
In fact,  by Theorem \ref{reduced-pure-ibvp-thm},  for $\frac12>s>-j-\frac12$, $b\in[0,\frac12)$ and  $\frac12<\alpha\le\frac{1}{2}+\frac1m(s+j+\frac12)$ we get 
\begin{align}
\label{ibvp-identity}
&\|S\big[0, G_0,\dots,G_{j-1}; 0\big]\|_{X^{s, b,\alpha}(\rr^+\times (0,T))} 
=
\|S\big[0, h_0,\dots,h_{j-1}; 0\big]\|_{X^{s, b,\alpha}(\rr^+\times(0,T))}
\\
\le&
\|S\big[0, h_0,\dots,h_{j-1}; 0\big]\|_{X^{s, b,\alpha}(\rr^+\times(0,2))}
\lesssim
\sum\limits_{\ell=0}^{j-1}\|h_\ell\|_{H_t^{\frac{s+j-\ell}{m}}(\mathbb R)}
\lesssim
\sum\limits_{\ell=0}^{j-1}\|G_\ell\|_{H_t^{\frac{s+j-\ell}{m}}(0,T)},
\nonumber
\end{align}
which is the desired estimate  \eqref{pure-ibvp-thm:se-bourgain}. 
This completes the proof  of Theorem \ref{pure-ibvp-thm}.
\,\,
$\Box$

%
%%%%%%%%%%%%%%%%%%%%%
%
%
%     Estimates for forced ivp  B$_1$
%
%
%
%%%%%%%%%%%%%%%%%%%%%
%
%
%
\begin{theorem}
[Estimates for forced ivp  B$_1$]
\label{forced-ivp-thm}
The solution  $W=S\big[0; w\big]$  of the forced ivp
\eqref{forced-ivp}
defined by equations 
\eqref{Duhamel-1}--\eqref{Duhamel-2} satisfies the following estimate in modified Bourgain spaces
\begin{equation}
\label{forced-ivp-bourgain-est}
\|\psi(t)S\big[0; w\big](x,t)\|_{X^{s,b,\alpha}}
\lesssim
c_\psi
\|w\|_{s,-b,\alpha-1},
\quad
s\in\rr,
\,\,
0< b<\frac12<\alpha<1,
\end{equation}
and the time estimate (\it needed to have boundary data in desired space)
\begin{equation}
\label{forced-ivp-te}
\sup\limits_{x\in\rr}
\|
\psi(t)
\p_x^\ell S[0,w]
\|_{H_t^{\frac{s+j-\ell}{m}}(\rr)}
\le
\begin{cases}
c_{s,m,\ell,b}\|w\|_{X^{s,-b}},
\quad
-1\le
s
\le
\frac12,
\\
c_{s,m,\ell,b}
(
\|w\|_{X^{s,-b}}
+
\|w\|_{Y^{s,-b}}
),
\quad
s\in\rr
,
\end{cases}
\end{equation}
where $X^{s,b}$ is the Bourgain space
\eqref{def-Bourgain} and 
$Y^{s,b}$ is a ``temporal" Bourgain space defined by \eqref{Ysb-def}.
\end{theorem}
\vskip-0.05in
\noindent
{\bf
Proof of Theorem \ref{forced-ivp-thm}.}
First, we prove estimate \eqref{forced-ivp-bourgain-est}. 
For the $X^{s,b}$ part of  
$\|\psi\cdot S\big[0; w\big]\|_{X^{s,b,\alpha}}$,
i.e.
$
\|\psi(t)S\big[0; w\big](x,t)\|_{X^{s,b}},
$
we  have the next basic estimate, whose proof can be found in \cite{fhy2020} 
\begin{align}
\label{T-bound}
\|
\psi(t)S\big[0; w\big](x,t)
\|_{X^{s,b}}^2
\lesssim
\|w\|_{X^{s,b-1}}^2
+
\int_{\mathbb{R}}(1+|\xi|)^{2s}
\left(
\int_{\mathbb{R}}
\frac{|\widehat{w}(\xi,\tau)|}{1+|\tau-\xi^{m}|}
d\tau
\right)^2d\xi,
\,\,
0< b<1.
\end{align}
Since $b-1<-\frac12<-b$ we get
$
\|w\|_{X^{s,b-1}}
\le
\|w\|_{X^{s,-b}}.
$
For the second term in \eqref{T-bound}, 
writing 
$
1+|\tau-\xi^{m}|
=
(1+|\tau-\xi^{m}|)^{1-b}(1+|\tau-\xi^{m}|)^b
$
and applying the Cauchy-Schwartz inequality for
 the  $\tau$-integral we obtain:
\begin{align}
\label{2nd-bilinear-est}
\hskip-0.05in
\int_{\mathbb{R}}(1+|\xi|)^{2s}
\Big(\int_{\mathbb{R}}
\frac{|\widehat{w}(\xi,\tau)|}{1+|\tau-\xi^{m}|}
d\tau\Big)^2d\xi
\le
c_b
\int_{\mathbb{R}}(1+|\xi|)^{2s}
\hskip-0.05in
\int_{\mathbb{R}}
\frac{|\widehat{w}(\xi,\tau)|^2}{(1+|\tau-\xi^{m}|)^{2b}}d\tau
d\xi
\simeq
\|w\|_{X^{s,-b}}^2.
\end{align}
For  the $\alpha$-part of  the norm
$\|\psi\cdot S\big[0; w\big]\|_{X^{s,b,\alpha}}$,
 using the fact that $|\xi|\le 1$, which
gives  
$
(1+|\tau|)^{2\alpha}
\simeq
(1+|\tau-\xi^m|)^{2\alpha},
$
we get
\begin{align}
\label{D-alpha-est}
\int_{\rr}
\int_{-1}^1
(1+|\tau|)^{2\alpha}
|
\reallywidehat{
\psi S\big[0; w\big]
}(\xi,\tau)|^2
d\xi
d\tau
\simeq
\int_{\rr}
\int_{\rr}
(1+|\tau-\xi^m|)^{2\alpha}
|
\chi_{|\xi|\le1}
\reallywidehat{
\psi S\big[0; w\big]
}(\xi,\tau)|^2
d\xi
d\tau.
\end{align}
Also, since 
$
\reallywidehat{
\psi S\big[0; w\big]
}^x(\xi,t)
=-i
\psi(t)
\int_{0}^t  e^{i\xi^m(t-t')}
\widehat w(\xi, t') dt',
$
we obtain
$$
\chi_{|\xi|\le 1}
\reallywidehat{
\psi S\big[0; w\big]
}(\xi,\tau)
=
\reallywidehat{
\psi S\big[0; w_1\big]
}(\xi,\tau),
$$
where  
$
\widehat{w_1}^x(\xi,t)
\doteq
\chi_{|\xi|\le 1}
\widehat{w}^x(\xi,t).
$
Using $w_1$ notation, from \eqref{D-alpha-est} we have
\begin{align}
\label{w1-express}
\int_{\rr}
\int_{-1}^1
(1+|\tau|)^{2\alpha}
|
\reallywidehat{
\psi S\big[0; w\big]
}(\xi,\tau)|^2
d\xi
d\tau
=
\|
\psi(t)S\big[0; w_1\big](x,t)
\|_{X^{0,\alpha}}^2.
\end{align}
For $\|
\psi(t)S\big[0; w_1\big](x,t)
\|_{X^{0,\alpha}}$, applying  estimate \eqref{T-bound} 
with $s=0$ and $b=\alpha>\frac12$,  we get 
\begin{align}
\label{w1-express-est}
\|
\psi(t)S\big[0; w_1\big](x,t)
\|_{X^{0,\alpha}}^2
\lesssim
\int_{\rr}
\int_{-1}^1
(1+|\tau-\xi^m|)^{2\alpha-2}
|
\reallywidehat{
w
}(\xi,\tau)|^2
d\xi
d\tau.
\end{align}
Using  estimates \eqref{w1-express}, \eqref{w1-express-est} and
the fact that $|\xi|\le 1$ again,
we obtain
\begin{align}
\label{alpha-est-fin}
\int_{\rr}
\int_{-1}^1
(1+|\tau|)^{2\alpha}
|
\reallywidehat{
\psi S\big[0; w\big]
}(\xi,\tau)|^2
d\xi
d\tau
\lesssim
\int_{\rr}
\int_{-1}^1
(1+|\tau|)^{2\alpha-2}
|
\reallywidehat{
w
}(\xi,\tau)|^2
d\xi
d\tau.
\end{align}
Combining  \eqref{T-bound} and \eqref{2nd-bilinear-est} with \eqref{alpha-est-fin} and taking into consideration 
that $b-1<-b$, we get
\begin{align*}
\|\psi(t)S\big[0; w\big](x,t)\|_{X^{s,b,\alpha}}^2
\lesssim&
\|w\|_{X^{s,b-1}}^2
+
\|w\|_{X^{s,-b}}^2
+
\int_{\rr}
\int_{-1}^1
(1+|\tau|)^{2\alpha-2}
|
\reallywidehat{
w
}(\xi,\tau)|^2
d\xi
d\tau
\\
\lesssim&
\|w\|_{s,-b,\alpha-1}^2,
\quad
s\in\rr,
\,\,
0\le b<\frac12<\alpha<1.
\end{align*}
This completes the proof of estimate \eqref{forced-ivp-bourgain-est}.

%
%%%%%%%%%%%%%%%%%%%%%
%
%
%    Time estimates for forced ivp  B$_1$
%
%
%
%%%%%%%%%%%%%%%%%%%%%
%
%
%
\noindent
{\it\bf Proof of estimate \eqref{forced-ivp-te}.}
For the KdV equation, this estimate was proved in \cite{h2006}
 (see Lemma 5.6). Also a similar estimate (for $s=0$)  was proved 
 in  \cite{ck2002} (see Lemma 5.5).   Differentiating the solution
 formula   \eqref{Duhamel-1}, i.e. 
 $
 S\big[0; w\big](x, t)
=
-\frac{i}{2\pi} \int_{\xi\in \mathbb R} \int_{t'=0}^t  e^{i\xi x+i\xi^m(t-t')}
\widehat w(\xi, t') dt' d\xi
 $
$\ell$ times with respect to $x$ and decomposing it (like in \cite{fhy2020}), we obtain the Bourgain writing
\begin{align}
\label{Tfg2-term-recall}
\psi(t)\p_x^\ell S[0,w](x,t)
\simeq&
\psi(t)
\int_{\mathbb{R}}\int_{\mathbb{R}}
e^{i(\xi x+\tau t)}
\frac{1-\psi(\tau-\xi^m)}{\tau-\xi^m}
\xi^\ell
\widehat{w}(\xi,\tau)  d\tau  d\xi \\
\label{Tfg3-term-recall}
-&
\psi(t)
\int_{\mathbb{R}}\int_{\mathbb{R}}
e^{i(\xi x+\xi^mt)}
\frac{1-\psi(\tau-\xi^m)}{\tau-\xi^m}
\xi^\ell
\widehat{w}(\xi,\tau)   d\tau  d\xi \\
\label{Tfg4-term-recall}
+&
\psi(t)
\int_{\mathbb{R}}\int_{\mathbb{R}}
e^{i(\xi x+\xi^m t)} 
\frac{\psi(\tau-\xi^m)[e^{i(\tau - \xi^m)t} - 1]}{\tau - \xi^m}
\xi^\ell
\widehat{w}(\xi,\tau)   d\tau  d\xi. 
\end{align}
 Next, we estimate each term above separately.  We start by estimating \eqref{Tfg3-term-recall}.
\vskip.05in
\nin
\underline{\it Estimate for \eqref{Tfg3-term-recall}.} For this term, we have 
$$
\eqref{Tfg3-term-recall}
\simeq
\psi(t)\p_x^\ell S[F_1,0],
\quad
\text{where }
\quad
\widehat{F}_1(\xi)
\simeq
\int_\rr
\frac{1-\psi(\tau-\xi^m)}{\tau-\xi^m}
\widehat{w}(\xi,\tau) d\tau.
$$
Using estimate \eqref{homo-ivp-te}, we get
\begin{align*}
\sup\limits_{x\in\rr}\|
\eqref{Tfg3-term-recall}
\|_{H_t^\frac{s+j-\ell}{m}(\rr)}^2
\lesssim&
\sup\limits_{x\in\rr}\|
\psi(t)
\p_x^\ell
S[F_1,0]
\|_{H_t^\frac{s+j-\ell}{m}(\rr)}^2
\lesssim
\|F_1\|_{H^s(\rr)}^2
\\
=&
\int_\rr
(1+|\xi|)^{2s}
\Big|
\int_\rr
\frac{1-\psi(\tau-\xi^m)}{\tau-\xi^m}
\widehat{w}(\xi,\tau)   
d\tau
\Big|^2
d\xi
\\
\lesssim&
\int_\rr
(1+|\xi|)^{2s}
\Big(
\int_\rr
\frac{|\widehat{w}(\xi,\tau)| }{1+|\tau-\xi^m|}
d\tau
\Big)^2
d\xi
\lesssim
\|w\|_{s,-b}^2,
\end{align*}
where in  the last step we used estimate \eqref{2nd-bilinear-est}.
This gives the desired estimate  \eqref{forced-ivp-te} for \eqref{Tfg3-term-recall}.

\nin
\underline{\it Estimate for \eqref{Tfg4-term-recall}.} For this term, using Taylor's expansion we have 
$$
\eqref{Tfg4-term-recall}
\simeq
\sum\limits_{k=1}^\infty
\frac{1}{k!}t^k
\p_x^\ell
\psi(t)S[c_k,0],
\,\,
\text{ where  }
\,\,
\widehat{c_k}(\xi)
\simeq
\int_{\rr}\psi(\tau-\xi^m)\cdot(\tau-\xi^m)^{k-1}
\widehat{w}(\xi,\tau)d\tau.
$$
Letting 
$\psi_k(t)\doteq t^k[\psi(t)]^{1/2}$ and 
using estimate \eqref{homo-ivp-te},  we get
\begin{align*}
\sup\limits_{x\in\rr}\|
&\eqref{Tfg4-term-recall}
\|_{H_t^\frac{s+j-\ell}{m}(\rr)}
\lesssim
\sum_{k=1}^{\infty}\frac{1}{k!}
\sup\limits_{x\in\rr}
\|
\psi_k(t)
\cdot
[\psi(t)]^{1/2}
\p_x^\ell
S[c_k,0]
\|_{H_t^\frac{s+j-\ell}{m}(\rr)}
\lesssim
\sum_{k=1}^{\infty}\frac{c_{\psi_k}}{k!}
\|c_k\|_{H^s}
\\
=&
\sum_{k=1}^{\infty}\frac{1}{k!}\left( \int_{\rr} (1+|\xi|)^{2s} 
\Big|\int_{\rr}\psi(\tau-\xi^m)\cdot(\tau-\xi^m)^{k-1}
\widehat{w}(\xi,\tau)d\tau
\Big|^2 d\xi \right)^{1/2},
\end{align*}
where $c_{\psi_k}\doteq \|\widehat{\psi_k}(\tau)(1+|\tau|)^\frac{s+j-\ell}{m}\|_{L^1}$ (like \eqref{c-eta-def} below).
Since the $\tau$-integration is over  $|\tau-\xi^m|\leq 1$ and  $|\psi(\tau-\xi^m)|\leq 1$
from the last relation we obtain that
\begin{align*}
\sup\limits_{x\in\rr}\|
\eqref{Tfg4-term-recall}
\|_{H_t^\frac{s+j-\ell}{m}(\rr)}
\lesssim&
\left(\int_{\mathbb{R}} (1+|\xi|)^{2s} \left(\int_{|\tau-\xi^m|\leq 1} 
|\widehat{w}(\xi,\tau)|d\tau\right)^2 d\xi \right)^{1/2}
\cdot
\sum_{k=1}^{\infty}\frac{1}{k!}
\\
\lesssim&
\left(\int_{\mathbb{R}} (1+|\xi|)^{2s} \left(\int_{\mathbb{R}} 
\frac{|\widehat{w}(\xi,\tau)|}{1+|\tau-\xi^m|}d\tau\right)^2 d\xi \right)^{1/2}
\lesssim
\|w\|_{s,-b},
\end{align*}
where last inequality  follows from  \eqref{2nd-bilinear-est}.
This gives the desired estimate \eqref{forced-ivp-te} for term \eqref{Tfg4-term-recall}. 

\medno
\underline{\it Estimate for \eqref{Tfg2-term-recall}.} 
We rewrite this term as
$$
\eqref{Tfg2-term-recall}
\simeq
\psi(t)
h(x,t),
\,\,
\text{where}
\,\,
h(x,t)
=
\int_{\mathbb{R}}\int_{\mathbb{R}}
e^{i(\xi x+\tau t)}
\frac{1-\psi(\tau-\xi^m)}{\tau-\xi^m}
\xi^\ell\widehat{w}(\xi,\tau)  d\tau  d\xi.
$$
Using the  property
$
\widehat{f\cdot g}
\simeq
\widehat{f}*\widehat{g}
$
we get
$
\widehat{\psi \cdot h}^t(x,\tau)
\simeq
\int_\rr
\widehat{\psi}(\tau-\tau_1)
\widehat{h}^t(x,\tau_1)
d\tau_1
$,
which combined with the  inequality
$
(1+|\tau|)^{\mu}
\le
(1+|\tau_1|)^{\mu}
(1+|\tau-\tau_1|)^{|\mu|},
$
for $\mu=\frac{s+j-\ell}{m}$, we get
\begin{align*}
\|\psi h\|_{H^\frac{s+j-\ell}{m}_t(\rr)}
\le&
\Big\|
\int_{\rr}
(1+|\tau_1|)^{\frac{s+j-\ell}{m}}(1+|\tau-\tau_1|)^{|\frac{s+j-\ell}{m}|}
\widehat{\psi}(\tau-\tau_1)
\widehat{h}^t(x,\tau_1)
d\tau_1
\Big\|_{L^2_\tau(\rr)}
\\
\le&
\|\widehat{\psi}(\tau)(1+|\tau|)^{|\frac{s+j-\ell}{m}|}\|_{L^1(\rr)}
\cdot
\Big\|
(1+|\tau|)^\frac{s+j-\ell}{m}
\widehat{h}^t(x,\tau)
\Big\|_{L^2_\tau(\rr)}
\\
=&
c_\psi
\Big\|
(1+|\tau|)^\frac{s+j-\ell}{m}
\widehat{h}^t(x,\tau)
\Big\|_{L^2_\tau(\rr)},
\end{align*}
where in the second step we use the Young's inequality for $r=q=2$, and $p=1$ and and we bound the constant
$c_\psi^2
\doteq
\|\widehat{\psi}(\tau)(1+|\tau|)^{|\frac{s+j-\ell}{m}|}\|_{L^1}^2$
 as follows
\begin{align}
\label{c-eta-def}
c_\psi^2
\leq
\!\!
\int_\rr
|\widehat{\psi}(\tau)|^2
(1+|\tau|)^{2|\frac{s+j-\ell}{m}|+2}
d\tau
\!\!
\int_\rr
(1+|\tau|)^{-2}
d\tau
\lesssim
\|\psi\|_{H_t^{|\frac{s+j-\ell}{2m}|+1}}^2.
\end{align}
Since
$
\widehat{h}^t(x,\tau)
\simeq
\int_{\mathbb{R}}
e^{i\xi x}
\frac{1-\psi(\tau-\xi^m)}{\tau-\xi^m}
\xi^\ell
\widehat{w}(\xi,\tau)  d\xi,
$
we have
\begin{align*}
\|\psi h\|_{H^\frac{s+j-\ell}{m}_t(\rr)}^2
\le&
c_\psi^2
\int_{\rr}
(1+|\tau|)^\frac{2(s+j-\ell)}{m}
\Big|
\int_{\mathbb{R}}
e^{i\xi x}
\frac{1-\psi(\tau-\xi^m)}{\tau-\xi^m}
\xi^\ell
\widehat{w}(\xi,\tau)  d\xi
\Big|^2
d\tau
\\
\le&
c_\psi^2
\int_{\rr}
(1+|\tau|)^\frac{2(s+j-\ell)}{m}
\Big(
\int_{\mathbb{R}}
\frac{1}{1+|\tau-\xi^m|}
\xi^\ell
|\widehat{w}(\xi,\tau)|  d\xi
\Big)^2
d\tau.
\end{align*}
Now we  consider the following two cases: 

\vskip0.05in
\noindent
$\bullet$ $-1\le s\le \frac12$
\quad
and 
\quad
$\bullet$ $s\not\in[-1,\frac12]$.

\vskip0.05in
\noindent
{\bf Case  $-1\le s\le \frac12$.} For this case, multiplying
and dividing $|\widehat{w}(\xi,\tau)|$ by $\frac{(1+|\xi|)^s}{(1+|\tau-\xi^m|)^{b}}$ and using Cauchy-Schwartz inequality for the integral of $d\xi$, we get
\begin{align*}
\|\psi h\|_{H^\frac{s+j-\ell}{m}_t(\rr)}^2
\le
c_\psi^2
\int_{\rr}
(1+|\tau|)^\frac{2(s+j-\ell)}{m}
G_1(\tau)
\int_{\mathbb{R}}
\frac{(1+|\xi|)^{2s}|\widehat{w}(\xi,\tau)|^2}{(1+|\tau-\xi^m|)^{2b}}
  d\xi
d\tau,
\end{align*}
where  
$
G_1(\tau)
\doteq
\int_{\rr}
\frac{\xi^{2\ell}}{(1+|\tau-\xi^m|)^{2-2b}(1+|\xi|)^{2s}}
d\xi.
$
Taking the sup norm in $\tau$ for $(1+|\tau|)^\frac{2(s+j-\ell)}{m}
G_1(\tau)$, we get
\begin{align*}
\|\psi h\|_{H^\frac{s+j-\ell}{m}_t(\rr)}^2
\le&
c_\psi^2
\Big\|
(1+|\tau|)^\frac{2(s+j-\ell)}{m}
G_1(\tau)
\Big\|_{L_\tau^\infty}
\int_{\rr}
\int_{\rr}
\frac{(1+|\xi|)^{2s}|\widehat{w}(\xi,\tau)|^2}{(1+|\tau-\xi^m|)^{2b}}
  d\xi
d\tau
\\
\le&
c_\psi^2
\Big\|
(1+|\tau|)^\frac{2(s+j-\ell)}{m}
G_1(\tau)
\Big\|_{L_\tau^\infty}
\|w\|_{s,-b}^2.
\end{align*}
For $G_1(\tau)$, we have the following result:
\begin{equation}
\label{G1-est}
\Big\|
(1+|\tau|)^\frac{2(s+j-\ell)}{m}
G_1(\tau)
\Big\|_{L_\tau^\infty}
\le
c_{s,b},
\qquad
-1\le s\le \frac12,
\quad
0\le b<\frac12,
\end{equation}
where $c_{s,b}$ is a constant depending on $s$ and $b$.
We shall prove estimate \eqref{G1-est} later. Now
using it we get the desired estimate  \eqref{forced-ivp-te} for term \eqref{Tfg2-term-recall} in the case $s\le\frac12$.

\noindent
{\bf Case $s\not\in[-1,\frac12]$.} For this case, multiplying 
and dividing $|\widehat{w}(\xi,\tau)|$ by 
$
\frac{(1+|\tau|)^{s/m}}{(1+|\tau-\xi^m|)^{b}}
$ 
and using the Cauchy-Schwartz inequality for the integral of $d\xi$, we get
\begin{align*}
\|\psi h\|_{H^\frac{s+j-\ell}{m}_t(\rr)}^2
\le
c_\psi^2
\int_{\rr}
(1+|\tau|)^\frac{2j-2\ell}{m}
G_2(\tau)
\int_{\mathbb{R}}
\frac{(1+|\tau|)^{2s/m}|\widehat{w}(\xi,\tau)|^2}{(1+|\tau-\xi^m|)^{2b}}
d\xi
d\tau,
\end{align*}
where  
$
G_2(\tau)
\doteq
\int
\frac{\xi^{2\ell}}{(1+|\tau-\xi^m|)^{2-2b}}
d\xi.
$
Like in  the case $s\in[-1,\frac12]$, we get
\begin{align*}
\|\psi h\|_{H^\frac{s+1}{m}_t(\rr)}^2
\leq
c_\psi^2
\Big\|
(1+|\tau|)^\frac{2j-2\ell}{m}
G_2(\tau)
\Big\|_{L_\tau^\infty}
\|w\|_{Y^{s,-b}}^2.
\end{align*}
For $G_2(\tau)$, we have the following result:
\begin{equation}
\label{G2-est}
\Big\|
(1+|\tau|)^\frac{2j-2\ell}{m}
G_2(\tau)
\Big\|_{L_\tau^\infty}
\le
c_{s,b},
\qquad
s\in\rr,
\quad
0<b<\frac12.
\end{equation}
Hence, we complete the proof of Theorem \ref{forced-ivp-thm} once we prove estimates  \eqref{G1-est} and \eqref{G2-est}.
The proof of estimate \eqref{G2-est} is similar to the proof of estimate \eqref{G1-est}. Here we prove  only estimate \eqref{G1-est}.

\nin
{\bf Proof of estimate \eqref{G1-est}.} To prove this estimate, we consider the following two cases:

\vskip0.05in
\noindent
$\bullet$ $|\xi|\le 1$
\quad
and
\quad
$\bullet$
$|\xi|> 1$

\vskip0.05in
\nin
{\bf Case $|\xi|\le 1$.} Since 
$
G_1(\tau)
=
\int_{-1}^1
\frac{\xi^{2\ell}}{(1+|\tau-\xi^m|)^{2-2b}(1+|\xi|)^{2s}}
d\xi,
$
we have
$$
(1+|\tau|)^\frac{2(s+j-\ell)}{m}
G_1(\tau)
\lesssim
\frac{(1+|\tau|)^{\frac{2(s+j-\ell)}{m}}}{(1+|\tau|)^{2-2b}}
=
(1+|\tau|)^{\frac{2s-2\ell+2j-2m+2mb}{m}},
$$
which is bounded since 
$\frac1m\cdot(2s-2\ell+2j-2m+2mb)\le 0$
when $s\le \frac12$, $\ell\ge0$ and $0\le b<\frac12$.

\vskip.05in
\nin
{\bf Case $|\xi|> 1$.} Since  $(1+|\xi|)^{2s} \simeq  |\xi|^{2s}$,
after making the change of variables $\xi_1=\xi^m$, we get
$
%\displaystyle
G_1(\tau)
\lesssim
\int_{|\xi_1|>1}
\frac{1}{(1+|\tau-\xi_1|)^{2-2b}|\xi_1|^{\frac{2s+2j-2\ell}{m}}}
d\xi_1.
$
Next, we consider the following two subcases:

\vskip.05in
\nin
$\bullet$ $|\xi_1|\le \frac12|\tau|$
\quad
and
\quad
$\bullet$ $|\xi_1|\ge \frac12|\tau|$

\vskip0.05in
\nin
{\bf Subcase $|\xi_1|\le \frac12|\tau|$.} Then, we have $(1+|\tau-\xi_1|)^{2-2b}\simeq (1+|\tau|)^{2-2b}$ and $|\tau|\ge 2$. Thus,
\begin{align*}
G_1(\tau)
\simeq&
\int_{|\xi_1|=1}^{\frac12|\tau|}\frac{1}{(1+|\tau-\xi_1|)^{2-2b}|\xi_1|^{\frac{2s+2j-2\ell}{m}}}
d\xi_1
\leq
2
(1+|\tau|)^{2b-2}
\int_{\xi_1=1}^{\frac12|\tau|}
\xi_1^{-\frac{2s+2j-2\ell}{m}}
d\xi_1,
\end{align*}
which implies that  
$
(1+|\tau|)^\frac{2(s+j-\ell)}{m}
G_1(\tau)
$
is bounded  if $s\le 1/2$ and $b<1/2$.

\vskip.05in
\nin
{\bf Subcase $|\xi_1|\ge \frac12|\tau|$.} Since $|\xi_1|\gtrsim 1+ |\tau|$, 
using $s\ge -1$, we get $\frac{2s+2j-2\ell}{m}$, which implies that 
$$
G_1(\tau)
\lesssim
(1+|\tau|)^{-\frac{2s+2j-2\ell}{m}}
\int_{|\xi_1|=\frac12|\tau|}^{\infty}\frac{1}{(1+|\tau-\xi_1|)^{2-2b}}
d\xi_1
\lesssim
(1+|\tau|)^{-\frac{2s+2j-2\ell}{m}}
\int_{0}^{\infty}\frac{1}{(1+x)^{2-2b}}
dx,
$$
where in the last step we make the change of variables  $x=\xi_1-\tau$.
Therefore, we get
$
(1+|\tau|)^\frac{2(s+j-\ell)}{m}
G_1(\tau)
\lesssim
\int_{0}^{\infty}
\frac{1}{(1+x)^{2-2b}}
dx$, which is bounded if $b<\frac12$.
\,\,
$\Box$

\nin
{\bf Estimates for pure ibvp  B$_2$.}
By the time estimate \eqref{forced-ivp-te} we have 
$-W_\ell(t)\in H_t^{\mu_\ell}(0,T)$.
 Thus,  the solution of problem  B$_2$  is like  that of problem A$_2$ and   is 
 estimated by using Theorem \ref{pure-ibvp-thm}.

\nin
{\bf Proof of Theorem \ref{forced-linear-kdvm-thm}.}
Now using the results above we can estimate the solution of
the forced linear  ibvp.  
For $x\ge 0$ and $0\le t\le T<\frac 12$ we have
\begin{align*}
S[u_0,g_0,\dots,g_{j-1};f]
=&
\psi(t)S[U_0;0]
+
S\big[0, G_0,\dots,G_{j-1}; 0\big]
\\
+&
\psi(t)S\big[0; w\big]
+
S\big[0, -W_0,\dots,-W_{j-1}; 0\big].
\end{align*}
This together with estimates \eqref{homo-ivp-te}, \eqref{homo-ivp-se-bourgain} and \eqref{pure-ibvp-thm:se-bourgain}-\eqref{forced-ivp-te}
gives the desired result \eqref{forced-linear-kdvm-est}.
\,\,
$\square$

%\newpage
%
%
%%%%%%%%%%%%%%%%%%%%%%%%%%%%%%%%
%
%
%
%     Proof of Bilinear Estimate for s >=0
%
%
%
%
%%%%%%%%%%%%%%%%%%%%%%%%%%%%%%%%%
%
%
%
\section{ Proof of Bilinear Estimate in modified Bourgain spaces
$X^{s,b,\alpha}$}
\label{sec:bilinear-estimate}
In this section, we prove the bilinear estimate in Bourgain spaces $X^{s,b,\alpha}$.
Following \cite{b1993-kdv} and \cite{kpv1996} we
begin the proof of the bilinear estimate \eqref{bilinear-est}
by first providing an equivalent $L^2$ formulation.
For this, using the fact that 
$a^2+b^2\simeq (|a|+|b|)^2$, we get
\begin{equation}
\label{Bourgain-like-norms-modified}
\|h\|_{s,b,\alpha}^2
\simeq
\int_{\rr}
\int_{\rr}
\left[
(1+|\xi|)^{s}
(1+|\tau-\xi^m|)^{b}
+
\chi_{|\xi|< 1}(1+|\tau|)^{\alpha}
\right]^2
|\widehat{h}(\xi,\tau)|^2
d\xi
d\tau,
\end{equation}
and if for a function $h$ 
we use the Bourgain type combination
\begin{equation}
\label{eq:c_h}
c_h(\xi,\tau)
\doteq
\left[
(1+|\xi|)^{s}
(1+|\tau-\xi^m|)^{b'}
+
\chi_{|\xi|< 1}(1+|\tau|)^{\alpha'}
\right]
|\widehat{h}(\xi,\tau)|,
\end{equation}
then  the modified   Bourgain norm   
of $h$ is equivalent to
the $L^2$ norm of $c_h$, that is
\begin{equation}
\label{eq:c_useful}
\|h\|_{s,b',\alpha'}=\|c_h(\xi,\tau)\|_{L^2_{\xi} L^2_{\tau}}.
\end{equation}
Next, we form the $\|\cdot\|_{s,-b,\alpha-1}$-norm of 
\begin{equation}
\label{w-fg-def}
w_{fg}\doteq\frac 12 \p_x[f\cdot g].
\end{equation}  
Using the definition of convolution and  the relation \eqref{eq:c_h}
we have
\begin{align*}
&|\widehat{w}_{fg}(\xi,\tau)|
\simeq 
\left| \xi
\iint_{\rr^2}
 \widehat{f}(\xi-\xi_1,\tau-\tau_1)\widehat{g}(\xi_1,\tau_1)d\xi_1 d\tau_1 \right|
\\
\leq&
|\xi  |
\iint_{\rr^2}
\frac{c_f(\xi-\xi_1,\tau-\tau_1) }
{
(1+|\xi-\xi_1|)^{s}
(1+|\tau-\tau_1-(\xi-\xi_1)^m|)^{b'}
+
\chi_{|\xi-\xi_1|\leq 1}(1+|\tau-\tau_1|)^{\alpha'}
}
\\
&
\times
\frac{c_g(\xi_1,\tau_1)}
{
(1+|\xi_1|)^{s}
(1+|\tau_1-\xi_1^m|)^{b'}
+
\chi_{|\xi_1|\leq 1}(1+|\tau_1|)^{\alpha'}
}
d\xi_1 d\tau_1.
\end{align*}
Then, the $\|\cdot\|_{s,-b,\alpha-1}$-norm reads as follows
\begin{align*}
\|w_{fg}\|_{s,-b,\alpha-1}^2
=&
\iint_{\rr^2}
\left\{
\iint_{\rr^2}
|\xi  |
\cdot
\left[
\frac{
(1+|\xi|)^s
}{
(1+|\tau-\xi^m|)^{b}
}
+
\frac{
\chi_{|\xi|< 1}
}
{
(1+|\tau|)^{1-\alpha}
}
\right]
\right.
\\
&
\times
\frac{c_f(\xi-\xi_1,\tau-\tau_1) }
{
(1+|\xi-\xi_1|)^{s}
(1+|\tau-\tau_1-(\xi-\xi_1)^m|)^{b'}
+
\chi_{|\xi-\xi_1|\leq 1}(1+|\tau-\tau_1|)^{\alpha'}
}
\\
&
\left.
\times
\frac{c_g(\xi_1,\tau_1)}
{
(1+|\xi_1|)^{s}
(1+|\tau_1-\xi_1^m|)^{b'}
+
\chi_{|\xi_1|\leq 1}(1+|\tau_1|)^{\alpha'}
}
d\xi_1 d\tau_1
\right\}^2
d\xi d\tau.
\end{align*}
Collecting all multipliers together to form the 
important quantity
\begin{subequations}
\label{Q-def}
\begin{align}
Q(\xi,\xi_1,\tau,\tau_1)
\doteq&
|\xi|
\nn
\\
\label{Q-def-term-a}
\times&\left[
\frac{
(1+|\xi|)^s
}{
(1+|\tau-\xi^m|)^{b}
}
+
\frac{
\chi_{|\xi|< 1}
}
{
(1+|\tau|)^{1-\alpha}
}
\right]
\\
\label{Q-def-term-b}
\times&\frac{1}{
(1+|\xi_1|)^{s}
(1+|\tau_1-\xi_1^m|)^{b'}
+
\chi_{|\xi_1|\leq 1}(1+|\tau_1|)^{\alpha'}
}
\\
\label{Q-def-term-c}
\times&
\frac{1}{
(1+|\xi-\xi_1|)^{s}
(1+|\tau-\tau_1-(\xi-\xi_1)^m|)^{b'}
+
\chi_{|\xi-\xi_1|\leq 1}(1+|\tau-\tau_1|)^{\alpha'}
}.
\end{align}
\end{subequations}
we see that to prove bilinear estimate  \eqref{bilinear-est} it suffices to show the $L^2$ inequality
\begin{align}
\label{bilinear-est-L2-form}
 \Bigg\| 
 \iint_{\rr^2}
 Q(\xi,\xi_1,\tau,\tau_1)
 c_f(\xi-\xi_1,\tau-\tau_1) c_g(\xi_1,\tau_1) 
 d\xi_1 d\tau_1 \Bigg\|_{L^2_\xi L^2_{\tau}}
\lesssim
\left\|c_f\right\|_{L^2_{\xi} L^2_{\tau}} \left\|c_g\right\|_{L^2_{\xi} L^2_{\tau}}.
\end{align}
Next, we   bound the term \eqref{Q-def-term-a} from above. Since
\begin{align}
\label{relation-alpha-b}
1-\alpha
\geq
b,
\end{align}
the second term that has numerator $\chi_{|\xi|< 1}$
(and so $\xi$ is bounded) is absorbed 
by the first term. Thus,
\begin{equation}
\label{Q-def-term-a-est}
\eqref{Q-def-term-a}
=
\frac{
(1+|\xi|)^s
}{
(1+|\tau-\xi^m|)^{b}
}
+
\frac{
\chi_{|\xi|< 1}
}
{
(1+|\tau|)^{1-\alpha}
}
\lesssim
\frac{
(1+|\xi|)^s
}{
(1+|\tau-\xi^m|)^{b}
}.
\end{equation}
Combining  estimate \eqref{Q-def-term-a-est} with  estimate \eqref{Q-def} we get the following form of $Q$,
which still involves $\alpha$ terms
\begin{subequations}
\label{Q-def-no1}
\begin{align}
\label{Q-def-no1-term-a}
Q(\xi,\xi_1,\tau,\tau_1)
\lesssim&
|\xi|
\times
\frac{
(1+|\xi|)^s
}{
(1+|\tau-\xi^m|)^{b}
}
\\
\label{Q-def-no1-term-b}
\times&\frac{1}{
(1+|\xi_1|)^{s}
(1+|\tau_1-\xi_1^m|)^{b'}
+
\chi_{|\xi_1|\leq 1}(1+|\tau_1|)^{\alpha'}
}
\\
\label{Q-def-no1-term-c}
\times&
\frac{1}{
(1+|\xi-\xi_1|)^{s}
(1+|\tau-\tau_1-(\xi-\xi_1)^m|)^{b'}
+
\chi_{|\xi-\xi_1|\leq 1}(1+|\tau-\tau_1|)^{\alpha'}
}.
\end{align}
\end{subequations}
So it suffices to prove the bilinear estimates by 
replacing $Q$ with the  right-hand side of the above inequality.
Note that all of the above is valid for any $s$. To make further reduction for $Q$, we need to consider the following two cases:

\vskip0.05in
\nin
$\bullet$ $s\ge 0$
\quad
\text{and}
\quad
$\bullet$ $-j+\frac14<s<0$

\vskip0.05in
\nin
{\bf Case $s\ge 0$.}
Collecting all the factors with $s$ power 
 we rewrite 
the above bound for $Q$ as
\begin{subequations}
\label{Q-modify-no1}
\begin{align}
\label{Q-modify-no1-term-a}
&Q(\xi,\xi_1,\tau,\tau_1)
\lesssim
\frac{|\xi|}{(1+|\tau-\xi^m|)^{b}}
\\
\label{Q-modify-no1-term-b}
\times&
\frac{(1+|\xi|)^s}{
(1+|\xi_1|)^{s}
(1+|\xi-\xi_1|)^{s}
}
\\
\label{Q-modify-no1-term-c}
\times&
\frac{1}{
(1+|\tau_1-\xi_1^m|)^{b'}
+
\chi_{|\xi_1|\leq 1}(1+|\tau_1|)^{\alpha'}
(1+|\xi_1|)^{-s}
}
\\
\label{Q-modify-no1-term-d}
\times&
\frac{1}
{
(1+|\tau-\tau_1-(\xi-\xi_1)^m|)^{b'}
+
\chi_{|\xi-\xi_1|\leq 1}(1+|\tau-\tau_1|)^{\alpha'}
(1+|\xi-\xi_1|)^{-s}
}.
\end{align}
\end{subequations}
Since $s\ge 0$, like on the line, we have the estimate
\begin{equation}
\label{triang-ineq-s}
(1+|\xi|)^s
\le
(1+|\xi-\xi_1|)^s
(1+|\xi_1|)^s
\iff
\frac{(1+|\xi|)^s}{(1+|\xi-\xi_1|)^s(1+|\xi_1|)^s}
\lesssim
 1,
\end{equation}
which helps us remove  term
\eqref{Q-modify-no1-term-b}.
Also, since $|\xi_1|$ and  $|\xi-\xi_1|$ are bounded
we can remove $(1+|\xi_1|)^{-s}$ and 
$(1+|\xi-\xi_1|)^{-s}$. 
Thus, for any $s\geq 0$, 
we have $Q(\xi,\xi_1,\tau,\tau_1)\leq Q_0(\xi,\xi_1,\tau,\tau_1)$, where
\begin{align}
\label{Q-def-s=0}
Q_0(\xi,\xi_1,\tau,\tau_1)
\doteq&
\frac{|\xi|}{(1+|\tau-\xi^m|)^{b}}
\frac{1}{
(1+|\tau_1-\xi_1^m|)^{b'}
+
\chi_{|\xi_1|\leq 1}(1+|\tau_1|)^{\alpha'}
}
\notag
\\
\times&
\frac{1}
{
(1+|\tau-\tau_1-(\xi-\xi_1)^m|)^{b'}
+
\chi_{|\xi-\xi_1|\leq 1}(1+|\tau-\tau_1|)^{\alpha'}
}.
\end{align}
Thus, for $s\geq 0$ to prove
our bilinear estimate \eqref{bilinear-est-L2-form}, it suffices to prove the following simpler one
\begin{align}
\label{bilinear-est-s-0}
\Bigg\| 
 \iint_{\rr^2}
Q_0(\xi,\xi_1,\tau,\tau_1)
c_f(\xi-\xi_1,\tau-\tau_1) c_g(\xi_1,\tau_1) 
d\xi_1 d\tau_1 \Bigg\|_{L^2_\xi L^2_{\tau}}
\lesssim
\left\|c_f\right\|_{L^2_{\xi} L^2_{\tau}} \left\|c_g\right\|_{L^2_{\xi} L^2_{\tau}},
\end{align}
which corresponds to proving the bilinear estimate when  $s=0$.
Moreover, by symmetry (in convolution writing), we may assume that 
\begin{equation}
\label{conv-symmetry-xi}
|\xi-\xi_1|
\leq
|\xi_1|.
\end{equation}
Then, we have\,
$
1=\chi_{|\xi_1|> 1}
\cdot
\chi_{|\xi-\xi_1|> 1}
+
\chi_{|\xi_1|\leq 1}
\cdot
\chi_{|\xi-\xi_1|\leq 1}
+
\chi_{|\xi_1|> 1}
\cdot
\chi_{|\xi-\xi_1|\leq 1}.
$
Therefore, we can rewrite $Q_0(\xi,\tau,\xi_1,\tau_1)$ as:
$
Q_0(\xi,\tau,\xi_1,\tau_1)
=
Q_1(\xi,\tau,\xi_1,\tau_1)
+
Q_2(\xi,\tau,\xi_1,\tau_1)
+
Q_3(\xi,\tau,\xi_1,\tau_1),
$
where
\begin{align}
\label{def-Q1}
Q_1(\xi,\tau,\xi_1,\tau_1)
\doteq&
\chi_{|\xi_1|> 1}
\cdot
\chi_{|\xi-\xi_1|> 1}
\cdot
Q_0(\xi,\tau,\xi_1,\tau_1)
\\
=&
\frac{|\xi|}{(1+|\tau-\xi^m|)^{b}}
\frac{
\chi_{|\xi_1|> 1}
}{
(1+|\tau_1-\xi_1^m|)^{b'}
}
\frac{
\chi_{|\xi-\xi_1|> 1}
}
{
(1+|\tau-\tau_1-(\xi-\xi_1)^m|)^{b'}
},
\nn
\end{align}
\begin{align}
\label{def-Q2}
&Q_2(\xi,\tau,\xi_1,\tau_1)
\doteq
\chi_{|\xi_1|\leq 1}
\cdot
\chi_{|\xi-\xi_1|\leq 1}
\cdot
Q_0(\xi,\tau,\xi_1,\tau_1)
\\
=&
\frac{|\xi|}{(1+|\tau-\xi^m|)^{b}}
\frac{
\chi_{|\xi_1|\leq 1}
}{
(1+|\tau_1-\xi_1^m|)^{b'}
+
(1+|\tau_1|)^{\alpha}
}
\frac{
\chi_{|\xi-\xi_1|\leq 1}
}
{
(1+|\tau-\tau_1-(\xi-\xi_1)^m|)^{b'}
+
(1+|\tau-\tau_1|)^{\alpha'}
},
\notag
\end{align}
and
\begin{align}
\label{def-Q3}
&Q_3(\xi,\tau,\xi_1,\tau_1)
\doteq
\chi_{|\xi_1|> 1}
\cdot
\chi_{|\xi-\xi_1|\leq 1}
\cdot
Q_0(\xi,\tau,\xi_1,\tau_1)
\\
=&
\frac{|\xi|}{(1+|\tau-\xi^m|)^{b}}
\frac{
\chi_{|\xi_1|> 1}
}{
(1+|\tau_1-\xi_1^m|)^{b'}
}
\frac{
\chi_{|\xi-\xi_1|\leq 1}
}
{
(1+|\tau-\tau_1-(\xi-\xi_1)^m|)^{b'}
+
(1+|\tau-\tau_1|)^{\alpha'}
}.
\notag
\end{align}
Note that $Q_1$ is like the regular Bourgain norm, because we do not have any term related to $\alpha$.

\nin
To prove the bilinear estimate \eqref{bilinear-est-s-0}, it suffices to prove that  for $\ell=1,2,3$
\begin{align}
\label{bilinear-est-s-l}
\Bigg\| 
 \iint_{\rr^2}
Q_\ell(\xi,\xi_1,\tau,\tau_1)
c_f(\xi-\xi_1,\tau-\tau_1) c_g(\xi_1,\tau_1) 
d\xi_1 d\tau_1 \Bigg\|_{L^2_\xi L^2_{\tau}}
\hskip-0.2cm
\lesssim
%\hskip-0.1cm
\left\|c_f\right\|_{L^2_{\xi} L^2_{\tau}} \left\|c_g\right\|_{L^2_{\xi} L^2_{\tau}}.
\end{align}
%
%
%%%%%%%%%%%%%%%%%%%%%%
%
%   Estimation when the multiplier is $Q_1$
%
%%%%%%%%%%%%%%%%%%%%%%
%
%
{\bf Estimation when the multiplier is $Q_1$.}
This is similar to the  whole line case (see \cite{fhy2020}).
We do it by considering the following two possibilities:

\vskip0.05in
\nin
$\bullet$ $|\xi|\le 1$
\quad
and 
\quad
$\bullet$ $|\xi|> 1$

\vskip0.05in
\nin
{\bf Case  $|\xi|\le 1$.}  In this case, we need to prove \eqref{bilinear-est-s-l} with  $Q_1$ replaced by $\chi_{|\xi|\le 1}Q_1$, where $\chi_{|\xi|\le 1}$ is the characteristic function of the region
$\{
(\xi,\xi_1,\tau,\tau_1):
|\xi|
\le
1
\}$.
This is done by applying duality 
and the Cauchy-Schwarz inequality  first in $(\xi_1, \tau_1)$ 
and then in $(\xi, \tau)$. Doing so, and after some manipulations,
the desired estimate takes the form
\begin{align*}
%\label{bilinear-est-s-0-c1}
\Big\| 
\iint_{\rr^2} 
(\chi_{|\xi|\le 1} Q_1)(\xi,\xi_1,\tau,\tau_1)
c_f(\xi-\xi_1,\tau-\tau_1) c_g(\xi_1,\tau_1) 
d\xi_1 d\tau_1 
\Big\|_{L^2_{\xi,\tau}}
\le
\|
\Theta_1
\|_{L^{\infty}_{\xi_1,\tau_1}}^{1/2}
\|c_f\|_{L^2_{\xi,\tau}}
\|c_g\|_{L^2_{\xi,\tau}},
\end{align*}
where  $\Theta_1$ is  as in the following lemma, which provides its estimate
\begin{lemma}
\label{caseI-1-kdvm} 
If $0<b<\frac12$, and $0<b'< \frac 12$, then there exists $c>0$ such that for $\xi_1,\tau_1\in\rr$
\begin{align*}
\Theta_1(\xi_1, \tau_1)
\doteq
\frac{1}{(1+|\tau_1-\xi_1^{m}|)^{2b'}}
\int_{|\xi|\le 1}
\int_{\rr}
\frac{\xi^2\,\,
d\tau d\xi
}{(1+|\tau-\tau_1-(\xi-\xi_1)^{m}|)^{2b'}(1+|\tau-\xi^{m}|)^{2b}}
\,
\lesssim 1.
\end{align*}
\end{lemma}
\nin
The proof  is similar to the proof of Lemma 7.1 in \cite{fhy2020} and we omit it.

\vskip0.05in
\nin
{\bf Case $|\xi|>1$.}
In this case, by symmetry (in convolution writing), we may assume that 
\begin{equation}
\label{conv-symmetry-tau}
|\tau-\tau_1-(\xi-\xi_1)^m|
\leq
|\tau_1-\xi_1^m|.
\end{equation}
Therefore, following \cite{b1993-kdv} and \cite{kpv1996},
to prove the bilinear estimate above we consider the following two microlocalizations:

\noindent
{\bf Microlocalization I. $|\tau_1-\xi_1^m|\leq|\tau-\xi^m|$}. In this case we define the domain $B_I$ to be
\begin{align}
\label{def-BI-domain}
B_I
\doteq
\big\{
(\xi,\tau,\xi_1,\tau_1) \in{\mathbb{R}}^4: 
&
|\tau-\tau_1-(\xi-\xi)^m|\leq|\tau_1-\xi_1^m|\leq|\tau-\xi^m|,
\\
&|\xi_1|>1,
|\xi-\xi_1|>1,
|\xi|>1
\big\}.
\nn
\end{align}
\nin
{\bf Microlocalization II. $|\tau-\xi^m|\leq|\tau_1-\xi_1^m|$}. 
In this case we define the domain $B_{II}$ to be
\begin{align}
\label{def-BII-domain}
B_{II}
\doteq
\big\{(\xi,\tau,\xi_1,\tau_1)\in{\mathbb{R}}^4: |\tau-\tau_1-(\xi-\xi)^m|\leq|\tau_1-\xi_1^m|,
\\
|\tau-\xi^m|\leq|\tau_1-\xi_1^m|,
|\xi_1|>1,
|\xi-\xi_1|>1,
|\xi|>1
\big\}.
\nn
\end{align}
{\bf Proof of bilinear estimate in Microlocalization I:} 
In this case,   $Q_1$ is replaced by the set $\chi_{B_I}Q_1$. Like before, using the Cauchy-Schwarz inequality with respect to 
$(\xi_1,\tau_1)$ and taking the supremum in  $(\xi,\tau)$ we
arrive at
\begin{align*}
%\label{bilinear-est-neg-BI-2}
\Big\| 
 \iint_{\rr^2}
 (\chi_{B_I} Q_1)(\xi,\xi_1,\tau,\tau_1)
 c_f(\xi-\xi_1,\tau-\tau_1) c_g(\xi_1,\tau_1) 
 d\xi_1 d\tau_1 \Big\|_{L^2_{\xi, \tau}}
\le
\|
\Theta_2
\|_{L^{\infty}_{\xi, \tau}}^{1/2}
\|c_f\|_{L^2_{\xi, \tau}} \|c_g\|_{L^2_{\xi, \tau}}.
\end{align*}
Thus, to prove our bilinear estimate in microlocalization I, it suffices to show 
the following result.

\begin{lemma} 
\label{s-0-case-1-lemma}
If $\frac{6+3m}{12m}\leq b'\leq b<1/2$, then  for $|\xi|>1$, and $\tau\in\rr$ we have
\begin{align}
\label{s-0-case1-theta-1-est-no1}
\Theta_2(\xi,\tau)
\doteq
\frac{\xi^2}{(1+|\tau-\xi^m|)^{2b}}
\iint_{\rr^2}
\frac{ \chi_{B_I}(\xi,\tau,\xi_1,\tau_1)\,\, d\tau_1 d\xi_1 
}{(1+|\tau-\tau_1-(\xi-\xi_1)^m|)^{2b'}(1+|\tau_1-\xi_1^m|)^{2b'}}
\lesssim 1.
\end{align}
\end{lemma}

\noindent
{\bf Proof of bilinear estimate in Microlocalization II:}
Using  duality and applying the Cauchy-Schwarz inequality twice, 
first in $(\xi_1, \tau_1)$ and then in $(\xi, \tau)$,
we get
\begin{align*}
%\label{bilinear-est-L2-form-leq0-3}
\Big\| 
\iint_{\rr^2}
(\chi_{B_{II}} Q_1)(\xi,\xi_1,\tau,\tau_1)
c_f(\xi-\xi_1,\tau-\tau_1) c_g(\xi_1,\tau_1) 
d\xi_1 d\tau_1 \Big\|_{L^2_{\xi,\tau}}
\le
\|
\Theta_{3}
\|_{L^{\infty}_{\xi_1,\tau_1}}^{1/2}
\|c_f\|_{L^2_{\xi,\tau}}\|c_g\|_{L^2_{\xi,\tau}}.
\end{align*}
Thus, to prove our bilinear estimate in microlocalization II, it suffices to show 
the following result.
\begin{lemma}
\label{s-0-case-2-lemma} 
If $\max\{\frac{4+3(m-1)}{12(m-1)},\frac{6+3m}{12m}\}\leq b'\leq b<1/2$, then  for $\xi_1,\tau_1\in\rr$ we have
\begin{align}
\label{s-0-case2-theta-2-est-no1}
\Theta_3(\xi_1, \tau_1)
\doteq
\frac{1}{(1+|\tau_1-\xi_1^m|)^{2b'}}
\iint_{\rr^2}
\frac{
\chi_{B_{II}}(\xi,\tau,\xi_1,\tau_1)
\,\,
\xi^2\,
\, d\tau d\xi}{(1+|\tau-\tau_1-(\xi-\xi_1)^m|)^{2b'}(1+|\tau-\xi^m|)^{2b}}
\lesssim 1.
\end{align}
\end{lemma}

\vskip.05in
\nin
Next, we shall prove the Lemmas \ref{s-0-case-1-lemma} and \ref{s-0-case-2-lemma}. We begin with the first one.

\nin
{\bf Proof of Lemma \ref{s-0-case-1-lemma}.}  To estimate the quantity $\Theta_{2}(\xi,\tau)$ (and  $\Theta_{3}(\xi_1,\tau_1)$ later),
we shall need the following calculus  estimates,
whose proof  can be found in \cite{kpv1996}, \cite{h2006} and \cite{fhy2020}.

\begin{lemma}
\label{lem:calc_ineq}  
If $1>\ell>1/2, l'>1/2$ then
\begin{equation}
\label{eq:calc_1}
\int_{\mathbb{R}}\frac{dx}{(1+|x-a|)^{2\ell}(1+|x-c|)^{2\ell}}\lesssim\frac{1}{(1+|a-c|)^{2\ell}},
\end{equation}

\begin{equation}
\label{eq:calc_2}
\int_{\mathbb{R}}\frac{dx}{(1+|x|)^{2\ell}|{a-x}|^{\frac{1}{2}}}\lesssim\frac{1}{(1+|a|)^{\frac{1}{2}}},
\end{equation}

\begin{equation}
\label{eq:calc_3}
\int_{\mathbb{R}}\frac{dx}{(1+|x-a|)^{2(1-\ell)}(1+|x-c|)^{2\ell'}}\lesssim\frac{1}{(1+|a-c|)^{2(1-\ell)}},
\end{equation}
\begin{equation}
\label{eq:calc_4}
\int_{|x|\leq c}\frac{dx}{(1+|x|)^{2(1-l)}|\sqrt{a-x}|}\lesssim\frac{(1+c)^{2(\ell-1/2)}}{(1+|a|)^{1/2}},
\end{equation}
and
\begin{equation}
\label{eq:calc_1a}
\int_{\mathbb{R}}\frac{dx}{(1+|x-a|)^{2\ell}(1+|x-c|)^{2\ell'}}\lesssim\frac{1}{(1+|a-c|)^{2\min\{\ell',\ell\}}}.
\end{equation}
In addition, if $\frac14<\ell'\leq \ell<\frac12$, then
\begin{equation}
\label{eq:calc_5}
\int_{\mathbb{R}}\frac{dx}{(1+|x-a|)^{2\ell}(1+|x-c|)^{2\ell'}}\lesssim\frac{1}{(1+|a-c|)^{2\ell+2\ell'-1}}.
\end{equation}
\end{lemma}
\nin
Also, in this case using the triangle inequality we have 
$$
|\tau-(\xi-\xi_1)^m-\xi_1^m|
=
|\tau-\tau_1-(\xi-\xi_1)^m+\tau_1-\xi_1^m|
\le
|\tau-\tau_1-(\xi-\xi_1)^m|+|\tau_1-\xi_1^m|
\le
2|\tau-\xi^m|.
$$
Furthermore, integrating  with respect to $\tau_1$ and 
applying estimate \eqref{eq:calc_5}  with $\ell=\ell'=b'$, $a=\tau-(\xi-\xi_1)^m$ and $c=\xi_1^m$, we get
\begin{equation}
\label{s-0-case1-theta-1-est-no3}
\Theta_2(\xi,\tau)
\lesssim
\frac{\xi^2}{(1+|\tau-\xi^m|)^{2b}}
\int_{\rr}
\frac{ 
 d\xi_1 
}{(1+|\tau-(\xi-\xi_1)^m-\xi_1^m|)^{4b'-1}}
=
\frac{\xi^2}{(1+|\tau-\xi^m|)^{2b}}
I(\xi,\tau),
\end{equation}
where $I$ is defined in the lemma below, where it is also estimated.

\begin{lemma}
\label{l-lem-s-0}
Let  $m=2j+1\ge 3$ and $\frac14<b'<1/2$.
If $|\xi^m-(\xi-\xi_1)^m-\xi_1^m|\lesssim |\tau-\xi^m|$, then for all $|\xi|>1$ and $\tau\in\rr$ we have 
\begin{align}
\label{I-est-s-0}
I(\xi, \tau)
\doteq
\int_{\rr}
\frac{
d\xi_1
}{(1+|\tau-(\xi-\xi_1)^m-\xi_1^m|)^{4b'-1}}
\lesssim
\frac{|\xi|^{-\frac12(m-2)}(1+|\tau-\xi^m|)^{2-4b'}}{(1+|\tau-2^{1-m}\xi^m|)^{\frac12}}.
\end{align}
\end{lemma}
The  proof of lemma \ref{l-lem-s-0}  is similar to the proof of Lemma 7.3 in \cite{fhy2020}. Combining  \eqref{I-est-s-0} with \eqref{s-0-case1-theta-1-est-no3}, we get the desired  estimate \eqref{s-0-case1-theta-1-est-no1} for  $\Theta_2$
 and this completes the proof of Lemma \ref{s-0-case-1-lemma}.
\,\,
$\square$

\nin
{\bf Proof of Lemma \ref{s-0-case-2-lemma}.} 
We recall that
$$
\Theta_3(\xi_1, \tau_1)
=
\frac{1}{(1+|\tau_1-\xi_1^{m}|)^{2b'}}
\iint_{\rr^2}
\frac{
\chi_{B_{II}}(\xi,\tau,\xi_1,\tau_1)
\,\,
\xi^2\,
\, d\tau d\xi}{(1+|\tau-\tau_1-(\xi-\xi_1)^{m}|)^{2b'}(1+|\tau-\xi^{m}|)^{2b}}.
$$
As before, using estimate \eqref{eq:calc_5} with $x=\tau$, $\ell=b$, $\ell'=b'$, $a=\xi^m$ and $c=\tau_1-(\xi_1-\xi)^m$, we get
\begin{align}
\label{s-0-case2-theta-2-est-no3}
\Theta_3(\xi_1,\tau_1)
\lesssim&
\frac{1}{(1+|\tau_1-\xi_1^m|)^{2b'}}
\int_{|\xi|>1}
\frac{ 
 \xi^2 \,\,d\xi 
}{(1+|\tau_1-(\xi_1-\xi)^m-\xi^m|)^{2b+2b'-1}}
\notag
\\
\le&
\frac{1}{(1+|\tau_1-\xi_1^m|)^{2b'}}
\int_{|\xi|>1}
\frac{ 
 \xi^2\,\,
  d\xi 
}{(1+|\tau_1-\xi_1^m-d_m(\xi_1,\xi)|)^{4b'-1}},
\end{align}
where $d_m(\xi_1,\xi)$ is defined as follows
\begin{align}
\label{dm-def-no1}
d_m(\xi_1,\xi)
\doteq
-\xi_1^m
+\xi^m
-(\xi-\xi_1)^m.
\end{align}
In order to show that $\Theta_3$ is bounded, 
 we need to consider the following two cases:

\vskip0.05in
\noindent
$\bullet$ $|\xi|\leq 10|\xi_1|$
\quad
and 
\quad
$\bullet$ $|\xi|> 10|\xi_1|$

\vskip0.05in
\noindent
{\bf Case $|\xi|\leq 10|\xi_1|$.} Then, 
\begin{align*}
%\label{s-0-case2-theta-3-est-no1}
\Theta_3(\xi_1,\tau_1)
\lesssim&
\frac{\xi_1^2}{(1+|\tau_1-\xi_1^m|)^{2b'}}
\int_{|\xi|>1}
\frac{ 
d\xi 
}{(1+|\tau_1-\xi_1^m-d_m(\xi_1,\xi)|)^{4b'-1}}
\lesssim
\frac{\xi_1^2}{(1+|\tau_1-\xi_1^m|)^{2b'}}
I(\xi_1,\tau_1),
\end{align*}
where $I$ is defined in \eqref{I-est-s-0}. Applying Lemma \ref{l-lem-s-0}, we get
\begin{align}
\label{s-0-case1-theta-3-est-no2}
\Theta_3(\xi_1,\tau_1)
\lesssim&
\frac{|\xi_1|^2}{(1+|\tau_1-\xi_1^m|)^{2b'}}
\frac{|\xi_1|^{-\frac12(m-2)}(1+|\tau_1-\xi_1^m|)^{2-4b'}}{(1+|\tau_1-2^{1-m}\xi_1^m|)^{\frac12}}
\nn
\\
=&
\frac{|\xi_1|^{3-\frac{m}{2}}}{(1+|\tau_1-\xi_1^m|)^{6b'-2}(1+|\tau_1-2^{1-m}\xi_1^m|)^{\frac12}}.
\end{align}
Like  KdVm  on the line (see \cite{fhy2020}), we  consider the following  two subcases:

\vskip0.05in
\nin
$
\bullet
\,\,
|\tau_1-\xi_1^m|\leq \frac12|\xi_1|^m
$
\quad
and
\quad
$
\bullet
\,\,
|\tau_1-\xi_1^m|> \frac12|\xi_1|^m
$

\vskip0.05in
\nin
{\bf Subcase $|\tau_1-\xi_1^m|\leq \frac12|\xi_1|^m$.}
Then, using the triangle inequality, we have 
$$
|\tau_1-2^{1-m}\xi_1^m|
=
|(\tau_1-\xi_1^m)+(1-2^{1-m})\xi_1^m|
\geq
\frac34|\xi_1^m|
-
|\tau_1-\xi_1^m|
\geq
\frac34|\xi_1^m|
-
\frac12|\xi_1^m|
=
\frac14|\xi_1^m|.
$$
Hence, from \eqref{s-0-case1-theta-3-est-no2} we get
$
\Theta_3(\xi_1,\tau_1)
\lesssim
\frac{|\xi_1|^{3-\frac{m}{2}}}{1}
\frac{1}{|\xi_1|^{m/2}}
=
|\xi_1|^{-m+3}
\overset{|\xi_1|>1}
{
\lesssim
}
1.
$

\vskip0.05in
\noindent
{\bf Subcase  $|\tau_1-\xi_1^m|> \frac12|\xi_1|^m$.} Then
$
|\tau_1-\xi_1^m|
\gtrsim
|\xi_1|^m
$
and therefore from \eqref{s-0-case1-theta-3-est-no2}  we get
\begin{align*}
\Theta_3(\xi_1,\tau_1)
\lesssim
\frac{|\xi_1|^{3-\frac{m}{2}}}{|\xi_1^m|^{6b'-2}}
\cdot
\frac{1}{1}
\leq
\frac{|\xi_1|^{3-\frac{m}{2}}}
{
|\xi_1|^{m(6b'-2)}
}
=
\frac{1}{|\xi_1|^{m(6b'-2)-(3-\frac{m}{2})}}.
\end{align*}
Since $|\xi_1|> 1$, the above quantity is bounded  if 
$
m(6b'-2)
-
(3-\frac{m}{2})
\ge
0,
$
which implies that
\begin{align}
\label{b-con-1}
b'
\geq
\frac{6+3m}{12m}.
\end{align} 
This completes the proof of Lemma \ref{s-0-case-2-lemma}  in this case.

\vskip0.05in
\nin
{\bf Case $|\xi|> 10|\xi_1|$.}  Then,  using the triangle inequality  $|\tau_1-\xi_1^m|$ is bounded from below as follows
\begin{align}
\label{est-tau1}
|\tau_1-\xi_1^m|
\geq&
\frac13
\Big[
|\tau-\xi^m|
+
|
(\tau_1-\xi_1^m)
|
+
|\tau -\tau_1 -(\xi-\xi_1)^m
|
\Big]
\notag
\\
\gtrsim&
|\tau-\xi^m
+
(\tau_1-\xi_1^m)
+
[
\tau -\tau_1 -(\xi-\xi_1)^m
]
|
=
|d_m(\xi,\xi_1)|,
\end{align}
which can be bounded by the following result.
\begin{lemma}
\label{lem-dm}
If $m$ is an odd positive integer,  then  there is a positive constant $c_m$ such that
\begin{align}
\label{dm-k}
|d_m(\xi,\xi_1)|
\geq
c_m|\xi|^{m-3}|\xi\xi_1(\xi-\xi_1)|,
\end{align}
\begin{align}
\label{dm-k1}
    \text{and}
    \quad
|d_m(\xi,\xi_1)|
\geq
c_m|\xi_1|^{m-3}|\xi\xi_1(\xi-\xi_1)|.
\end{align}
Also if $|\xi|\geq 1$, $|\xi_1|\geq 1$ and $|\xi-\xi_1|\geq 1$, then 
\begin{equation}
\label{basic-est-2}
|\xi\xi_1(\xi-\xi_1)|
\geq
\frac13
\xi^2
\quad
\text{and}
\quad
|\xi\xi_1(\xi-\xi_1)|
\geq
\frac13\xi_1^2.
\end{equation}
\end{lemma}
The proof of Lemma \ref{lem-dm} can be found in \cite{fhy2020}, now we complete the proof in this case. Since $|\xi|>10|\xi_1|$, we have
$
|\xi-\xi_1|
\geq
|\xi|
-
|\xi_1|
\geq
\frac{9}{10}|\xi|
$
and
$
|\xi-\xi_1|
\leq
|\xi|
+
|\xi_1|
\geq
\frac{11}{10}|\xi|,
$
which gives us that 
\begin{equation}
\label{est-xi-xi1}
|\xi-\xi_1|
\simeq
|\xi|.
\end{equation}
Combining  estimates  \eqref{est-tau1} and  \eqref{dm-k} with \eqref{est-xi-xi1}, we get
$
|\tau_1-\xi_1^m|
\gtrsim
|\xi|^{m-1}|\xi_1|
$
or
$
|\xi|
\leq
(|\tau_1-\xi_1^m\|\xi_1|^{-1})^{\frac{1}{m-1}}.
$
In addition, using estimate \eqref{s-0-case2-theta-2-est-no3} we have
\begin{align}
\label{s-0-case2-theta-3-est-no1b}
\Theta_3(\xi_1,\tau_1)
\lesssim&
\frac{|\xi_1|^{-\frac{2}{m-1}}|\tau_1-\xi_1^m|^{\frac{2}{m-1}}}{(1+|\tau_1-\xi_1^m|)^{2b'}}
\int_{|\xi|>1}
\frac{ 
d\xi 
}{(1+|\tau_1-\xi_1^m-d_m(\xi_1,\xi)|)^{4b'-1}}
\\
\leq&
|\xi_1|^{-\frac{2}{m-1}}(1+|\tau_1-\xi_1^m|)^{\frac{2}{m-1}-2b'}
I(\xi_1,\tau_1),
\nn
\end{align}
where $I$ is defined in \eqref{I-est-s-0}. Applying Lemma \ref{l-lem-s-0}, we get
\begin{align}
\label{s-0-case1-theta-3-est-no2b}
\Theta_3(\xi_1,\tau_1)
\lesssim&
|\xi_1|^{-\frac{2}{m-1}}(1+|\tau_1-\xi_1^m|)^{\frac{2}{m-1}-2b'}
\frac{|\xi_1|^{-\frac12(m-2)}(1+|\tau_1-\xi_1^m|)^{2-4b'}}{(1+|\tau_1-2^{1-m}\xi_1^m|)^{\frac12}}
\notag
\\
=&
|\xi_1|^{-\frac{2}{m-1}-\frac{1}{2}(m-2)}
\frac{(1+|\tau_1-\xi_1^m|)^{\frac{2}{m-1}+2-6b'}}{(1+|\tau_1-2^{1-m}\xi_1^m|)^{\frac12}}.
\end{align}
Now using the triangle inequality, we can bound $|\tau_1-2^{1-m}\xi_1^m|$ from below, that is
\begin{align*}
|\tau_1-2^{1-m}\xi_1^m|
=&
|\tau_1-\xi_1^m+(1-2^{1-m})\xi_1|
\geq
|\tau_1-\xi_1^m|
-
|\xi_1^m|
=
\frac12|\tau_1-\xi_1^m|
+
\left(\frac12|\tau_1-\xi_1^m|
-
|\xi_1^m|
\right)
\\
\overset{\eqref{est-tau1}}
{\geq}&
\frac12|\tau_1-\xi_1^m|
+
\left(\frac12|\xi_1|\xi^{m-1}
-
|\xi_1^m|
\right)
\overset{|\xi|>10|\xi_1|}
{
\geq
}
\frac12|\tau_1-\xi_1^m|.
\end{align*}
Combining the above estimate with \eqref{s-0-case1-theta-3-est-no2b}, we obtain
\begin{align}
\label{s-0-case1-theta-3-est-no3b}
\Theta_3(\xi_1,\tau_1)
\lesssim
|\xi_1|^{-\frac{1}{m-1}-\frac{1}{2}(m-2)}
(1+|\tau_1-\xi_1^m|)^{\frac{2}{m-1}+\frac{3}{2}-6b'}.
\end{align}
Since $|\xi_1|>1$, the above quantity is bounded if and only if
\begin{equation}
\label{b'-cond-2}
b'
\geq
\frac{4+3(m-1)}{12(m-1)}.
\end{equation}
This completes the proof of Lemma \ref{s-0-case-2-lemma}.
\,\,
$\Box$

%
%
%%%%%%%%%%%%%%%%%%%%%%
%
%   Estimation when the multiplier is $Q_2$
%
%%%%%%%%%%%%%%%%%%%%%%
%
%
\vskip0.1in
\nin
{\bf Estimation when the multiplier is $Q_2$.}
In this case, applying Cauchy-Schwarz inequality with respect to 
$\xi_1,\tau_1$, nd taking the super norm over  $(\xi,\tau)$ we get
\begin{align}
\label{bilinear-est-s-Q2}
\Bigg\| 
\iint_{\rr^2}
Q_2(\xi,\xi_1,\tau,\tau_1)
c_f(\xi-\xi_1,\tau-\tau_1) c_g(\xi_1,\tau_1) 
d\xi_1 d\tau_1 \Bigg\|_{L^2_\xi L^2_{\tau}}
\lesssim
\|
\Theta_4
\|_{L^{\infty}_{\xi, \tau}}^{\frac12}
\left\|c_f\right\|_{L^2_{\xi} L^2_{\tau}} \left\|c_g\right\|_{L^2_{\xi} L^2_{\tau}
},
\end{align}
which shows that the proof of the bilinear estimate  follows from the next result.
\begin{lemma} 
\label{s-0-theta-4-lemma}
If $b\ge 0$ and $\alpha'>\frac12$, then for $\xi,\tau\in\rr$, we have
\begin{align}
\label{s-0-theta-4}
\Theta_4(\xi,\tau)
\doteq&
\frac{\xi^2}{(1+|\tau-\xi^m|)^{2b}}
\iint_{\rr^2}
\frac{
\chi_{|\xi_1|\leq 1}
}{
[
(1+|\tau_1-\xi_1^m|)^{b'}
+
(1+|\tau_1|)^{\alpha'}
]^2
}
\notag
\\
&\frac{
\chi_{|\xi-\xi_1|\leq 1}
\quad
d\xi_1
d\tau_1
}
{
[
(1+|\tau-\tau_1-(\xi-\xi_1)^m|)^{b'}
+
(1+|\tau-\tau_1|)^{\alpha'}
]^2
}
\lesssim 1.
\end{align}
\end{lemma}
The proof of this lemma is straightforward.
Using the fact that $|\xi_1|$ and $|\xi|$ are bounded and applying estimate \eqref{eq:calc_1} with $\ell=\alpha'$, $x=\tau_1$, $a=0$ and $c=\tau_1$, we get the desired estimate \eqref{s-0-theta-4}.

%
%
%%%%%%%%%%%%%%%%%%%%%%
%
%
%
%   Estimation when the multiplier is $Q_3$
%
%
%
%%%%%%%%%%%%%%%%%%%%%%
%
%
\nin
{\bf Estimation when the multiplier is $Q_3$.}
To prove the estimate \eqref{bilinear-est-s-l} for $Q_3$, we will consider two possible microlocalizations:
\vskip.05in
\noindent
{\bf Microlocalization III. $|\tau_1-\xi_1^m|\leq|\tau-\xi^m|$}. In this case we define the domain $B_{III}$ to be
\begin{align}
\label{def-BIII-domain}
B_{III}
\doteq
\big\{
(\xi,\tau,\xi_1,\tau_1) \in{\mathbb{R}}^4: 
|\tau_1-\xi_1^m|\leq|\tau-\xi^m|, |\xi_1|>1,
|\xi-\xi_1|\leq1
\big\}.
\end{align}
\noindent
{\bf Microlocalization IV. $|\tau-\xi^m|\leq|\tau_1-\xi_1^m|$}. 
In this case we define the domain $B_{IV}$ to be
\begin{align}
\label{def-BIV-domain}
B_{IV}
\doteq
\big\{(\xi,\tau,\xi_1,\tau_1)\in{\mathbb{R}}^4: 
|\tau-\xi^m|\leq|\tau_1-\xi_1^m|,
|\xi_1|>1,
|\xi-\xi_1|\leq1
\big\}.
\end{align}
\nin
{\bf Proof of bilinear estimate in Microlocalization III.}
As before, using   the Cauchy-Schwarz inequality with respect to $(\xi_1,\tau_1)$  and taking the superemum over  $(\xi,\tau)$ we
arrive at
\begin{align*}
%\label{bilinear-est-BIII-2-Q3}
\Bigg\| 
 \iint_{\rr^2}
 \chi_{B_{III}}
Q_3(\xi,\xi_1,\tau,\tau_1)
 c_f(\xi-\xi_1,\tau-\tau_1) c_g(\xi_1,\tau_1) 
 d\xi d\tau \Bigg\|_{L^2_\xi L^2_{\tau}}
\lesssim
\|
\Theta_{5}
\|_{L^{\infty}_{\xi,\tau}}^{1/2}
\|c_f\|_{L^2_{\xi}L^2_{\tau}}
\|c_g\|_{L^2_{\xi}L^2_{\tau}}.
\end{align*}
Thus, to prove our bilinear estimate in microlocalization III, it suffices to show 
the following result.
\begin{lemma} 
\label{s-0-theta-5-lemma}
If   $\frac13\leq b'\leq b<\frac12<\alpha'$, then for $\xi,\tau\in\rr$ we have
\begin{align}
\label{s-0-theta-5}
\Theta_5(\xi,\tau)
\doteq&
\frac{
\xi^2
}{
(1+|\tau-\xi^m|)^{2b}
}
\iint_{\rr^2}
\frac{
\chi_{|\xi_1|>1}
\chi_{B_{III}}(\xi,\tau,\xi_1,\tau_1)
}{
(1+|\tau_1-\xi_1^m|)^{2b'}
}
\notag
\\
&\frac{
\chi_{|\xi-\xi_1|\leq 1}
\quad
d\xi_1
d\tau_1
}
{
[
(1+|\tau-\tau_1-(\xi-\xi_1)^m|)^{b'}
+
(1+|\tau-\tau_1|)^{\alpha'}
]^2
}
\lesssim 1.
\end{align}
\end{lemma}

\nin
{\bf Proof of bilinear estimate in Microlocalization IV.} As before, using duality and  the Cauchy-Schwarz inequality twice,  first in $(\xi_1, \tau_1)$ and then in $(\xi, \tau)$, we get
\begin{align*}
%\label{bilinear-est-BIV-2-Q3}
\Bigg\| 
\iint_{\rr^2}
\chi_{B_{IV}}
Q_3(\xi,\xi_1,\tau,\tau_1)
c_f(\xi-\xi_1,\tau-\tau_1) c_g(\xi_1,\tau_1) 
d\xi_1 d\tau_1 \Bigg\|_{L^2_\xi L^2_{\tau}}
\lesssim
\|
\Theta_{6}
\|_{L^{\infty}_{\xi_1,\tau_1}}^{1/2}
\|c_f\|_{L^2_{\xi}L^2_{\tau}}\|c_g\|_{L^2_{\xi}L^2_{\tau}}.
\end{align*}
Thus, to prove our bilinear estimate in microlocalization IV, it suffices to show 
the following result.
\begin{lemma} 
\label{s-0-theta-6-lemma}
If   $\frac13\leq b'\leq b<\frac12<\alpha'$, then for $\xi,\tau\in\rr$ we have
\begin{align}
\label{s-0-theta-6}
\Theta_6(\xi_1,\tau_1)
\doteq&
\frac{
\chi_{|\xi_1|>1}
}{
(1+|\tau_1-\xi_1^m|)^{2b'}
}
\iint_{\rr^2}
\frac{
\chi_{B_{IV}}(\xi,\tau,\xi_1,\tau_1)
\,\,
\xi^2
}{
(1+|\tau-\xi^m|)^{2b}
}
\notag
\\
&\frac{
\chi_{|\xi-\xi_1|\leq 1}
\quad
d\xi
d\tau
}
{
[
(1+|\tau-\tau_1-(\xi-\xi_1)^m|)^{b'}
+
(1+|\tau-\tau_1|)^{\alpha'}
]^2
}
\lesssim 1.
\end{align}
\end{lemma}
The proof of Lemma \ref{s-0-theta-6-lemma} is similar to the proof of Lemma \ref{s-0-theta-5-lemma}. So, here we provide only the proof of Lemma \ref{s-0-theta-5-lemma}.

\medno
{\bf Proof of Lemma \ref{s-0-theta-5-lemma}.}
Since 
$
|\tau-\xi^m|
\geq
|\tau_1-\xi_1^m|
$,
and using $\alpha'>\frac12>b'$, $|\xi-\xi_1|\le 1$
we get
\begin{align}
\label{s-0-theta-5-est-no1}
\Theta_5(\xi,\tau)
\lesssim
\frac{
\xi^2
}{
(1+|\tau-\xi^m|)^{4b+2b'-2}
}
\iint_{\rr^2}
\frac{
\chi_{|\xi_1|>1}
}{
(1+|\tau_1-\xi_1^m|)^{2-2b}
}
\frac{
\chi_{|\xi-\xi_1|\leq 1}
}
{
(1+|\tau-\tau_1|)^{2\alpha'}
}
d\xi_1
d\tau_1.
\end{align}
Now, using estimate \eqref{eq:calc_1a} with $\ell'=\alpha'$, $\ell=1-b$, $x=\tau_1$, $a=\tau$ and $c=\xi_1^m$
\begin{equation}
\label{s-0-theta-5-est-no1}
\Theta_5(\xi,\tau)
\lesssim
\frac{
\xi^2
}{
(1+|\tau-\xi^m|)^{4b+2b'-2}
}
\int_\rr
\frac{
\chi_{|\xi_1|>1}
\chi_{|\xi-\xi_1|\leq 1}
}{
(1+|\tau-\xi_1^m|)^{\min\{2\alpha',2-2b\}}
}
d\xi_1.
\end{equation}
For $|\xi|\le 20$, it is obvious that $Q_5\lesssim 1$. For $|\xi|\ge 20$, using the fact that $\xi_1\simeq \xi$ and making the change of variables $\mu=\xi^m$ for the integral of $d\xi$, for $4b+2b'-2>0$ we also get $Q_5\lesssim 1$.
Thus, we complete the proof of Lemma \ref{s-0-theta-5-lemma}.
\,\,
$\square$

%
%%%%%%%%%%%%%%%%%%%%%%%%%%%%%%%%
%
%
%     Proof of Bilinear Estimate s positive
%
%
%
%%%%%%%%%%%%%%%%%%%%%%%%%%%%%%%%%
%
%
%
\vskip0.05in
\nin
{\bf Case $-j+\frac14<s< 0$.}
We recall that 
in order to prove the bilinear estimate \eqref{bilinear-est}, it suffices to prove $L^2$ inequality \eqref{bilinear-est-L2-form} with $Q$ is estimated in \eqref{Q-def-no1}. Also, similar to the case $s\ge 0$, by further reduction we get
\begin{align}
\label{Q-est-1-neg}
Q(\xi,\xi_1,\tau,\tau_1)
\lesssim&
\frac{|\xi|}{(1+|\tau-\xi^m|)^{b}}
\frac{1}{
(1+|\tau_1-\xi_1^m|)^{b'}
+
\chi_{|\xi_1|\leq 1}(1+|\tau_1|)^{\alpha'}
}
\\
\times&
\frac{(1+|\xi|)^s}{
(1+|\xi_1|)^{s}
(1+|\xi-\xi_1|)^{s}
}
\cdot
\frac{1}
{
(1+|\tau-\tau_1-(\xi-\xi_1)^m|)^{b'}
+
\chi_{|\xi-\xi_1|\leq 1}(1+|\tau-\tau_1|)^{\alpha'}
}.
\notag
\end{align}
Furthermore, like KdVm  on the line \cite{fhy2020} we can restrict our estimations into the set
\begin{equation}
\label{set-E-def}
E
=
\left\{
(\xi,\xi_1,\tau,\tau_1)\in \rr^4:
|\xi_1-\xi|
>
1
\,\,
\text{and}
\,\,
|\xi_1|
> 
1
\right\}.
\end{equation}
Hence, the $\alpha$ terms in the denominators of $Q$
can be dropped and 
and the quantity $Q$ given in \eqref{Q-est-1-neg} becomes
\begin{align}
\label{Q1-modify-no1-s-neg}
Q_4(\xi,\xi_1,\tau,\tau_1)
\lesssim
\frac{|\xi|(1+|\xi|)^s|\xi_1(\xi-\xi_1)|^{-s}}{(1+|\tau-\xi^{m}|)^{b}}
\frac{1}{
(1+|\tau_1-\xi_1^{m}|)^{b'}}
\frac{1}
{
(1+|\tau-\tau_1-(\xi-\xi_1)^{m}|)^{b'}
}.
\end{align}
Moreover, by symmetry (in convolution writing), we may assume that 
\begin{equation}
\label{conv-symmetry-tau-s-neg}
|\tau-\tau_1-(\xi-\xi_1)^{m}|
\leq
|\tau_1-\xi_1^{m}|.
\end{equation}
Finally, following \cite{b1993-kdv}, \cite{kpv1996} and \cite{fhy2020},  in order to prove \eqref{bilinear-est-L2-form} we distinguish two cases (microlocalization):

\noindent
{\bf Microlocalization I. $|\tau_1-\xi_1^{m}|\leq|\tau-\xi^{m}|$}. In this case we define the domain $E_I$ to be
\begin{align}
\label{def-EI-domain}
E_I
\doteq
\big\{
(\xi,\tau,\xi_1,\tau_1) \in{\mathbb{R}}^4: 
&
|\tau-\tau_1-(\xi-\xi)^{m}|\leq|\tau_1-\xi_1^{m}|\leq|\tau-\xi^{m}|,
|\xi_1|>1,
|\xi-\xi_1|>1
\big\}.
\end{align}
{\bf Microlocalization II. $|\tau-\xi^{m}|\leq|\tau_1-\xi_1^{m}|$}. 
In this case we define the domain $E_{II}$ to be
\begin{align}
\label{def-EII-domain}
E_{II}
\doteq
\big\{(\xi,\tau,\xi_1,\tau_1)\in{\mathbb{R}}^4: |\tau-\tau_1-(\xi-\xi)^{m}|\leq|\tau_1-\xi_1^{m}|,
\notag
\\
|\tau-\xi^{m}|\leq|\tau_1-\xi_1^{m}|,
|\xi_1|>1,
|\xi-\xi_1|>1
\big\}.
\end{align}
{\bf Proof of bilinear estimate in Microlocalization I.} 
Here $Q$ is replaced with the $\chi_{E_I}Q$ and our $L^2$ inequality \eqref{bilinear-est-L2-form} reads as
\begin{align}
\label{bilinear-est-neg-EI}
 \Bigg\| 
 \iint_{\rr^2}
 (\chi_{E_I} Q_4)(\xi,\xi_1,\tau,\tau_1)
 c_f(\xi-\xi_1,\tau-\tau_1) c_g(\xi_1,\tau_1) 
 d\xi_1 d\tau_1 \Bigg\|_{L^2_\xi L^2_{\tau}}
\lesssim
\left\|c_f\right\|_{L^2_{\xi} L^2_{\tau}} \left\|c_g\right\|_{L^2_{\xi} L^2_{\tau}}.
\end{align}
As before using the Cauchy-Schwarz inequality with respect to 
$(\xi_1,\tau_1)$ and taking the supremum in  $(\xi,\tau)$ we get
\begin{align*}
%\label{bilinear-est-neg-BI-2}
\Bigg\| 
\iint_{\rr^2}
 (\chi_{E_I} Q_4)(\xi,\xi_1,\tau,\tau_1)
 c_f(\xi-\xi_1,\tau-\tau_1) c_g(\xi_1,\tau_1) 
 d\xi_1 d\tau_1 \Bigg\|_{L^2_\xi L^2_{\tau}}
\lesssim
\|
\Theta_I
\|_{L^{\infty}_{\xi, \tau}}^{1/2}
\|c_f\|_{L^2_{\xi} L^2_{\tau}} 
\|c_g\|_{L^2_{\xi} L^2_{\tau}},
\end{align*}
where $\Theta_I$ is defined and estimated in the following lemma.
\begin{lemma} 
\label{s-neg-case-1-lemma}
If $\max\{\frac12
-
\frac{s-(-j+\frac14)}{12j},\frac{5}{12},\frac{-2s+2+2j}{12j}\}\leq b'<\frac12$ and $-j+\frac14<s<0$, then  for $\xi,\tau\in\rr$
\begin{align}
\label{s-neg-case1-theta-est}
\Theta_I(\xi,\tau)
\doteq&
\frac{\xi^2(1+|\xi|)^{2s}}{(1+|\tau-\xi^{m}|)^{2b}}
\iint_{\rr^2}
\frac{ \chi_{E_I}(\xi,\tau,\xi_1,\tau_1)|\xi_1(\xi-\xi_1)|^{-2s}\,\,\,
 \quad d\tau_1 d\xi_1 
}{(1+|\tau-\tau_1-(\xi-\xi_1)^{m}|)^{2b'}(1+|\tau_1-\xi_1^{m}|)^{2b'}}
\lesssim
1.
\end{align}
\end{lemma}
\nin
The proof of Lemma \ref{s-neg-case-1-lemma} is omited since it is
similar to the proof of Lemma 7.4 in \cite{fhy2020}.   In fact, if we choose
$\frac12-\frac14\beta_1
\le
b'
\le
b
<
\frac12$, where $\beta_1=\beta$, which is defined  in Theorem 2.1 in \cite{fhy2020}, then Lemma \ref{s-neg-case-1-lemma} is reduced to the Lemma 7.4 in \cite{fhy2020}.

\nin
{\bf Proof of bilinear estimate in Microlocalization II.}
Using  duality and Cauchy-Schwarz inequality twice, 
first in $(\xi_1, \tau_1)$ and then in $(\xi, \tau)$, as before,
we have
\begin{align*}
%\label{bilinear-est-neg-BII-2}
\Bigg\| 
\iint_{\rr^2}
 (\chi_{E_{II}} Q_4)(\xi,\xi_1,\tau,\tau_1)
 c_f(\xi-\xi_1,\tau-\tau_1) c_g(\xi_1,\tau_1) 
 d\xi_1 d\tau_1 \Bigg\|_{L^2_\xi L^2_{\tau}}
\lesssim
\|
\Theta_{II}
\|_{L^{\infty}_{\xi_1,\tau_1}}^{1/2}
\|c_f\|_{L^2_{\xi}L^2_{\tau}}
\|c_g\|_{L^2_{\xi}L^2_{\tau}} \cdot 
\end{align*}
where $\Theta_{II}$ is defined and estimated in the following lemma.
\begin{lemma}
\label{s-neg-case-2-lemma} 
If $\max\{\frac12
-
\frac{1}{12j+6},\frac{5}{12},\frac{-4s+m+3}{6m}\}\leq b'\leq b<1/2$ and $-j+\frac14<s<0$, then  for $\xi_1,\tau_1\in\rr$
\begin{align}
\label{s-neg-case2-theta-est}
\Theta_{II}(\xi_1, \tau_1)
\doteq&
\frac{1}{(1+|\tau_1-\xi_1^{m}|)^{2b'}}
\iint_{\rr^2}
\frac{
\chi_{E_{II}}(\xi,\tau,\xi_1,\tau_1)
\,\,
\xi^2
(1+|\xi|)^{2s}
|\xi_1(\xi-\xi_1)|^{-2s}
\,
\, d\tau d\xi}{(1+|\tau-\tau_1-(\xi-\xi_1)^{m}|)^{2b'}(1+|\tau-\xi^{m}|)^{2b}}
\lesssim
1.
\end{align}
\end{lemma}
\nin
{\bf Proof of Lemma \ref{s-neg-case-2-lemma}.} 
Since $0<b'<b<\frac12$, by 
applying  calculus estimate
 \eqref{eq:calc_5}  with $\ell=b$, $\ell'=b'$ $\alpha=\tau_1+(\xi-\xi_1)^m$, $\beta=\xi^m$  and $x=\tau$,  we get
\begin{align}
\label{theta2-est-neg}
\Theta_{II}(\xi_1,\tau_1)
\le&
\frac{1}{(1+|\tau_1-\xi_1^m|)^{2b'}}
\int_\rr
\frac{|\xi|^2(1+|\xi|)^{2s}(\xi_1|\xi-\xi_1|)^{-2s}}{(1+|\tau_1+(\xi-\xi_1)^{m}-\xi^m|)^{2b+2b'-1}}d\xi
\nonumber
\\
=&
\frac{1}{(1+|\tau_1-\xi_1^m|)^{2b'}}
\int_\rr
\frac{|\xi|^2(1+|\xi|)^{2s}(\xi_1|\xi-\xi_1|)^{-2s}}{(1+|\tau_1-\xi_1^m+(\xi-\xi_1)^{m}-\xi^m+\xi_1^m|)^{2b+2b'-1}}d\xi
\nonumber
\\
=&
\frac{1}{(1+|\tau_1-\xi_1^m|)^{2b'}}
\int_\rr
\frac{|\xi|^2(1+|\xi|)^{2s}(\xi_1|\xi-\xi_1|)^{-2s}}{(1+|\tau_1-\xi_1^m+d_m(\xi,\xi_1)|)^{2b+2b'-1}}d\xi.
\end{align}
Then, we complete the proof by following argument  similar to those used  in the proof of Lemma 7.5 in \cite{fhy2020}. 
\,\, $\square$

%%%%%%%%%%%%%%%%%%%%%%%
%
%
%
%     Bilinear Estimates in Y^{s,b} spaces
%
%
%
%%%%%%%%%%%%%%%%%%%%%%%
%
%
%
\section{Proof of Bilinear estimates in temporal $Y^{s,b}$ spaces} 
In this section we  prove Theorem  \ref{bi-est-Y-thm},
that are the bilinear estimates in the spaces $Y^{s,b}$.
These  appears in the basic linear 
 estimate via the time estimate of the forced ivp with zero data, i.e. estimate \eqref{forced-ivp-te}. 
Since the proof of estimates \eqref{bi-est-Y-1} is similar to that of estimate \eqref{bi-est-Y} and for $m=3$ estimate \eqref{bi-est-Y} is proved in \cite{h2006}. 
Here we only provide an outline of the proof for  estimate \eqref{bi-est-Y} with $s\ge 0$.
For  $s\ge 0$ we have the following inequality
\begin{align}
\label{Y-est-1}
\|
w_{fg}
\|_{Y^{s,-b}}^2
\lesssim
\iint_{\rr^2}
\chi_{|\tau|>10^m|\xi|^m}
(1+|\tau|)^{\frac{2s}{m}}
(1+|\tau-\xi^m|)^{-2b}
|\widehat{w}_{fg}(\xi,\tau)|^2
d\xi d\tau
+
\|w_{fg}\|_{X^{s,-b}}^2,
\end{align}
where $w_{fg}=\p_x(f\cdot g)$.
So, to prove the ``temporal" bilinear estimate  \eqref{bi-est-Y} it suffices to show that
\begin{align}
\label{bi-est-Y-reduced}
\left(
\iint_{\rr^2}
\chi_{|\tau|>10^m|\xi|^m}
(1+|\tau|)^{\frac{2s}{m}}
(1+|\tau-\xi^m|)^{-2b}
|\widehat{w}_{fg}(\xi,\tau)|^2
d\xi d\tau
\right)^{1/2}
\lesssim
\|f\|_{X^{s,b',\alpha'}}
\|g\|_{X^{s,b',\alpha'}}.
\end{align}
Like before, writting the  $\|\cdot\|_{s,b,\alpha'}$-norm of $h$ 
as the $L^2$ norm of $c_h$, that is
$
\|h\|_{s,b',\alpha'}\simeq\|c_h(\xi,\tau)\|_{L^2_{\xi,\tau}},
$
where $c_h$ is defined in \eqref{eq:c_h},
the estimate  \eqref{bi-est-Y-reduced} reads as follows
\begin{align}
\label{bilinear-est-c-notation}
\left(
\iint_{\rr^2}
\chi_{|\tau|>10^m|\xi|^m}
(1+|\tau|)^{\frac{2s}{m}}
(1+|\tau-\xi^m|)^{-2b}
|\widehat{w}_{fg}(\xi,\tau)|^2
d\xi d\tau
\right)^{1/2}
\lesssim
\|c_f \|_{L^2_{\xi,\tau}}
\|c_g \|_{L^2_{\xi,\tau}}.
\end{align}
 Next, expressing 
$
\widehat{w}_{fg}(\xi,\tau)
\simeq 
\xi 
\iint_{\rr^2}
\widehat{f}(\xi-\xi_1,\tau-\tau_1)\widehat{g}(\xi_1,\tau_1)d\xi_1 d\tau_1 
$
in terms of $c_f$ and $c_g$
we see that the
inequality \eqref{bilinear-est-c-notation}  takes the following $L^2$  formulation
\begin{align}
\label{bilinear-est-L2-form-Y}
\Big\| 
\iint_{\rr^2}
 Q(\xi,\xi_1,\tau,\tau_1)
 c_f(\xi-\xi_1,\tau-\tau_1) c_g(\xi_1,\tau_1) 
 d\xi_1 d\tau_1 \Big\|_{L^2_{\xi,\tau}}
\lesssim
\|c_f\|_{L^2_{\xi,\tau}} \|c_g\|_{L^2_{\xi,\tau}},
\end{align}
where
\begin{subequations}
\label{Q-def-Y}
\begin{align}
\label{Q-def-term-a-Y}
Q(\xi,\xi_1,\tau,\tau_1)
\doteq&
\chi_{|\tau|>10^m|\xi|^m}
\frac{|\xi|}
{
(1+|\tau-\xi^m|)^{b}
}
\times
(1+|\tau|)^{\frac{s}{m}}
\\
\label{Q-def-term-b-Y}
\times&\frac{1}{
(1+|\xi_1|)^{s}
(1+|\tau_1-\xi_1^m|)^{b'}
+
\chi_{|\xi_1|\leq 1}(1+|\tau_1|)^{\alpha'}
}
\\
\label{Q-def-term-c-Y}
\times&
\frac{1}{
(1+|\xi-\xi_1|)^{s}
(1+|\tau-\tau_1-(\xi-\xi_1)^m|)^{b'}
+
\chi_{|\xi-\xi_1|\leq 1}(1+|\tau-\tau_1|)^{\alpha'}
}.
\end{align}
\end{subequations}

\noindent
In additon, like the bilinear estimate in the space $X^{s,b,\alpha}$, collecting all the factors with $s$ power, making further reduction, 
and observing that if $|\tau|\le 10^m|\xi_1|^m$, then 
$$
\frac{(1+|\tau|)^{s/m}}{
(1+|\xi_1|)^{s}(1+|\xi-\xi_1|)^{s}
}
\lesssim
\frac{(1+10^m|\xi_1|^m)^{s/m}}{
(1+|\xi_1|)^{s}(1+|\xi-\xi_1|)^{s}
}
\lesssim
1,
$$
which implies
$
Q(\xi,\xi_1,\tau,\tau_1)
\lesssim
Q_0(\xi,\xi_1,\tau,\tau_1),
$
where  $Q_0$ is given by \eqref{Q-def-s=0}, 
we reduce bilinear estimate \eqref{bilinear-est-L2-form-Y} to  the bilinear estimate in $X^{s,b,\alpha}$ with $s=0$, i.e. estimate \eqref{bilinear-est-s-0}.
Thus, we assume
$
|\tau|>10^m|\xi_1|^m
$
and   $Q(\xi,\xi_1,\tau,\tau_1)$ becomes $Q_1(\xi,\xi_1,\tau,\tau_1)$, which is given by
\begin{align}
\label{Q1-def}
Q_1(\xi,\xi_1,\tau,\tau_1)
\doteq&
\frac{\chi_{|\tau|>10^m|\xi|^m} |\xi|}{(1+|\tau-\xi^{m}|)^{b}}
\frac{\chi_{|\tau|>10^m|\xi_1|^m}(1+|\tau|)^{s/m}}{
(1+|\xi_1|)^{s}(1+|\xi-\xi_1|)^{s}
}
\\
\cdot&
\frac{1}
{
(1+|\tau-\tau_1-(\xi-\xi_1)^{m}|)^{b'}
(1+|\tau_1-\xi_1^{m}|)^{b'}
}.
\nn
\end{align}
Now, in order to prove estimate \eqref{bilinear-est-L2-form-Y}, it suffices to show that
\begin{align}
\label{bilinear-est-L2-form-Q1}
 \Big\| 
 \iint_{\rr^2}
 Q_1(\xi,\xi_1,\tau,\tau_1)
 c_f(\xi-\xi_1,\tau-\tau_1) c_g(\xi_1,\tau_1) 
 d\xi_1 d\tau_1 \Big\|_{L^2_{\xi,\tau}}
\lesssim
\|c_f\|_{L^2_{\xi,\tau}} \|c_g\|_{L^2_{\xi,\tau}}.
\end{align}
To show this, like before using the Cauchy-Schwarz inequality with respect to 
$(\xi_1,\tau_1)$ and taking the supremum over  $(\xi,\tau)$ we get
\begin{align}
\label{bilinear-est-Q1-Y}
\Big\| 
\iint_{\rr^2}
Q_1(\xi,\xi_1,\tau,\tau_1)
c_f(\xi \hskip-0.03in - \hskip-0.03in \xi_1,\tau \hskip-0.03in - \hskip-0.03in  \tau_1) c_g(\xi_1,\tau_1) 
d\xi_1 d\tau_1 
\Big\|_{L^2_{\xi,\tau}}
\le
\|
\Theta_1
\|_{L^{\infty}_{\xi,\tau}}^{1/2}
\|c_f\|_{L^2_{\xi,\tau}}
\|c_g\|_{L^2_{\xi,\tau}} ,
\end{align}
where $\Theta_1$ is defined and estimated in the following result:
\begin{lemma}
\label{lem: for est-Q1-Y}
If $\textcolor{red}{0}\le s< m$ and $\max\{\frac{2s+m}{6m},\frac{m+2}{6m}\}\le b'<b<\frac12$ satisfy , then we have
\begin{align}
\label{ine: for est-Q1-Y}
\Theta_1(\xi,\tau)
\doteq&
\frac{\chi_{|\tau|>10^m|\xi|^m} |\xi|^2(1+|\tau|)^{2s/m}}{(1+|\tau-\xi^{m}|)^{2b}}
\iint_{\rr^2}
\frac{\chi_{|\tau|>10^m|\xi_1|^m}}{
(1+|\xi_1|)^{2s}(1+|\xi-\xi_1|)^{2s}
}
\\
&
\cdot
\frac{1}
{
(1+|\tau-\tau_1-(\xi-\xi_1)^{m}|)^{2b'}
(1+|\tau_1-\xi_1^{m}|)^{2b'}
}
\xi_1 d\tau_1 
\lesssim
1.
\notag
\end{align}
\end{lemma}

\nin
{\bf Proof of Lemma \ref{lem: for est-Q1-Y}.}
Applying calculus estimate \eqref{eq:calc_5} with $\ell=\ell'=b'$, $a=\tau-(\xi-\xi_1)^m$, $c=\xi_1^m$ and $x=\tau_1$
we arrive at the following estimate
\begin{align}
\label{theta1-est-0}
\Theta_1(\xi,\tau)
\le&
\frac{\chi_{|\tau|>10^m|\xi|^m} |\xi|^2(1+|\tau|)^{2s/m}}{(1+|\tau-\xi^{m}|)^{2b}}
\\
&\int_{\rr}
\frac{\chi_{|\tau|>10^m|\xi_1|^m}}{
(1+|\xi_1|)^{2s}(1+|\xi-\xi_1|)^{2s}
}
\cdot
\frac{1}
{
(1+|\tau-(\xi-\xi_1)^{m}-\xi_1^m|)^{4b'-1}
}
d\xi_1.
\nn
\end{align}
Then,  for $\textcolor{red}{0}\le s< m$ and $\max\{\frac{2s+m}{6m},\frac{m+2}{6m}\}\le b'<b<\frac12$, using $|\tau|>10^m|\xi|^m$ and $ |\tau|>10^m|\xi_1|^m$ we get the desired estimate \eqref{ine: for est-Q1-Y}. Here we omit the detail of the proof. 
\,\,
$\square$

%
%%%%%%%%%%%%%%%%%%%%%%%%%%%%%%
%
%
%
%Proof of  Well-posedness Theorem
%
%
%
%
%
%%%%%%%%%%%%%%%%%%%%%%%%%%%%%%%
\section{Well-posedness in modified Bourgain spaces -- Proof of Theorem \ref{thm-kdvm-half-line}
}
We only prove  well-posedness for $-1\le s< \frac12$. The proof for  $s\in(-j+\frac14,-1)\cup (\frac12,j+1)$, $s\neq \frac32, \frac52, \dots,j-\frac12$, is similar. Also, we assume that 
$$
0
<
T
<
1/2.
$$
\nin
 {\bf Small data.}
First, we  prove Theorem \ref{thm-kdvm-half-line}
for  initial and boundary data such that
\begin{equation}
\label{smallness} 
\|u_0\|_{H^s(\rr^+)}
+
\sum\limits_{\ell=0}^{j-1}\|g_\ell\|_{H^{\frac{s+j-\ell}m}(0,T)}
\le
\frac{1}{144C^2},
\,\,
\text{with}
\,\,
C=c_{s,b,\alpha}+\frac 12 c_{s,b,\alpha}^2,
\end{equation}
where $c_{s,b,\alpha}$  
is the constant appearing in the estimate \eqref{forced-linear-kdvm-est} and the  bilinear estimates 
\eqref{bilinear-est}.
Under the above smallness condition \eqref{smallness}, we prove that the the  integral equation 
\begin{equation}
\label{iteration-map}
u
=
\Phi u
\doteq
S
\Big[
u_0,g_0,\dots,g_{j-1};-\frac12\p_x(u^2)
\Big],
\end{equation}
has a unique  solution in the space $X^{s,b,\alpha}(\rr^+\times (0,T))$. For this, we shall prove that  the iteration map 
$\Phi$  has a fixed point in $X^{s,b,\alpha}(\rr^+\times (0,T))$. 
In fact,  for any  $u$ in the (closed) ball 
\begin{equation}
\label{smallness-ball} 
B=\Big\{u\in X^{s,b,\alpha}(\rr^+\times(0,T)):\|u\|_{X^{s,b,\alpha}(\rr^+\times(0,T))}
\le
\frac{1}{24C}\Big\},
\end{equation}
using linear  estimate  \eqref{forced-linear-kdvm-est}
with forcing replaced by $-\frac12\p_x(u^2)$
and bilinear estimates {\eqref{bilinear-est}
we get
\begin{align*}
\|\Phi u\|_{X^{s,b,\alpha}(\rr^+\times (0,T))}
\le &
c_{s,b,\alpha}
\Big(  \|u_0\|_{H^s(\rr^+)} +\sum\limits_{\ell=0}^{j-1} \|g_\ell\|_{H^{\frac{s+j-\ell}{m}}(0,T)}
+  
\frac{1}2
\|\p_x(u^2)\|_{X^{s,-b,\alpha-1}(\rr^+\times(0,T))}
\Big)
\\
\le&
c_{s,b,\alpha}
 \Big( \|u_0\|_{H^s(\rr^+)} 
+ 
\sum\limits_{\ell=0}^{j-1}\|g_\ell\|_{H^{\frac{s+j-\ell}{m}}(0,T)}
+  
\frac{1}2
\|\p_x(\tilde u^2)\|_{X^{s,-b,\alpha-1}(\rr^2)}
\Big)
\\
\overset{\eqref{bilinear-est}}{\le}&
C \Big( \|u_0\|_{H^s(\rr^+)} 
+ 
\sum\limits_{\ell=0}^{j-1}\|g_\ell\|_{H^{\frac{s+j-\ell}{m}}(0,T)}
+ 
 \|\tilde u\|^2_{X^{s,b,\alpha}(\rr^2)}\Big),
\end{align*}
where $\tilde u$
is an extension of  $u$ from $\rr^+\times(0,T)$ to $\rr^2$ such that 
\begin{equation}
\label{u-extension}
\|\tilde u\|_{X^{s,-b,\alpha-1}(\rr^2)}
\le
2
\|u\|_{X^{s,b,\alpha}(\rr^+\times(0,T))}.
\end{equation}
Furthermore, using estimate \eqref{u-extension} we get
\begin{align}
\label{contranction-est}
\|\Phi u\|_{X^{s,b,\alpha}(\rr^+\times (0,T))}
\le&
C \Big( \|u_0\|_{H^s(\rr^+)} +
\sum\limits_{\ell=0}^{j-1}\|g_\ell\|_{H^{\frac{s+j-\ell}{m}}(0,T)}
+  4\|u\|^2_{X^{s,b,\alpha}(\rr^+\times(0,T))}\Big).
\end{align}
And, since $u\in B$, we have
$
\|\Phi u\|_{X^{s,b,\alpha}(\rr^+\times (0,T))}
\le
C
\Big(
\frac{1}{144C^2}+\frac{1}{144C^2}
\Big)
\le
\frac{1}{24C}.
$
Thus  $\Phi$ maps the  ball $B$ into itself.

To show that $\Phi$ is a contraction,  for any $u,v\in B$,
using   linear estimate \eqref{forced-linear-kdvm-est}
with forcing replaced by $-\frac12\p_x(u^2-v^2)$
we get
\begin{align*}
\|\Phi u-\Phi v\|_{X^{s,b,\alpha}(\rr^+\times\rr)}
\le& 
\frac{c_{s,b,\alpha}}2
\|\p_x(u^2-v^2)\|_{X^{s,-b,\alpha-1}(\rr^+\times\rr)}
\\
\le&
\frac{c_{s,b,\alpha}}2
\|\p_x(\tilde u^2-\tilde v^2)\|_{X^{s,-b,\alpha-1}(\rr^2)},
\end{align*}
where $\tilde u$ is the extension of $u$ from $\rr^+\times(0,T)$ to $\rr^2$ satisfying \eqref{u-extension}.
The extension of $v$  is obtained as follows. 
First, we extend
$w=v-u$ from $\rr^+\times(0,T)$ to $\rr^2$ such that 
\begin{equation}
\label{v-extension}
\|\tilde w\|_{X^{s,b,\alpha}(\rr^2)}
\le
2
\|v-u\|_{X^{s,b,\alpha}(\rr^+\times(0,T))}.
\end{equation}
Then defining $\tilde v\doteq \tilde w+\tilde u$, we see that 
$\tilde v$ extends $v$ from $\rr^+\times(0,T)$ to $\rr^2$.
In addition, using the triangle inequality, we get
\begin{align}
\label{v1-bound}
\|\tilde v\|_{X^{s,b,\alpha}(\rr^2)}
\le&
2
\|v-u\|_{X^{s,b,\alpha}(\rr^+\times(0,T))}
+
2
\|u\|_{X^{s,b,\alpha}(\rr^+\times(0,T))}
\le
\frac{1}{4C}.
\end{align}
Combining estimate  \eqref{u-extension} and  estimate \eqref{v-extension} with bilinear estimate \eqref{bilinear-est} again, we get
\begin{align}
\label{contraction-small-est-2}
\|\Phi u-\Phi v\|_{X^{s,b,\alpha}(\rr^+\times(0,T))}
\le& 
C\left(\|\tilde u\|_{s,b,\alpha}+\|\tilde v\|_{s,b,\alpha}\right)\cdot\|\tilde u-\tilde v\|_{s,b,\alpha}
\\
\le&
C\cdot\frac{4}{12C}\cdot\|\tilde u-\tilde v\|_{s,b,\alpha}
\le
\frac{2}{3}\|u-v\|_{X^{s,b,\alpha}(\rr^+\times(0,T))},
\nn
\end{align}
which shows that  $\Phi$ is a contraction on $B$.  
Since $B$ is a complete Banach space, 
by the  contraction mapping theorem
there is a unique $u\in B$ such that $\Phi u=u$.

\vskip0.05in
\nin
{\it Lip-continuous  dependence on data.}
Let   $u_0(x)$, $g_\ell(t)$ and $v_0(x)$, $h_\ell(t)$
be two sets of data satisfying  the smallness condition \eqref{smallness}.
If $u$ is the solution that corresponds to  
 $u_0(x)$, $g_\ell(t)$, which we denote by
$
u(x, t)
=
\psi(t)
S
\Big[
u_0,g_0,\dots,g_{j-1};-\frac12\p_x(u^2)
\Big],
$
and  $v$ is the solution that corresponds to 
$v_0(x)$, $h_\ell(t)$, that is
$
v(x, t)
=
\psi(t)
S
\Big[
v_0,h_0,\dots,h_{j-1};-\frac12\p_x(u^2)
\Big]
$
then 
\begin{equation}
\label{u-v-eqn}
u(x, t)-v(x, t)
=
\psi(t)
S
\Big[
u_0-v_0,g_0-h_0,\dots,g_{j-1}-h_{j-1};-\frac12\p_x(u^2-v^2)
\Big].
\end{equation}
Using linear estimate \eqref{forced-linear-kdvm-est}, extensions $\tilde u$, $\tilde v$ (as above), and bilinear estimate 
 \eqref{bilinear-est},  we have
\begin{align*}
%\label{Tu-Tv-est} 
\|u-v\|_{X^{s,b,\alpha}(\rr^+\times(0,T))}
&\le
C \big(
\|u_0-v_0\|_{H^s(\rr^+)}
+
\sum\limits_{\ell=0}^{j-1}
\|g_\ell-h_\ell\|_{H^{\frac{s+j-\ell}m}(0,T)}
\big)
+
C \|\tilde u+\tilde v\|_{s,b,\alpha} \|\tilde u-\tilde v\|_{s,b,\alpha}
\\
&\le
C\big(
\|u_0-v_0\|_{H^s(\rr^+)}
+
\sum\limits_{\ell=0}^{j-1}
\|g_\ell-h_\ell\|_{H^{\frac{s+j-\ell}m}(0,T)}
\big)
+
\frac{2}{3}\|u-v\|_{X^{s,b,\alpha}(\rr^+\times(0,T))}.
\end{align*}
Moving  all $ \|u-v\|_{X^{s,b,\alpha}(\rr^+\times(0,T))}$  to  the lhs gives
\begin{equation}
\label{lip-small-est}
\|u-v\|_{X^{s,b,\alpha}(\rr^+\times(0,T))}
\le
3C
\big(
\|u_0-v_0\|_{H^s(\rr^+)}
+
\sum\limits_{\ell=0}^{j-1}
\|g_\ell-h_\ell\|_{H^{\frac{s+j-\ell}m}(0,T)}
\big),
\end{equation}
which completes the proof of 
Lip-continuous  dependence on data.

%%%%%%%%%%%%%%%%%%%%%%
%
%
%
%
%      Arbitrary size data
%
%
%
%%%%%%%%%%%%%%%%%%%%%
%
%
%
\nin
{\bf Large data.}
For any size  initial data $u_0\in H^s$, boundary data $g$, and
for $T^*$ such that 
\begin{equation}
\label{T-less-than-1/2}
0<T^*\le T<1/2,
\end{equation}
we replace the integral equation  \eqref{iteration-map}
with its following localization 
\begin{align}
\label{iteration-map-T-loc}
u(x, t)
=
\Phi_{T^*} u
\doteq
S
\Big[
u_0,g_0,\dots,g_{j-1};
-\frac12\p_x(\psi_{2T^*}\cdot u^2)
\Big],
\quad   |t|\le T^*,
\end{align}
where   $\psi_{T^*}(t)=\psi(t/{T^*})$ with $\psi(t)$  being 
our familiar cutoff function in $\in C_0^\infty(-1, 1)$ with 
$0\le \psi(t)\le 1$, $\psi(t)=1$ for $|t|\le 1/2$.
First,  we notice that for $|t|\le {T^*}$, the fixed point of the iteration map \eqref{iteration-map-T-loc} is the solution to the KdVm ibvp \eqref{KdVm}.  
Thus,
$\Phi_{T^*}(u)=\Phi(u)$
if
$|t|\le {T^*}$,
i.e.  when $|t|\le {T^*}$, then $\Phi_{T^*}(u)$  becomes 
the iteration map \eqref{iteration-map}
.
Next, we shall choose appropriate $ {T^*}$ and  use the contraction mapping theorem to show that there is a fixed point of the iteration map \eqref{iteration-map-T-loc} in the ball $B(r)\subseteq X^{s,b,\alpha}(\rr^+\times(0,T))$. 
In fact, using  the linear estimate \eqref{forced-linear-kdvm-est} with forcing replaced by $-\frac12\p_x(u^2)$, for $b<b_1$, which are given below in \eqref{b-b'-choice-wp},   we get
\begin{align}
\label{onto-map-fin-est-large}
&\|\Phi_ {T^*}(u)\|_{X^{s,b,\alpha}(\rr^+\times(0,T))}
\le
\|\Phi_ {T^*}(u)\|_{X^{s,b_1,\alpha}(\rr^+\times(0,T))}
\\
\le&
c_{s,b,\alpha}
\Big(
\|
u_0
\|_{H^s(\rr^+)}
+
\sum\limits_{\ell=0}^{j-1}\|g_\ell\|_{H_t^\frac{s+j-\ell}{m}(0,T)}
+
\frac{1}2
\left\|
\psi_{2{T^*}}(t)
\p_x(u^2(t))
\right\|_{X^{s,-b_1,\alpha-1}(\rr^+\times(0,T))}
\Big)
\nn
\\
\le&
c_{s,b,\alpha}
\Big(
\|
u_0
\|_{H^s(\rr^+)}
+
\sum\limits_{\ell=0}^{j-1}\|g_\ell\|_{H_t^\frac{s+j-\ell}{m}(0,T)}
+
\frac{1}2
\left\|
\psi_{2{T^*}}(t)
\p_x(\tilde u^2(t))
\right\|_{s,-b_1,\alpha-1}
\Big)
,
\nn
\end{align}
where $\tilde{u}$ is the extension of $u$ from $\rr^+\times(0,T)$ to $\rr^2$, 
which  satisfies \eqref{u-extension}.
Now, we estimate the $\|\cdot\|_{s,-b_1,\alpha-1}$. For this we shall need the following result.
\begin{lemma}
\label{tao-lemma-modified}
Let $\eta(t)$ be a function in the Schwartz space $\mathcal{S}(\rr)$. If 
$
-\frac12
<
b'
\le 
b
<
\frac12
$
and 
$
-\frac12
<
\alpha'
-
1
\le 
\alpha
-
1
<
\frac12
$
(%
$
\frac12
<
\alpha'
\le
\alpha
<
1
$ is sufficient condition)
then for any $0<{T^*}\le 1$ we have
\begin{equation}
\label{tao-est-modified}
\|\eta(t/{T^*})u\|_{X^{s,b',\alpha'-1}}
\le
c_1(\eta,b,b',\alpha,\alpha')
\,\,
\max\{
{T^*}^{b-b'}, {T^*}^{\alpha-\alpha'}
\}
\|u\|_{X^{s,b,\alpha-1}}.
\end{equation}
\end{lemma}
\noindent
The proof of this result is based on 
the following multiplier estimate in $X^{s, b}$ spaces,
which can be found in  \cite{tao-book} (see page 101,
 Lemma 2.11), i.e.
\begin{equation}
\label{tao-est}
\|\eta(t/{T^*})u\|_{X^{s,b'}}
\le
c_1(\eta,b,b')
\,\,
{T^*}^{b-b'}
\|u\|_{X^{s,b}}.
\end{equation}
Applying estimate \eqref{tao-est-modified} with 
the following choice
\begin{equation}
\label{b-b'-choice-wp}
b
=
\frac12-\beta
\text{  (in place of $b'$) }
\quad
\text{and}
\quad
b_1
=
\frac12-\frac12\beta
\text{  (in place of $b$)},
\end{equation}
and 
\begin{equation}
\label{a-a1-choose}
\alpha
=
\frac12+\frac12\beta
\text{  (in place of $\alpha'$) }
\quad
\text{and}
\quad
\alpha_1
=
\frac12+\beta
\text{  (in place of $\alpha$) },
\end{equation}
where $\beta$ is defined in \eqref{beta-choice} and it is only depending on $s$ for fixed $m$. From \eqref{onto-map-fin-est-large}
we obtain
\begin{align*}
%\label{onto-map-fin-est}
\|\Phi_ {T^*}(u)\|_{X^{s,b,\alpha}(\rr^+\times(0,T))}
\le
c_{s,b,\alpha}
\Big(
\|
u_0
\|_{H^s(\rr^+)}
+
\sum\limits_{\ell=0}^{j-1}\|g_\ell\|_{H_t^\frac{s+j-\ell}{m}(0,T)}
+
\frac{c_1}2
{T^*}^{\frac12\beta}
\left\|
\p_x(\tilde u^2(t))
\right\|_{s,-b,\alpha_1-1}
\Big).
\nn
\end{align*}
Then  the bilinear estimates \eqref{bilinear-est}
reads as follows
\begin{equation}
\label{bilinear-est-special}
\|\p_x(f\cdot g) \|_{s,-b,\alpha_1-1}
\le
c_{s,b,\alpha}
\| f \|_{s,b,\alpha} \| g \|_{s,b,\alpha},
\quad
f, g \in X^{s,b,\alpha},
\end{equation}
and we will use it in this form.
Therefore, we get
\begin{align}
\label{onto-map-fin-est-1}
\|\Phi_ {T^*}(u)\|_{X^{s,b,\alpha}(\rr^+\times(0,T))}
\le &
c_2
\Big(
\|
u_0
\|_{H^s(\rr^+)}
+
\sum\limits_{\ell=0}^{j-1}\|g_\ell\|_{H_t^\frac{s+j-\ell}{m}(0,T)}
+
{T^*}^{\frac12\beta}
\|
\tilde u
\|_{s,b,\alpha}^2
\Big)
\\
\le&
c_2
\Big(
\|
u_0
\|_{H^s(\rr^+)}
+
\sum\limits_{\ell=0}^{j-1}\|g_\ell\|_{H_t^\frac{s+j-\ell}{m}(0,T)}
+
4 {T^*}^{\frac12\beta}
\|
u
\|^2_{X^{s,b,\alpha}(\rr^+\times(0,T))}
\Big),
\nn
\end{align}
where $c_2=c_2(s,b,\alpha)\doteq c_{s,b,\alpha}+\frac12c_1\cdot c_{s,b,\alpha}^2$.
From \eqref{onto-map-fin-est-1} we see that
for  the map $\Phi_{T^*}$ \eqref{iteration-map-T-loc} to be  onto, it suffices to have
$$
c_2
\Big(
\|
u_0
\|_{H^s(\rr^+)}
+
\sum\limits_{\ell=0}^{j-1}\|g_\ell\|_{H_t^\frac{s+j-\ell}{m}(0,T)}
+
4 {T^*}^{\frac12\beta}
\|
u
\|^2_{X^{s,b,\alpha}(\rr^+\times(0,T))}
\Big)
\le 
r.
$$
And, since $u\in B(r)$ it suffices to have
\begin{align}
\label{onto-map-condition}
c_2
\Big(
\|
u_0
\|_{H^s(\rr^+)}
+
\sum\limits_{\ell=0}^{j-1}\|g_\ell\|_{H_t^\frac{s+j-\ell}{m}(0,T)}
\Big)
+
4
c_2
{T^*}^{\frac12\beta}
r^2
\le 
r.
\end{align}
To show that $\Phi_{T^*}$ is a contraction, again,
using   linear estimate \eqref{forced-linear-kdvm-est}
with forcing replaced by $-\frac12\p_x(u^2-v^2)$, extensions $\tilde u$, $\tilde v$ (as above), 
for $b\le b_1$  we have
\begin{align}
\label{contraction-1}
\|\Phi_{T^*}(u)-\Phi_{T^*}(v)\|_{X^{s,b,\alpha}(\rr^+\times(0,T))}
{\le }&
\frac{c_{s,b,\alpha}}2
\left\|
\psi_{2{T^*}}(t)
\p_x(u^2(t)-v^2(t))
\right\|_{X^{s,-b_1,\alpha-1}(\rr^+\times(0,T))}
\\
\le&
\frac{c_{s,b,\alpha}}2
\left\|
\psi_{2{T^*}}(t)
\p_x(\tilde{u}^2(t)-\tilde{v}^2(t))
\right\|_{X^{s,-b_1,\alpha-1}},
\nn
\end{align}
Applying estimate \eqref{tao-est-modified} with $b$, $b_1$ given in \eqref{b-b'-choice-wp} and $\alpha$, $\alpha_1$ given by \eqref{a-a1-choose}, we get
\begin{align}
\label{contraction-map-fin-est}
\|\Phi_{T^*}(u)-\Phi_{T^*}(v)\|_{X^{s,b,\alpha}(\rr^+\times(0,T))}
\le
\frac{c_1\cdot c_{s,b,\alpha}}2
{T^*}^{\frac12\beta}
\left\|
\p_x[(\tilde u(t)+\tilde v(t))(\tilde u(t)-\tilde v(t))]
\right\|_{s,-b,\alpha_1-1}.
\end{align}
Next, using the bilinear estimate \eqref{bilinear-est-special},  from \eqref{contraction-map-fin-est} we get
\begin{align}
\label{contraction-map-fin-est-1}
\|\Phi_{T^*}(u)-\Phi_{T^*}(v)\|_{X^{s,b,\alpha}(\rr^+\times(0,T))}
\le&
c_2
{T^*}^{\frac12\beta}
\|\tilde u+\tilde v\|_{s,b,\alpha}
\|\tilde u-\tilde v\|_{s,b,\alpha}
\\
\le&
16
c_2
{T^*}^{\frac12\beta}
 r
\|u-v\|_{X^{s,b,\alpha}(\rr^+\times(0,T))}.
\nn
\end{align}
Thus,
in order to make the iteration map $\Phi_ {T^*}$ a contraction map, it suffices to have
\begin{equation}
\label{contraction-T-condition}
16
c_2
{T^*}^{\frac12\beta}
 r
\le 
\frac12.
\end{equation}
Combining conditions \eqref{onto-map-condition} with \eqref{contraction-T-condition}, we see that it suffices to have
$$
c_2
\Big(
\|
u_0
\|_{H^s(\rr^+)}
+
\sum\limits_{\ell=0}^{j-1}\|g_\ell\|_{H_t^\frac{s+j-\ell}{m}(0,T)}
\Big)
+
\frac18r
\le 
r
\iff
r
\ge
\frac87
c_2
\Big(
\|
u_0
\|_{H^s(\rr^+)}
+
\sum\limits_{\ell=0}^{j-1}\|g_\ell\|_{H_t^\frac{s+j-\ell}{m}(0,T)}
\Big)
.
$$
So, we choose  the radius to be
\begin{equation}
\label{ball-radius-choice}
r
\doteq
2c_2
\Big(
\|
u_0
\|_{H^s(\rr^+)}
+
\sum\limits_{\ell=0}^{j-1}\|g_\ell\|_{H_t^\frac{s+j-\ell}{m}(0,T)}
\Big).
\end{equation}
Then, from \eqref{contraction-T-condition}  it suffices to have 
$
{T^*}^{\frac12\beta}
\le
(32c_2r)^{-1},
$
which follows from choosing
\begin{equation}
\label{T-choice-wp-1}
{T^*}
=
\frac 12
(1+32c_2r)^{-\frac{2}{\beta}}
<
\frac 12.
\end{equation}
Combining this choice of ${T^*}$ together with choice  \eqref{ball-radius-choice} for $r$  we get
\begin{align}
\label{T-choice-wp-2}
{T^*}
=&
\frac 12
\Big[
1+64c_2^2
\big(
\|u_0\|_{H^s(\rr^+)}+\sum\limits_{\ell=0}^{j-1}\|g_\ell\|_{H_t^\frac{s+j-\ell}{m}(0,T)}
\big)
\Big]^{-\frac{2}{\beta}}
\\
\ge&
c_0\cdot
\Big(
1+\|u_0\|_{H^s(\rr^+)}+\sum\limits_{\ell=0}^{j-1}\|g_\ell\|_{H_t^\frac{s+j-\ell}{m}(0,T)}
\Big)^{-\frac{4}{\beta}}
\doteq T_0,
\nn
\end{align}
for some $c_0$ depending on $c_2(s,b,\alpha)$ and $\beta=\beta(s)$,
that is  $c_0=c_0(s,b,\alpha)$.
Thus we choose the  lifespan as stated in 
\eqref{lifespan-est}.
This completes the proof of well-posedness for $-\frac12<s<\frac12$. \,\, $\square$

Lip continuity of the data to solution map and uniquenees
is similar to the well-posedness on the line
described in \cite{fhy2020}.

%\newpage
%%
%
%%%%%%%%%%%%%%%%%%%%%%%
%
%
%
%   Derivation of the Fokas Solution Formula
%
%
%
%%%%%%%%%%%%%%%%%%%%%%%%
%
%
%
\section{Derivation of the Fokas Solution Formula}
\label{kdvm-sln-derivation}
Here, we provide an outline 
of UTM for the solution to the forced linear KdVm ibvp in three steps. 
First, we  use the Fourier transform on the half-line to 
get a solution formula  to  ibvp \eqref{LKdVm} via the Fourier inversion  formula
on the  real line. 
Then, we deform the contour
via the Cauchy's Theorem and derive a  formula for the solution 
integrating over the contours $\p D_{2p}^+$, $p=1,2,\cdots,j$ in 
the upper half of the complex plane.  Finally, we eliminate the unknown boundary data and get the desired solution formula \eqref{UTM-sln-compact}.

\vskip.05in
\nin
{\bf Step 1: Solving  KdVm ibvp \eqref{LKdVm} via half-line Fourier transform.}
If $\tilde u$ is a solution to the  LKdVm formal adjoint equation
\begin{align}
\label{adjoint-j eqn}
\p_t \tilde u+(-1)^{j+1}\p^{2j+1}_x\tilde u
=
0,
\end{align}
then multiplying it by $u$ and equation  \eqref{LKdVm eqn} by 
$\tilde u$, and adding the resulting equations gives
\begin{align*}
\tilde u\p_t u+u\p_t\tilde u+(-1)^{j+1}(\tilde u \p_x^mu+u\p_x^m\tilde u)
=
\tilde uf,
\end{align*}
or
\begin{align}
\label{combine-j eqn}
(\tilde uu)_t+(-1)^{j+1}[\tilde u\p_x^{2j} u+\cdots+(-1)^{n}\p_x^n\tilde u\p_x^{2j-n}u+\cdots+\p_x^{2j}\tilde uu]_x
=
\tilde uf.
\end{align}
Then, choosing as $\tilde u$ the  exponential solutions 
to transpose equation \eqref{adjoint-j eqn}
\begin{align}
\label{adjoint-j sln}
\tilde u
=
e^{-i\xi x-i\xi^mt}, 
\quad
\xi\in\cc,
\end{align}
and  substituting them into identity
\eqref{combine-j eqn} we get the {\bf divergence form:}
\begin{align}
\label{divergence form-j}
%\boxed{
(e^{-i\xi x-i\xi^mt}u)_t+(-1)^{j+1}(e^{-i\xi x-i\xi^mt})[\p_x^{2j} u+\cdots+i^n\xi^n\p_x^{2j-n}u+\cdots+(-1)^j\xi^{2j} u])_x
=
\tilde uf.
%}
\end{align}
Integrating the divergence form \eqref{divergence form-j} from $x=0$ to $\infty$ gives the $t$-equation
\begin{align}
\label{tode-j}
(e^{-i\xi^mt}\widehat u(\xi,t))_t
=
e^{-i\xi^mt}\widehat f(\xi,t)+(-1)^{j+1}e^{-i\xi^mt}g(\xi,t),
\end{align}
where $\widehat u$ and $\widehat f$ are the  half-line
Fourier transforms  of $u$ and $f$, which are defined in \eqref{FT-halfline}, and $g$ is the following combination of $m=2j+1$ boundary data (some of which are not given)
\begin{align}
g(\xi,t)
\doteq
\p_x^{2j}u(0,t)+\cdots+i^n\xi^n\p_x^{2j-n}u(0,t)+\dots+(-1)^{j}\xi^{2j}u(0,t).
\end{align}
Integrating \eqref{tode-j} from $0$ to $t$, $0\leq t\leq T$, we obtain the so called  {\bf global relation:}
\begin{align}
\label{global relation-j}
e^{-i\xi^mt}\widehat u(\xi,t)
=&
\widehat u_0(\xi)
+
F(\xi,t)
\\
+&(-1)^{j+1}[\tilde g_{2j}(\xi^m,t)+\cdots+i^n\xi^n\tilde g_{2j-n}(\xi^m,t)
+\cdots+(-1)^j\xi^{2j}\tilde g_0(\xi^m,t)]
,
\quad
\text{Im}(\xi)\leq 0,
\nonumber
\end{align}
where $F(\xi,t)$ and $\tilde g_{\ell}(\xi,t)$ are given in \eqref{F-time-transform} and \eqref{g-time-transform} respectively.
Now, inverting  \eqref{global relation-j} we get
\begin{align}
\label{sln-line-j}
u(x,t)
&=
\frac{1}{2\pi}\int_{-\infty}^\infty e^{i\xi x+i\xi^mt}
[\widehat u_0(\xi)+F(\xi,t)]d\xi
\\
&+
\frac{(-1)^{j+1}}{2\pi}\int_{-\infty}^\infty e^{i\xi x+i\xi^mt}
[\tilde g_{2j}(\xi^m,t)+\cdots+(i\xi)^\ell\tilde g_{2j-\ell}(\xi^m,t)
+\cdots+(-1)^j\xi^{2j}\tilde g_0(\xi^m,t)]d\xi.
\nonumber
\end{align}
{\bf Step 2:  Deforming integration over the contours 
$\p D_{2p}^+$ in the upper half-pane.}
Formula \eqref{sln-line-j} contains $(j+1)$ unknown data.
To eliminate them, we  deform the integration contour from $\rr$ to $\p D_{2p}^+$, $p=1,2,\cdots,j$. This is expressed by the following result.
\begin{lemma}
\label{UTM-lem}
The solution $u(x,t)$ to ibvp  \eqref{LKdVm} can be written in the form
\begin{align}
\label{sln1-j}
u(x,t)
&=
\frac{1}{2\pi}\int_{-\infty}^\infty e^{i\xi x+i\xi^mt}
[\widehat u_0(\xi)+F(\xi,t)]d\xi
\\
&+
\frac{(-1)^{j+1}}{2\pi}\sum\limits_{p=1}^j\int_{\p D_{2p}^+} e^{i\xi x+i\xi^mt}
[\tilde g_{2j}(\xi^m,t)+\cdots+(i\xi)^\ell\tilde g_{2j-\ell}(\xi^m,t)
+\cdots+(-1)^j\xi^{2j}\tilde g_0(\xi^m,t)]d\xi,
\nonumber
\end{align}
where the domains $D^+_2, \cdots,  D_{2j}^+$ are shown in Figure 1.1 (if $j$ is odd) or Figure 1.2   (if $j$ is even), and the orientation of the boundary $\p D^+_{2p}$ is given by the left-hand rule.
\end{lemma}
\noindent
The proof of above lemma can be found in \cite{y2020}. Also, a similar lemma for the KdV equation can be found in \cite{fhm2016}.

\vskip.05in
\noindent
{\bf Step 3: Eliminating the unknown boundary data.}
For each $p=1,2,\cdots,j$, we construct a linear system with $j+1$
 equations (as many as the unkown data). For this we apply the invariant transformations of $\xi^m$, i.e. $\xi\longrightarrow\alpha_{p,n}\xi$,  $n=1,2,\cdots,j+1$, where $\alpha_{p,n}$ are defined in \eqref{a-rot-angles}, and use the global relation \eqref{global relation-j}.
Thus, we obtain the linear system of the following $j+1$ equations
\begin{align*}
\begin{cases}
e^{-i\xi^mt}\widehat u(\alpha_{p,1}\xi,t)
=
\widehat u_0(\alpha_{p,1}\xi)
+
F(\alpha_{p,1}\xi,t)
\nonumber
\\
\hskip.7in
+(-1)^{j+1}[\tilde g_{2j}(\xi^m,t)+\cdots+(\alpha_{p,1})^\ell(i\xi)^\ell\tilde g_{2j-\ell}(\xi^m,t)
+\cdots+(\alpha_{p,1})^{2j}(i\xi)^{2j}\tilde g_0(\xi^m,t)],
\\
\hskip1in \cdots
\\
e^{-i\xi^mt}\widehat u(\alpha_{p,j+1}\xi,t)
=
\widehat u_0(\alpha_{p,j+1}\xi)
+
F(\alpha_{p,j+1}\xi,t)
\nonumber
\\
\hskip.7in
+(-1)^{j+1}[\tilde g_{2j}(\xi^m,t)+\cdots+(\alpha_{p,j+1})^\ell(i\xi)^\ell\tilde g_{2j-\ell}(\xi^m,t)
+\cdots+(\alpha_{p,j+1})^{2j}(i\xi)^{2j}\tilde g_0(\xi^m,t)].
\end{cases}
\end{align*}
Solving these equations  for $(i\xi)^\ell\tilde g_{2j-\ell}(\xi^m,t)
$, $\ell=0,1,\cdots,j$, and substituting the obtained solutions into formula \eqref{sln1-j}, we get the desired Fokas solution formula \eqref{UTM-sln-compact} involving only the given data
and no unknown data.

%
%%%%%%%%%%%%%%%%%%%%%%
%
%
%		Acknowledgements
%
%
%%%%%%%%%%%%%%%%%%%%%%  
%
\vspace*{1mm}
\noindent
\textbf{Acknowledgements.} The first author was partially supported by a grant from the Simons Foundation (\#524469 to Alex Himonas).

%
%\vspace*{1mm}
%\noindent
%\textbf{Conflicts of interest.} On behalf of both authors, the corresponding author states that there is no conflict of interest.
%
%

%
%%%%%%%%%%%%%%%%%%%%%%%%  
%
%			Bibliography
%
%%%%%%%%%%%%%%%%%%%%%%%%  
%
%
%

%

\vspace{3mm}
\noindent
A. Alexandrou Himonas  \hfill Fangchi Yan\\
Department of Mathematics  \hfill Department of Mathematics\\
University of Notre Dame  \hfill West Virginia University\\
Notre Dame, IN 46556  \hfill Morgantown, WV 26506 \\
E-mail: \textit{himonas.1$@$nd.edu}  \hfill E-mail: \textit{fangchi.yan@mail.wvu.edu} 

\end{document}